# LETTERE SCRITTE DA **THOMAS ARCHER HIRST**
# A **LUIGI CREMONA**
# DAL 1865 AL 1892
# CONSERVATE PRESSO
# L'ISTITUTO MAZZINIANO DI GENOVA

A CURA DI GIOVANNA DIMITOLO

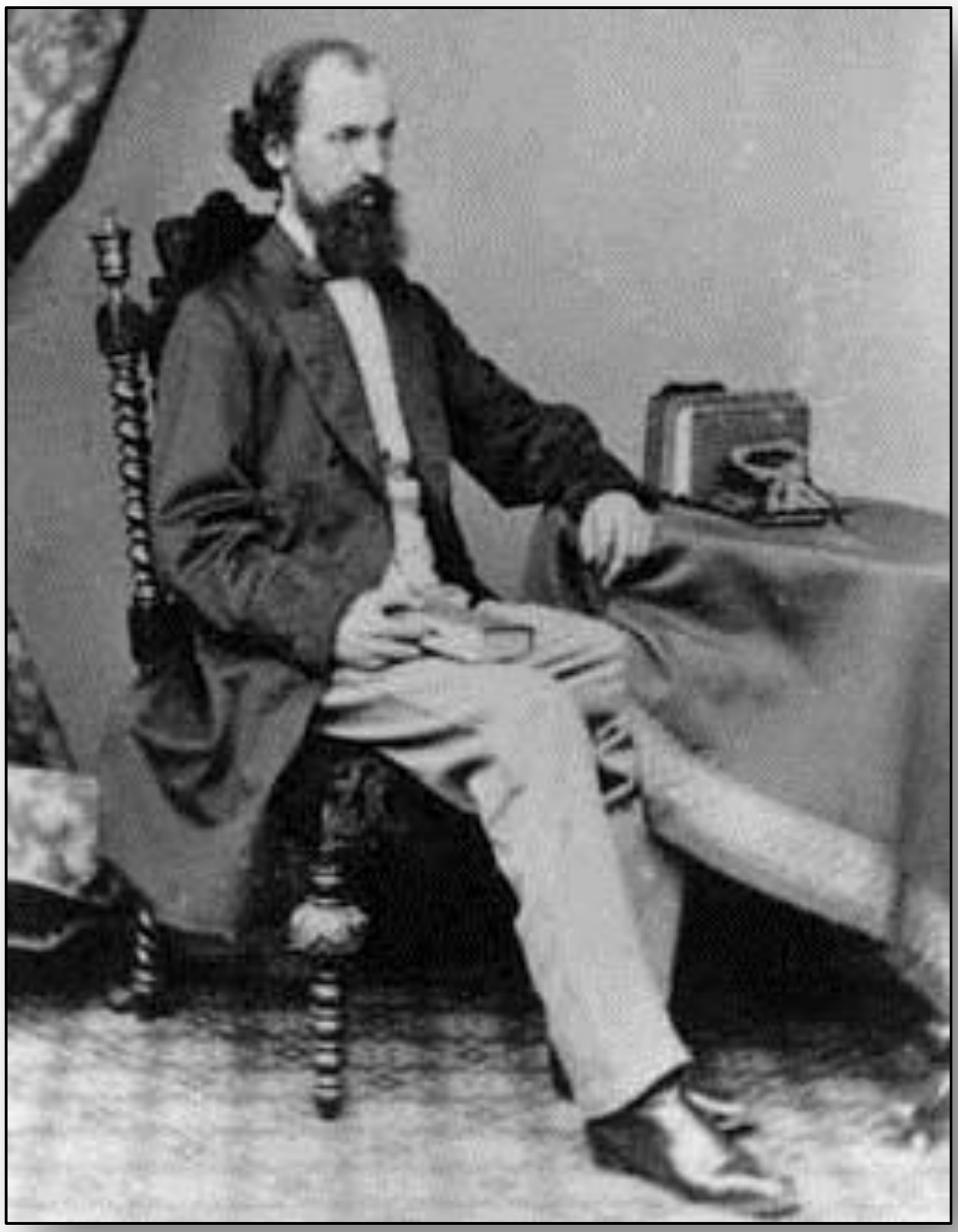

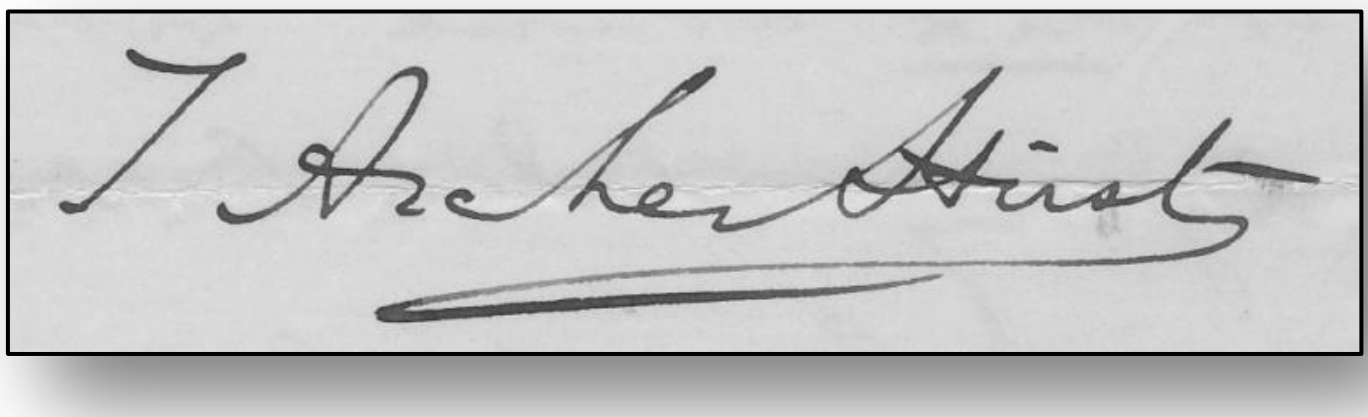

## Indice

## Presentazione della corrispondenza

La corrispondenza qui trascritta è composta da 84 lettere di Thomas Archer Hirst a Luigi Cremona e da 2 di Cremona a Hirst (una delle quali è probabilmente una bozza).
La corrispondenza tra Hirst e Cremona racconta l'evoluzione di un profondo rapporto di amicizia tra i due matematici, della viva passione intellettuale che li accomuna e dell'ammirazione di Hirst per il collega italiano.
Se nelle prime lettere vi sono molti richiami alle questioni matematiche, alle ricerche, alle pubblicazioni sulle riviste di matematiche, insieme a un costante rammarico per il poco tempo che Hirst può dedicare allo studio e un amichevole rimbrotto per gli incarichi politici di Cremona che lo distolgono dalla ricerca, man mano che passano gli anni, la corrispondenza diventa sempre più intima, vi compaiono le preoccupazioni per la salute declinante di Hirst, il dolore per la scomparsa prematura di persone care (il fratello e la prima moglie di Cremona, la nipote di Hirst), le vite degli amici comuni, i racconti dei viaggi e il desiderio sempre vivo di incontrarsi durante le vacanze estive o in occasione di conferenze internazionali.
Fa da sfondo alle vicende professionali e personali dei due matematici la comunità scientifica europea, vivace e in fermento, animata da rapporti umani e vere e proprie amicizie, da curiosità e dispute, da corrispondenze epistolari, scambi di riviste e libri e viaggi di conoscenza.
Hirst, in particolare, è un appassionato viaggiatore. Si reca diverse volte in Italia, in Germania, Svizzera e Francia, e almeno una volta in Grecia, Spagna ed Egitto, sempre alla ricerca di un clima ideale per rinfrancarsi dalla malattia cronica che lo colpisce a più riprese e di fecondi scambi intellettuali con uomini di scienza.

Più di metà di queste lettere ricevono risposta in molte missive di Cremona conservate negli Archivi della London Mathematical Society e trascritte nel IV volume di *Per l'Archivio della Corrispondenza dei Matematici Italiani*[1]. Il dialogo tra i due scienziati si interrompe spesso, infatti molte lettere di entrambi cominciano con espressioni di rammarico per il proprio o l'altrui prolungato silenzio; alcune lettere poi sono andate probabilmente perdute, ma la lettura dei due carteggi permette una ricostruzione quasi completa di una corrispondenza lunga 37 anni. Per facilitare il lettore nella tabella a p. 82 vengono anche riportati la cronologia e i riferimenti delle lettere dell'Archivio della London Mathematical Society, trascritte nel volume già citato, che rispondono o che suscitano una risposta in quelle di Hirst qui trascritte.

## Criteri di edizione

La trascrizione delle lettere ha rispettato il più possibile il testo originale sia nella grafia, sia nella punteggiatura.
Si è mantenuta la sottolineatura "a parola" del manoscritto.
Si sono mantenuti termini che oggi si scrivono diversamente ("to day" anziché "today" e "wearysome" anziché "wearisome"), sono stati segnalati invece con [sic!] termini che appaiono scorretti e, per le parole straniere, in nota è stata indicata la grafia corretta.
I nomi propri sono stati riprodotti come scritti dall'autore, anche se inesatti e sono stati corretti soltanto in nota.
Non è stato possibile decifrare alcune parole e in tali casi si è lasciata tra parentesi quadre un'interpretazione plausibile oppure si è deciso di sostituire le parti non recuperabili con [...].
Tutte le scritte in corsivo tra parentesi quadre sono state inserite dalla curatrice.

[1] L. Nurzia (a cura di), *Per l'archivio della corrispondenza dei Matematici Italiani – La corrispondenza di Luigi Cremona (1830-1903)*, Vol. IV, Quaderni P.RI.ST.EM.

Necrologio di Thomas Archer Hirst
Da "Obituary Notices of Fellows Deceased." In: *Proceedings of the Royal Society of London,* 52 (1892): pp. xii-xviii.

**THOMAS ARCHER HIRST**, the third and youngest son of Mr. Thomas Hirst, a wool-stapler, was born at Heckmondwike, in Yorkshire, on 22nd April, 1830. He was educated at the West Riding Proprietary School; and, in 1844 he became an articled pupil of surveyor at Halifax. It was in this office that he made the acquaintance of John Tyndall, who became a life-long friend of Hirst, and exercised a deep influence on his scientific career; and, in particular, it was the example of Tyndall which led him to give up the pursuit of the profession at first chosen for him. Tyndall had left England to study chemistry, under Bunsen, at Marburg; and thither Hirst followed him, in 1849, to study mathematics, physics, and chemistry. After three years at that University he obtained the degree of Ph.D., on examination in his three subjects and an approved dissertation in analytical geometry. Subsequently, a short time was spent in Göttingen with Gauss and Weber, and then he went to Berlin where he attended lectures by Dirichlet, Steiner, and Joachimsthal. His intercourse with Steiner did much to determine the ultimate bent of his mathematical investigations; but some years elapsed before it was fully indicated, as the majority of his earlier papers are devoted to researches in mathematical physics.

He succeeded Tyndall, at Queenwood, in 1853, on the appointment of the latter at the Royal Institution; but this post was resigned in 1856 on account of the delicate health of his wife, whom he had married in 1854. The succeeding winter was spent at Biarritz and Pau, but without permanent good results, as Mrs. Hirst died in Paris, in 1857. Tyndall took him to Switzerland for six weeks; and, on their return, left him in Paris, where the next winter was occupied in attending lectures by Chasles, Lamé, and others, and in making the acquaintance of the leading French mathematicians. In the winter of 1858-59 he was in Rome; and subsequently he travelled in Italy, making the acquaintance of many mathematicians, especially of his friend Cremona.

Returning to England in 1860, he was appointed mathematical master at University College School; and he held this post for a period of five years, during which he made his first developments in teaching geometry. The experience of the results then attained led him to join the Association for the Improvement of Geometrical Teaching, when it was formed in 1871; and for the first seven years of its existence he was its President and took an active part in its work.

He had been elected a Fellow of the Royal Society in 1861, and it was from this date onwards that his researches are devoted to the various branches of pure geometry which proved of most absorbing interest to him. He was appointed Professor of Physics at University College in 1865, and, on the death of Professor De Morgan in 1867, he succeeded to the Professorship of Pure Mathematics; but the latter chair he resigned in 1870 to become Assistant Registrar in the University of London. In 1873, the date of the establishment of the Royal Naval College at Greenwich, he was appointed Director of Naval Studies; and he continued to discharge the duties of that office for ten years, the precarious condition of his health then compelled him to resign: and he subsequently lived in comparative retirement, spending most of his winters abroad, until his death on 16th February, 1892.

Hirst took a prominent part in the foundation of the London Mathematical Society in 1865, served as its President from 1872 to 1874, and was a member of its Council for over twenty years. His active co-operation with the Society did much to extend its influence: and it was largely to the pages of its Proceedings that his papers on pure geometry were contributed. These are the papers containing the particular researches for which a Royal Medal was awarded to him in 1883 - the year which practically marks the termination of his public life. He had served on the Council of the Royal Society in the years 1864-66, 1871-73, 1880-82: after 1883, the only sign of activity was the production of several papers one of which in continuation of his earlier researches is of considerable importance.

The amount of Hirst's published work is not great; but his work is valuable, and an appreciable portion of it has been translated into French and into Italian. His papers are singularly clear: and side issues, that might lead him away from the main line of development of his subjects, are severely excluded. It is evident that he bestowed great care not merely in carrying out his investigations, but in considering the form in which they are expressed: and his reluctance to premature publication of incomplete work may be gathered from the fact that, though, on quitting the Presidential Chair of the London Mathematical Society on 12th November, 1874, he made a brief statement of some results which he had obtained in the theory of Correlation of Space, the full exposition of his results was not communicated to the Society until 9th January, 1890.

The work with which his name as a mathematician will be most definitely associated is contained in his papers on the Correlation of Planes and the Correlation of Space. The simplest space of two dimensions is a plane, and the elements of such a space are a point and a straight line; a correspondence is established between two planes when all the elements of one plane are connected by a relation or relations with all the elements of the other. When the relations are such that, in general, one element of one plane is associated with one (and with only one) element of the other plane, and *vice versâ*, the correspondence is unique. If each point in one plane corresponds uniquely with a point in the other, and likewise each line in one plane uniquely with a line in the other, the correspondence is called a homography: the theory of homography is considered at length by Chasles in his "Traité de Géométrie Supérieure". If each point in one plane corresponds uniquely with a line in the other, and likewise each line in the one plane with a point in the other, the correspondence is a correlation. A few properties of correlative planes are proved by Chasles in the treatise quoted: it is Hirst's distinction to have constructed the theory of correlation of planes and to have developed it to a great degree of perfection. The extension of the theory of correlation so that it may be applied to space of three dimensions was adverted to by Chasles in his "Aperçu Historique"; the full extension was carried out by Hirst, whose investigations in this subject, together with those of his friends Rudolf Sturm, Cremona, and others, have resulted in important and substantial additions to the theory of pure geometry.

The following memoranda of Dr. Hirst are due to Professor Tyndall; the present state of his health is sufficient to account for their brevity.

A.R.F.[1]

[*Memoranda concerning Dr. Hirst*]

The " railway mania " was at its height, and profitable employment was in prospect for young men trained to the use of the theodolite, spirit-level, chain, and drawing-pen. Youths of well-to-do families were articled in numbers to a profession offering so many attractions. Under such circumstances, Thomas Archer Hirst was articled, in 1846, to Mr. Richard Carter, then resident in Halifax, in whose office, at the time, I happened to be principal assistant. He was then about sixteen years of age. His father had been engaged in the wool trade, his mother was a widow, and he was the youngest of a family of three sons.

[1] Probabilmente si tratta di Andrew Russel Forsyth (1858-1942), matematico, eletto membro della Royal Society of London nel 1886.

The West Riding of Yorkshire was then the battle-field of two great rival railway companies - the West Riding Junction and the West Riding Union. Our duty at the time consisted in the preparation of plans and sections of the lines proposed by the latter company. Save in the office, I saw but little of Hirst at the commencement - he being told off in the field to a party different from mine. But an intimacy gradually grow up between us, and in due time, though I was ten years his senior, we became steadfast friends. Influenced by the writings of Carlyle, Emerson, Fichte, and other philosophers, I held, in those days, very serious views of human conduct and duty. After some time, I noticed that my conversations with Hirst were producing a similar shade of earnestness in his mind. In 1848, I quitted England for Germany, choosing the University of Marburg, in Hesse Cassel, where, in regard to science, Bunsen was then the leading star. Hirst, instead of pursuing the profession chosen for him, soon resolved to follow me to that University. He paid me a preliminary visit in the summer of 1849. It was associated with a pathetic incident. Prior to the examination for the Doctor's Degree, it was customary for the candidates to visit the Professors, and to invite them to be present at the examination. I was on my way to the rooms of one of the Professors, when the postman, meeting me in the street, placed a letter in my hands. It was from a young colleague of Hirst's, who worked in the same office with him in Halifax. I was stunned by the perusal of its first few lines. Hirst had left his mother in good health, and this letter informed me of her sudden death. The writer told me that he had also written to Hirst, but that knowing his strong attachment to his mother, he was afraid to let him know the worst. He trusted to my discretion to disclose it to him in the gentlest manner possible. On returning to our lodgings, I found that Hirst had received the letter announcing his mother's illness, and was making preparations for his immediate return to England. He was far from well. In those days he frequently suffered from a malady of the throat, which entirely quitted him in later years. Everything being prepared for the journey, he had his trunk taken to the coach office, where, after waiting some time, he entered the *Post Wagen*. I was in great perplexity; for, while shrinking from imparting to him the knowledge in my possession, I could not bear to allow him to return, cherishing the delusion that his mother still lived. On squeezing his hand for the last time, I said to him, “Dear Tom, prepare your mind for the worst". His reply was a startled, steadfast look, and, in a moment, he added, “John, my mother is dead. Tell me all; I can bear it". “Yes" I replied, "she is dead”, and he drove away.

After my severance from Halifax, Hirst was accustomed to hold a little weekly symposium in his lodgings. A group of young fellows desiring intellectual intercourse used thus to meet, mainly for the discussion of questions touching upon religion. One of them, I remember, was the son of a Congregationalist minister; another a young author of considerable ability - at that time an ardent admirer of Carlyle, but who afterwards became an equally ardent Roman Catholic. Hirst had now grown into a tall man, with a singularly noble countenance. It was interesting, indeed, to observe how this nobility of expression increased as thought and aspiration mingled more and more with his physical conformation. His hair was dark, his forehead finely formed, his nose and general features well chiselled. The only exception that could be taken to the beauty of his countenance was a certain looseness of lip, which seemed to indicate a lack of firmness of character. But the indication was deceptive, for Hirst could be immovable when circumstances called forth the exercise of firmness.

There was, at that time, near Halifax, a tract of heath-land called Skircoat Moor, at one corner of which stood a little cottage called " The Birdcage". The widow who occupied this cottage eked out a livelihood by selling sweets to the children who came to the moor. She had one son, who had begun life as an errand boy in a printer's office, but who, by good conduct and intelligence, had risen to a highly respectable position in one of the mills. This youth, whose name was Booth, attended Hirst's symposium. His health began to fail, and Hirst observed with anxiety the increasing pallor of his countenance as he walked to and from his work. Medical advice was resorted to, and Booth's malady was pronounced to be consumption. His weakness increased, his usefulness as a clerk diminished, until at last, Hirst insisted that he must cease working and direct the whole of his attention to the care of his health. As to his salary, he (Hirst) undertook to make that good. The poor youth lingered long. Hirst had followed me to Marburg, had quitted that University, and had become a favourite pupil of the illustrious Steiner at Berlin. The *Semester* had begun; its busiest time had set in. One morning, however, he made his appearance in London, and told me that Booth, who was obviously dying, had written, imploring his benefactor to visit him. In response to that letter, Hirst had quitted his studies and had come over to England. Travelling down to Halifax, he found that Booth's chief anxiety related to his mother. “What will become of her," he exclaimed, “when I am gone?" The means at Hirst's command at the time were very moderate. Through the accidents of trade, or the misconduct of individuals, the property left by his father had, in great part, disappeared. Still, without a moment’s hesitation, he gave the dying youth the comforting assurance that his mother would be properly taken care of. For some years afterwards, while Hirst continued in Germany, I was myself the intermediary through whom the widow's allowance passed to her hands.

The sum paid when Hirst was articled to Mr. Carter, together with a professional education of five years, was naturally, on the part of his relatives, expected to produce some tangible result; and when they found that he had resolved to relinquish it all for the sake of the cultivation of his intellect, they thought his resolution a wild one. But he never wavered for a moment, and, except on occasions when his health caused me anxiety, I never wavered in the conviction that the step he had taken was a wise one. Throughout the winter and spring in Marburg, our days were spent in labour, attending lectures, working in the chemical laboratory, and studying at home. At 10 P.M., the stroke of a piano by Hirst gave notice that the work of the day was ended. We had music and light reading afterwards, the latter including the “Essays of Montaigne,” which proved to us a source of strength as well as of delight. At 11 P.M. we went to bed. I was earliest up, for, soon observing that hard study was telling upon the younger and less vigorous constitution, I persuaded him to give more time to repose.

After I had quitted Marburg, leaving Hirst behind me, Dr. Simpson, a medical man with a passion for chemistry, afterwards Professor of Chemistry at Queen's College, Cork, came to the University to pursue his studies in Bunson's laboratory. He brought with him his wife and family and a sister of his wife. They were both sisters of John Martin, the pure-minded Young Irelander. Acquaintance ripened into friendship between the young people, friendship into love, and, after Hirst had taken my place at Queenwood College in 1853, he married Miss Anna Martin. A few years of unalloyed happiness were wound up by an attack of tuberculosis, which ended in her death in 1857. He had taken her to Biarritz, thence to Pau, whence he had returned to Paris, where she died. I was on my way to Switzerland when intelligence of the calamity reached me. Ignorant of their address in Paris, I sought him in his old quarters. Failing to find him there, and guided by a vague indication, I sought him in the Rue Marbeuf. Here, though baulked at first, I discovered where he lodged. On turning, afterwards, the corner of a street, I met him face to face, looking as white as marble. He was returning from ordering his wife's coffin.

I stood beside him in the cemetery of Père Ia Chaise when she was lowered into the grave, and afterwards carried him with me to the Alps. Prior to joining me at Chamounix, he had three days of lonely wandering; communing with himself under his changed conditions. We made the little *auberge* at the Montanvert, which was then a very small affair, our permanent residence. For six weeks his life was filled with healthy exercise, under conditions where constant attention was necessary to his personal safety. Nothing could have been better calculated to divert his mind from the grief which weighed upon it. A few days after our arrival, we were joined by Professor Huxley. As we assembled night after night round our pine-wood fire, life became to all of us more and more a thing to be enjoyed. We returned from the Alps, Hirst halting in Paris, where for some time he took up his abode. Here he made the acquaintance and secured the friendship of the foremost mathematicians. He went afterwards to ltaly; and was at the village of Solferino, helping the wounded, on the day after the battle. In Italy he met Professor Cremona, who remained his friend in a very special sense to the end of his life. On his return to London, Hirst became mathematical master in University College School, then Professor of Applied Mathematics in University College. For the sake of his health, he afterwards accepted the Assistant Registrarship of the University. His final appointment was to the Post of Director of Studies in the Royal Naval College, Greenwich, under the presidency of Admiral Sir Cooper Key, who was succeeded by Admiral Fanshawe; Mr. Goschen was at the time of his appointment First Lord of the Admiralty. Hirst never forgot the Minister's high-minded sympathy with his mathematical studies, or his willingness so to arrange matters as to reconcile the prosecution of those studies with his duties as Director.

During his later years, the state of Hirst's health frequently caused his friends the gravest anxiety. He relinquished in succession the posts he had occupied, retiring finally from Greenwich with a Government pension. During the present calamitous year, he was smitten with influenza, to which he finally succumbed.

I have already given an example of Hirst’s kindness of heart. To this it may be added that, apart from his scientific labours, his life throughout was one of wise beneficence.

J. T.[1]

[1] John Tyndall (1820-1893).

**Traduzione del necrologio di Thomas Archer Hirst** (a cura di G. Dimitolo)

**THOMAS ARCHER HIRST**, il terzo e più giovane figlio di Thomas Hirst, un commerciante di lana, nacque a Heckmondwike, nello Yorkshire, il 22 aprile del 1830. Studiò alla West Riding Proprietary School e, nel 1844, divenne apprendista ingegnere in uno studio ad Halifax. Fu in questo ufficio che fece la conoscenza di John Tyndall, che divenne amico di lunga data di Hirst ed esercitò una profonda influenza sulla sua carriera scientifica; in particolare, fu l'esempio di Tyndall che lo portò ad abbandonare la professione a cui era stato destinato inizialmente. Tyndall aveva lasciato l'Inghilterra per studiare chimica a Marburgo, sotto la guida di Bunsen; e là lo seguì Hirst, nel 1849, per studiare matematica, fisica e chimica. Dopo tre anni di studio presso quella Università egli ottenne il diploma di Ph.D., con esami nelle tre materie e una dissertazione di geometria analitica. In seguito trascorse un breve periodo a Gottinga con Gauss e Weber e poi si recò a Berlino dove frequentò le lezioni di Dirichlet, Steiner, e Joachimsthal. Fu il suo rapporto con Steiner che determinò l'orientamento definitivo delle sue ricerche matematiche; ma trascorsero alcuni anni prima che questo si manifestasse pienamente, infatti la gran parte dei suoi primi studi sono dedicati alla ricerca in fisica matematica.

Nel 1853 Hirst succedette a Tyndall al Queenwood College, quando questi venne nominato alla Royal Institution; lasciò, però, l'impiego nel 1856 a causa delle precarie condizioni di salute di sua moglie che aveva sposato nel 1854. Trascorse l'inverno seguente a Biarritz e Pau, ma non vi furono miglioramenti stabili nella salute della Sig.a Hirst che, infatti, morì a Parigi nel 1857. Tyndall portò Hirst con sé in Svizzera per sei settimane; e, sulla strada del ritorno, lo lasciò a Parigi, dove trascorse l'inverno seguente frequentando le lezioni di Chasles, Lamè e altri e facendo la conoscenza dei più eminenti matematici francesi. Nell'inverno del 1858-59 fu a Roma e successivamente viaggiò per l'Italia, incontrando molti matematici, in particolare il suo amico Cremona.

Al suo ritorno in Inghilterra, nel 1860, fu nominato docente di matematica all'University College School; mantenne l'incarico per cinque anni, durante i quali fece i suoi primi progressi nell'insegnamento della geometria. L'esperienza e i risultati raggiunti lo spinsero a entrare a far parte dell'Associazione per lo Sviluppo dell'Insegnamento della Geometria (Association for the Improvement of Geometrical Teaching), quando venne fondata nel 1871; e per i primi sette anni di attività ne fu il presidente e partecipò attivamente ai suoi lavori.

Venne eletto membro della Royal Society nel 1861, e da questa data in poi le sue ricerche furono dedicate alle diverse branche della geometria pura, materia che lo appassionava più di ogni altra. Nel 1865 venne nominato professore di fisica all'University College e, nel 1867, alla morte del professor De Morgan, gli succedette sulla cattedra di Matematica Pura; ma si dimise nel 1870 per diventare Assistant Registrar[1] all'University of London. Nel 1873, quando venne istituito il Royal Naval College a Greenwich, Hirst ne fu nominato Direttore; continuò ad assolvere i compiti di quel ruolo per dieci anni, poi le sue precarie condizioni di salute lo obbligarono a dimettersi: in seguito visse in relativo isolamento, trascorrendo quasi tutti gli inverni all'estero, fino alla sua morte il 16 febbraio 1892.

Hirst ebbe un ruolo fondamentale nella fondazione della London Mathematical Society nel 1865 di cui fu Presidente dal 1872 al 1874 e membro del Consiglio per oltre vent'anni. La sua attiva cooperazione accrebbe l'autorità della Società: i suoi scritti di geometria pura furono pubblicati principalmente negli Atti della Società. Quegli scritti riguardavano le minuziose ricerche che gli valsero una Royal Medal nel 1883 – anno che segna anche la fine della sua vita pubblica. Fece parte del Consiglio della Royal Society negli anni 1864-66, 1871-73 e 1880-82; dopo il 1883 l'unica testimonianza della sua attività fu la produzione di alcuni scritti dei quali uno, che riprende le sue prime ricerche, è degno di nota. Le pubblicazioni di Hirst non sono molto numerose, ma il suo lavoro è di grande valore e una gran parte è stata tradotta in francese e in italiano. I suoi scritti sono straordinariamente chiari: le questioni secondarie, che potrebbero distogliere dalla principale linea di sviluppo dell'argomento, vengono rigorosamente escluse. È evidente che attribuiva grande importanza non solo allo sviluppo della sua analisi ma anche alla forma in cui questa veniva espressa: la sua riluttanza a pubblicare in anticipo lavori incompleti si può dedurre dal fatto che, quando lasciò la carica di Presidente della London Mathematical Society, il 12 novembre del 1874, espose brevemente alcuni dei risultati che aveva raggiunto nella Teoria della Correlazione dello Spazio, ma, nonostante ciò, l'esposizione completa dei risultati venne comunicata alla Società solo il 9 gennaio del 1890.

[1] Carica amministrativa.

Lo studio al quale verrà certamente associato il suo nome di matematico è contenuto negli scritti sulla Correlazione dei Piani e la Correlazione dello Spazio. Il più semplice spazio nelle due dimensioni è il piano, e gli elementi di tale spazio sono il punto e la linea retta; esiste una corrispondenza tra due piani quando tutti gli elementi di un piano sono collegati da una o più relazioni con tutti gli elementi dell'altro. Quando le relazioni sono tali che, in generale, un elemento di un piano è associato con (e uno solo) elemento dell'altro piano, e viceversa, la corrispondenza è univoca. Se ogni punto di un piano corrisponde in modo univoco a un punto dell'altro piano, e allo stesso modo ogni retta di un piano corrisponde in modo univoco a una retta dell'altro piano, la corrispondenza è detta omografia: la teoria dell'omografia viene esaurientemente trattata da Chasles nel suo "Traité de Géométrie Supérieure". Se ogni punto di un piano corrisponde in modo univoco a una retta dell'altro, e allo stesso modo ogni retta di un piano corrisponde a un punto dell'altro piano, la corrispondenza è una correlazione. Chasles dimostrò poche proprietà dei piani correlati nel trattato già citato: Hirst ha il merito di avere elaborato la teoria della correlazione dei piani e di averla sviluppata fino a un alto grado di perfezione. Fu Chasles nel suo "Aperçu Historique" a estendere la teoria della correlazione allo spazio in tre dimensioni; l'estensione completa fu opera di Hirst, le sue ricerche in questo campo, insieme a quelle dei suoi amici Rudolf Sturm, Cremona e altri, portarono a importanti e notevoli ampliamenti della teoria della geometria pura.

La seguente commemorazione del Dr. Hirst è a cura del Prof. Tyndall; a causa del suo stato attuale di salute sarà breve.

A.R.F.

[*Commemorazione del Dott. Hirst*]

La "ferrovia-mania" era al suo culmine e i giovani addestrati all'uso del teodolite, della livella a bolla, della catena e del tiralinee avevano la prospettiva di un impiego redditizio. I giovani delle famiglie benestanti venivano destinati in gran numero all'apprendistato di una professione che offriva così tante attrattive. In tali circostanze, nel 1846 Thomas Archer Hirst divenne apprendista del sig. Richard Carter, che risiedeva allora ad Halifax, e del quale, all'epoca, si dava il caso che io fossi il vice. Hirst aveva allora quasi sedici anni. Suo padre era stato un commerciante di lana, sua madre era vedova e lui era il più giovane di tre fratelli.

Il West Riding of Yorkshire[1] era allora il campo di battaglia di due grandi compagnie ferroviarie rivali – la West Riding Junction e la West Riding Union. Il nostro compito all'epoca consisteva nel preparare i disegni e le sezioni delle tratte progettate dalla seconda compagnia. Tranne che in ufficio, incontrai molto raramente Hirst all'inizio – essendo egli assegnato al lavoro di campo in un gruppo diverso dal mio. Ma gradualmente crebbe tra di noi una certa intimità e, a tempo debito, si sviluppò una salda amicizia, nonostante io fossi maggiore di dieci anni. In quei giorni, sotto l'influenza degli scritti di Carlyle, Emerson, Fichte e di altri filosofi, avevo un'idea molto rigorosa della condotta e del dovere. Qualche tempo dopo mi resi conto che le nostre conversazioni stavano indirizzando il suo pensiero verso un'analoga serietà. Nel 1848 lasciai l'Inghilterra per la Germania e scelsi l'Università di Marburgo, a Hesse Cassel[1], dove, in campo scientifico, Bunsen era allora la celebrità indiscussa. Hirst, invece di continuare la professione che era stata scelta per lui, presto si decise a seguirmi in quella stessa Università. Mi fece una visita in via esplorativa nell'estate del 1849 e, in quell'occasione, si verificò un commovente episodio. Era consuetudine che i candidati all'esame di dottorato andassero a fare visita ai professori e li invitassero a presenziare all'esame. Mi stavo recando nelle stanze di uno dei professori, quando il postino, incontrandomi per strada, mi piazzò in mano una lettera. Era di un giovane collega di Hirst che lavorava nel suo stesso ufficio ad Halifax. Lessi attentamente le prime poche righe e rimasi sbalordito. Hirst aveva lasciato sua madre in buona salute e quella lettera mi comunicava la sua morte improvvisa. Chi scriveva mi informava che aveva scritto anche a Hirst, ma che, conoscendo il suo forte attaccamento verso la madre, aveva avuto timore di fargli sapere il peggio. Contava sulla mia discrezione per rivelarglielo nella maniera più delicata possibile. Tornando verso i nostri alloggi, scoprii che Hirst aveva ricevuto la lettera che annunciava la malattia di sua madre e che stava facendo i preparativi per il suo immediato rientro in Inghilterra. Stava tutt'altro che bene. In quei giorni soffriva di una malattia della gola che lo abbandonò soltanto negli anni successivi. Tutto era pronto per il viaggio, aveva fatto portare il suo baule nell'ufficio del servizio postale e da lì, dopo una breve attesa, salì sulla carrozza. Mi trovavo in un grande imbarazzo dato che, se rifuggivo dallo svelare l'informazione in mio possesso, non potevo però permettere che tornasse a casa covando nell'animo l'illusione che sua madre fosse ancora in vita. Stringendogli la mano per l'ultima volta gli dissi: "Caro Tom, prepara il tuo animo al peggio" lui reagì con uno sguardo allarmato e fermo e, subito dopo rispose: "John, mia madre è morta. Dimmi la verità, posso sopportarla". "Sì" risposi "è morta" ed egli partì.

[1] Storica area amministrativa dello Yorkshire, suddiviso dal 1889 al 1974 in North, West e East Riding.

Dopo che lasciai Halifax, Hirst aveva preso l'abitudine di tenere un piccolo simposio settimanale nei suoi alloggi. Un gruppo di giovani desiderosi di relazioni intellettuali erano soliti incontrarsi, principalmente per discutere di questioni che riguardavano la religione. Uno di loro, ricordo, era il figlio di un ministro congregazionalista; un altro era un giovane scrittore di notevole talento – a quell'epoca un fervente ammiratore di Carlyle, ma che divenne in seguito un altrettanto fervente cattolico romano. Hirst era allora un uomo alto, con un'espressione straordinariamente nobile. Era oltremodo interessante osservare come questa nobiltà di espressione aumentasse man mano che i pensieri e le aspirazioni modellavano il suo aspetto fisico. Aveva capelli scuri, fronte finemente formata, naso e lineamenti cesellati. L'unica obiezione che si poteva muovere alla bellezza dei suoi lineamenti era un'indubbia mollezza delle labbra, che sembrava l'indizio di una debolezza di carattere. Ma si trattava di un indizio ingannevole perché Hirst sapeva essere irremovibile quando le circostanze richiedevano l'esercizio della fermezza.

C'era, a quei tempi, vicino ad Halifax, un tratto di brughiera chiamato Skircoat Moor, dove, in un angolo, si ergeva una casupola detta "La gabbia degli uccelli". La vedova che occupava la casupola sbarcava il lunario vendendo dolci ai bambini che andavano nella brughiera. Aveva un figlio, che aveva cominciato come fattorino di un tipografo ma che, grazie a buona condotta e a intelligenza, aveva ottenuto una posizione rispettabile in uno dei mulini. Questo giovane, che si chiamava

[1] Così si chiamava fino al 1926 la città di Kassel, in Germania

Booth, frequentava il simposio di Hirst. La sua salute cominciò a venir meno e Hirst osservava con preoccupazione il suo pallore aumentare mentre andava o tornava dal lavoro. Si fece riscorso a un consulto medico che rivelò come la malattia di Booth fosse consunzione. Divenne sempre più debole e meno efficiente al lavoro, alla fine Hirst insistette perché smettesse di lavorare e rivolgesse completamente la sua attenzione alla cura della salute. Per quanto riguardava lo stipendio, egli (Hirst) gli garantì che lo avrebbe rimborsato. Il povero giovane indugiò a lungo. Hirst mi aveva seguito a Marburgo, aveva lasciato quell'università ed era diventato il pupillo dell'illustre Steiner a Berlino. Il *Semestre* era cominciato, era iniziato il periodo più impegnativo. Un giorno, tuttavia, Hirst fece la sua comparsa a Londra, e mi disse che Booth, che evidentemente stava morendo, aveva scritto implorando il suo benefattore di andare a fargli visita. In risposta a quella lettera Hirst aveva lasciato i suoi studi ed era venuto in Inghilterra. Arrivato ad Halifax scoprì che la principale angoscia di Booth riguardava sua madre. "Cosa ne sarà di lei" esclamò "quando me ne sarò andato?" I mezzi a disposizione di Hirst a quei tempi erano piuttosto modesti. A causa di vicende commerciali o di una cattiva gestione, il patrimonio lasciato da suo padre era, in gran parte, svanito. Purtuttavia, senza un momento di esitazione, egli diede al giovane moribondo la confortante assicurazione che si sarebbe preso decorosamente cura di sua madre. Per alcuni anni da allora, mentre Hirst rimaneva in Germania, fui proprio io l'intermediario attraverso il quale l'indennità per la vedova arrivava nelle sue mani.

Naturalmente la famiglia di Hirst si era aspettata che le somme pagate quando era apprendista dal sig. Carter e i cinque anni di educazione professionale avrebbero prodotto dei risultati tangibili; e quando scoprì che aveva deciso di abbandonare tutto per il desiderio di coltivare il proprio intelletto, giudicò la sua decisione avventata. Ma egli non tentennò un solo momento e, con l'eccezione delle occasioni in cui ero in ansia per la sua salute, io non abbandonai mai la convinzione che il cammino che aveva intrapreso fosse quello giusto. A Marburgo, durante l'inverno e la primavera, le nostre giornate erano molto impegnative: frequentavamo le lezioni, lavoravamo nel laboratorio di chimica e studiavamo a casa. Alle dieci di sera il suono del pianoforte di Hirst annunciava che il lavoro per quel giorno era terminato. Poi ci dedicavamo alla musica e a letture leggere, tra queste i "Saggi di Montaigne" che si rivelavano una fonte di ispirazione e di divertimento. Alle undici di sera andavamo a letto. Io mi alzavo per primo, dato che lo avevo persuaso a concedere più tempo al riposo, avendo osservato da subito che lo studio intenso aveva un effetto negativo sulle costituzioni più giovani e meno vigorose.

Dopo che ebbi abbandonato Marburgo, lasciando Hirst, il dott. Simpson, un medico appassionato di chimica, che divenne poi professore di Chimica al Queen's College a Cork, venne all'Università per portare avanti i suoi studi nel laboratorio di Bunsen. Portò con sé sua moglie, la famiglia e una sorella di sua moglie. Erano entrambe sorelle di John Martin, quel nobile animo di Giovane Irlandese[1]. La conoscenza tra i giovani maturò in amicizia, l'amicizia in amore e, dopo che nel 1853 ebbe preso il mio posto al Queenwood College, Hirst sposò la signorina Anna Martin. Pochi anni di pura felicità furono interrotti dall'attacco della tubercolosi, che si concluse con la morte di Anna nel 1857. Hirst l'aveva portata a Biarritz, poi a Pau, da dove era tornato a Parigi, dove ella morì. Io ero in viaggio verso la Svizzera quando mi raggiunse la notizia della disgrazia. Non conoscendo il loro indirizzo, lo cercai nel suo vecchio quartiere. Non trovandolo lì e guidato da un vago indizio, lo cercai in Rue Marbeuf. Qui, nonostante qualche difficoltà iniziale, scoprii dove alloggiava. Più tardi, svoltando l'angolo di una strada, mi imbattei in lui, bianco come il marmo. Stava tornando dall'aver ordinato la bara di sua moglie.

[1] John Martin (1812 –1875) militò nella *Giovane Irlanda*, movimento nazionalistico rivoluzionario irlandese che operò intorno alla metà degli anni Quaranta del XIX secolo.

Stetti al suo fianco nel cimitero di Père la Chaise mentre la bara scendeva nella tomba e poi lo portai con me sulle Alpi. Prima di raggiungermi a Chamounix, trascorse tre giorni in solitari vagabondaggi, meditando sul cambiamento occorso nella sua vita. Facemmo della piccola locanda di Montavert, che era allora cosa assai modesta, la nostra residenza stabile. Per sei settimane le sue giornate furono dedicate all'attività fisica, e, date le circostanze, era necessaria una costante vigilanza sulla sua salute personale. Non avremmo potuto organizzare niente di meglio per distogliere la sua mente dal lutto che la gravava. Pochi giorni dopo il nostro arrivo fummo raggiunti dal professor Huxley. Riunendoci, notte dopo notte, attorno al fuoco di legna di pino, ritrovammo ogni giorno di più, la voglia di godere della vita. Tornammo dalle Alpi, Hirst si fermò a Parigi, dove prese dimora per qualche tempo. Qui fece la conoscenza e si assicurò l'amicizia dei più eminenti matematici. Si recò successivamente in Italia ed era presente nel villaggio di Solferino, il giorno successivo alla battaglia, ad aiutare i feriti. In Italia egli incontrò il professor Cremona al quale lo legò un'amicizia profonda fino alla fine della sua vita. Tornato a Londra Hirst divenne professore di Matematica alla University College School, poi professore di Matematica Applicata all'University College. Successivamente, per riguardo alla propria salute, accettò il ruolo di Assistant Registrarship all'Università. La sua ultima nomina fu a Direttore degli Studi nel Royal Naval College a Greenwich, sotto la presidenza dell'Ammiraglio Sir Cooper Key, a cui succedette l'Ammiraglio Fanshawe; alla nomina di Hirst, Goschen era Primo Lord dell'Ammiragliato. Hirst non dimenticò mai l'intelligente sostegno del ministro ai suoi studi matematici o la sua disponibilità a organizzare il lavoro in modo da permettere a Hirst di conciliare la prosecuzione dei suoi studi con i suoi compiti di Direttore.

Negli ultimi anni lo stato di salute di Hirst preoccupò frequentemente e gravemente i suoi amici. Egli lasciò uno dopo l'altro i posti di lavoro che aveva occupato, ritirandosi infine da Greenwich con una pensione governativa. In questo anno sfortunato fu colpito dall'influenza della quale, alla fine, morì.

Ho già dato un esempio della gentilezza di cuore di Hirst. A questo vorrei aggiungere che, a parte che dal suo lavoro scientifico, tutta la sua vita fu caratterizzata da una assennata beneficenza.

J.T.

## Il carteggio

1

August 24th 1865

My dear Cremona.
A few days ago on arriving at Geneva I received your welcome letter
Your colleague Prof. Capellini has kindly offered to convey to you my reply and I hasten to accept his offer. I hope that before this reaches you you will have completely recovered from the prostration which your arduous duties at the University have caused. I trust you will give yourself a complete rest before you commence working at your own investigations. The meeting here of the Swiss Savants[1] has been a most agreeable one but uninteresting from a mathematical point of view. I have had the pleasure however of meeting and conversing with four mathematicians Christoffel, Prym and Geiser of Zurich and Sidler of Berna.
Prym has been lately appointed Professor at Zurich and I am told he is an able mathematician. He is a pupil of Riemann's and has lately published in the Transactions of the Academy of Vienna[2] an important memoir on Riemann's theory. I hope to read it when I return to England (in a few days) He speaks in great praise of the late work by Durège on the Theory of Complex Functions[3]; this theory it appears from a part of Riemann's theory of which everyone speaks in the highest praise I have a great desire to study it. You will be glad to hear that at last I have been appointed Professor of Mathematical Physics at University College London. The appointment has given me great pleasure I regard it as a proof of the esteem of my London friends. I should have greatly preferred a Professorship of pure mathematics and particularly of Geometry, for I shall be compelled to turn aside from my favourite studying for some time to come. However having at length secured a position in the College I have little doubt that before long I shall be able to improve my position and perhaps make some arrangement with De Morgan the Prof. of Pure Mathematics in order to obtain permission to give a course of lectures on Geometry. The remuneration I shall receive will be very poor at first, the professorship are not endowed and I shall receive but two thirds of the fees which my students pay. At the commencement I do not expect to have more than 40 students but I trust that the number will be increased in a year or two. 
Will you have the kindness to pay the 44 francs to Tortolini I have no doubt that the amount is accurate. I am glad to see that Miller has inserted your solutions (in French) in the Educational Times[4] They will be read with great pleasure
On the 6th of September I hope to attend the meeting of the British Association at Birmingham. If you would like me to communicate anything for you I shall be delighted to do so. (Address as normal 14 Waverly Place)
I have just received a letter from Beltrami to which I enclose a short reply. As it partly concerns yourself perhaps you will have the kindness to read it and then forward it to Beltrami (presso l'Avocato *[sic!]* Nicolò Barozzi, Venezia, S.Maria Formosa, Casa Barozzi)
With sincerest wishes for your health and happiness I remain

my dear Cremona
Ever truly yours
T.A. Hirst

---

[1] 50° Conferenza della Società Svizzera di Scienze Naturali.
[2] F.E.F. Prym "Neue Theorie der ultraelliptischen Functionen. (Mit 3 Tafeln.)", *Denkschr.Akad.Wiss.Wien*, 24_2, 1865, pp. 1–104.
[3] H. Durège *Elemente der Theorie der Functionen einer complexen veränderlichen Grösse. Mit besonderer Berücksichtigung der Schöpfungen Riemanns*, Teubner, Leipzig, 1864.
[4] *The Educational Times* era un giornale rivolto a docenti e studenti e fondato nel 1847 dal College of Preceptors di Londra. Fin dalla sua prima uscita pubblicò problemi matematici proposti dai lettori, con le soluzioni. Questa rubrica ebbe un tale successo che Miller, direttore responsabile della rivista, lanciò una pubblicazione separata: *Mathematical Questions with Their Solutions from The Educational Times*.

2

14 Waverley Place
St. John's Wood
Oct 24th 1865

My dear Cremona.

I am anxious on account of your long silence. It is some months over since I heard from you and I begin to fear that you are not well. The cholera I know visited Bologna and once or twice I have thought with horror that either you or your family might have suffered from it. Do, I pray, write to me at once and put an end to my suspense. I hope, however, to hear that the cause of your long silence is that you have been occupying your vacation with geometrical studies and have been too busy to write to me. I met your colleague Capellini in Geneva and I requested him to take charge of a letter I there wrote to you. It contained nothing of any interest but I hope it reached you. Soon after my return to England I attended the meeting of the British Association at Birmingham and there made one or two communications; one on my own work and another on Chasles' late paper on Conics in Space[1]. I sent you a journal containing a short notice of these communications merely to indicate to you where I was and what I was doing.

I think I told you in one of my former letters that I have been appointed Professor of Mathematical Physics at University-College London. I have now fully entered on my new duties which require me to give six lectures a week, one every morning. I have put aside all other work, and I shall not be able to resume my own studies for at least a year; for my lectures involve considerable preparation. I regret this forced separation from my favourite studies, but I have no alternative; I cannot live on pure geometry.

Cayley is more fortunate; his professional duties occupy him but little and with his usual marvellous industry he is accomplishing much. He is now engaged on the Transformation of Curves after the manner of Clebsch, and is obtaining beautiful results. A few days ago he gave us an admirable paper at the London Mathematical Society. It will interest you greatly and as soon as it is published I will take care that you have a copy. I had an interesting interview with Chasles when I passed through Paris. He was just writing his last paper on Conics in space (which is curiously interesting but perhaps not very useful) and you will have seen that at the end of this paper he examines the question of the independence of the three characteristics of a system of surfaces of the second order subject to eight conditions. He disproves the relation you thought you had established, but curiously enough he does not mention your name, nor does he settle the question whether they are dependent or not.

I tried to obtain information on this point from him but could not do so; his reply was "C'est une question à examiner. Je n'en veux pas parler ". He evidently contemplates writing on the whole subject himself, and it appeared to me that he wished that other people (De Jonquières, yourself and others) would not meddle with his subject. Whenever I attempted to speak on the subject he interrupted me with the words "Ne me parlez pas des surfaces". It is a slight weakness in Chasles, I am afraid, that he wishes to retain for his special cultivation the magnificent field whose gates he has himself opened. If so, however, he will be doomed to disappointment and deservedly so. Prof. Young has been disputing with Sylvester as to the priority of the discovery of the demonstration of Newton's Rule. He has just corrected his own error, however, and in a very noble way confessed that his own demonstration was inaccurate, and that Sylvester's is the only correct one.

Write to me soon, my dear friend, to assure me that you are still well and active

Ever your sincerely
T.A. Hirst

[1] M. Chasles "Systèmes de coniques qui satisfont à sept conditions dans l'espace", *Comptes Rendus hebdomadaires des séances de l'Académie des Sciences*, LXI, 1865, pp. 389-397.

3

14 Waverley Place
St. John's Wood
April 5th 1866

My dear Cremona.
My long silence is not due, I can assure you, to any feeling of anger at your own silence. Its sole cause is the too great preoccupation of mind which is now imposed upon me for my professional duties. I have been waiting too with the hope that I should be able to forward to you the memoirs which you requested me to obtain from Sylvester and from Cayley. The latter gentleman has sent me several memoirs which I forward by this post Sylvester however tells me that he has none at present ready to send to you. Your excellent translation of my paper on Quadric Inversion[1] gave me very great pleasure and the few flattering lines with which you introduce it gratified me still more[2]. Pray send copies to all my friends in Italy not forgetting Bellavitis, Brioschi, Betti, Beltrami, Dini, Battaglini, Trudi, Ianni[3], Padula, Fergola, Gasparis[4], Tardy &c[5] *[nota inserita dall'autore:* Chelini, Nuovi Lincei and the Academies of Bologna and Naples &c*]*
I was so displeased with the second article by Transon in the Nouvelles Annales[6] that, notwithstanding the fact that he speaks respectfully of my memoir, I have already written a small article and sent it to Prouhet for insertion in the Annales. I have of course claimed for Steiner the discovery of projection gauche. The very name, in fact, (schiefe projection) is Steiner's. I told Transon so when I was in Paris last Christmas and notwithstanding this, and the fact that I quote Steiner in my memoir on Quadric Inversion, which Transon has read, he has not had the honesty to correct his oversight himself. You know that I have long had a memoir ready on quadric transformation, had it not been for my professorship it would have been published long ago; in my letter to Prouhet I make several extracts from it and I only regret that Transon has compelled me thus to publish prematurely, for I dislike above all things to publish an incomplete research.
Please thank Beltrami for his kindness in sending me Betti's book. I should be delighted to come to see you in Naples, but I fear I shall not be able; however, if you will not come to see me next year I shall make every effort to come to see you; for I long to resume the pleasant intercourse we had in Bologna nearly two years ago.
Sylvester was greatly pleased with his visit to Bologna and with yourself. I have ordered Neumann's work and hope to be able to understand something of this famous theory of Riemann's. Therefrom. Many thanks for your reply to my question regarding the generalisation of the question in Conics. I am not able to pursue it. I am also (alas) unable to study thoroughly the last splendid memoir by Chasles on Surfaces of the second order and on the extension to curves of all orders of his method of characteristics.
Salmon is filled with wonder as to how Chasles has obtained his results on surfaces. I do hope my friend you will not punish me again by your long silence your letters are a perfect Godsend to me now that I am temporarily shut out from the studies I love most. Do write to me as often as you have leisure and believe me to be Ever yours sincerely

T.A. Hirst

[1] T.A. Hirst, "On the quadric inversion of plane curves", *Proceedings of the Royal Society*, 1865, XIV, pp. 91-106.
L'articolo fu tradotto in italiano da Luigi Cremona:
T.A. Hirst, "Sull'inversione quadrica delle curve piane" (traduzione di L. Cremona), *Annali di Matematica Pura e Applicata*, VII, 1865, pp. 49-65;
T.A. Hirst, "Sull'inversione quadrica delle curve piane" (Memoria estratta dagli *Annali di Matematica Pura e Applicata*), *Giornale di Matematiche*, IV, 1866, pp. 278-293.
[2] Si riferisce alla nota introduttiva di Cremona: "Stimiamo cosa buona e utile il far conoscere ai lettori degli *Annali* questo importante ed elegantissimo lavoro del nostro amico, il Sig. Hirst".
[3] Probabilmente si riferisce a Vincenzo **J**anni.
[4] Si riferisce, probabilmente, ad Annibale De Gasparis.
[5] Eccetera.
[6] M.A. Transon, "De la projection gauche", *Nouvelles annales de mathématiques*, IV, 1866, pp. 63-70.

4

Bristol
Dec. 26th 1866

My dear Cremona.
Your very welcome letter was forwarded to me here yesterday and I lose no time in answering the question which at Prof. Brioschi's request you put to me. I do not know whether the French Government has applied to our Government or not to furnish Reports on the Progress of Science in England during the last 20 years. But I do know that the Royal Society has received no communication from the English Government to that effect and I also know that none of our mathematicians Cayley, Sylvester &c have been requested to prepare such a report. I was speaking a few days ago to Cayley in the very subject. I was telling him that you and Brioschi and Betti where engaged in writing such reports and I was asking him if he had heard anything of a similar project in England. He answered that he had not heard of any such intention on the part of English men of science, and with me he was at a loss to explain why we had not been requested to take a part in the project. I am inclined to believe, therefore, either that the English Government have not been requested to contribute any such reports or, if they have been asked, that they have declined. I have likewise never heard that any other nation except Italy is similarly engaged. The Italian report will unquestionably be of great interest to all European Mathematicians but to have rendered the project a successful one all other nations ought to have cooperated with France and Italy.
I am delighted to hear of your promotion to Milan[1]. It must be a great pleasure to you to have Brioschi as a colleague and to have so good an opportunity of forming an Italian school of pure Geometry. In England pure Geometry is not yet appreciated. I told you in my last letter that I was about to give a course of lectures on the subject at University College. When the time of commencement arrived however (last October) there were 20 few students present that it was not worth any while to continue the lectures and for the present the intention is abandoned. The reason of this is that in the public university examinations a knowledge of pure geometry is not required, and naturally enough students only wish to learn what is absolutely necessary for them. In a few years I hope to see this state of things altered. I am so much engaged with my professorial lectures on Mathematical Physics that I have no time for my own researches. Next year a change will take place. Prof. De Morgan has resigned his Professorship of Mathematics and I intend to become a candidate for the Chair. Should I obtain it I shall be in a better position to promote the study of my favourite science Geometry.
The London Mathematical Society is in a very prosperous state
All our best mathematicians have joined it and very good papers are presented and published
Our proceedings would interest you and if the Istituto Lombardo would desire to exchange their published proceedings the Mathematical Society, although it could not give perhaps a full equivalent, would I am sure be glad to send its proceedings in return for yours. If you approve of my suggestion please let me know and I will propose it to the Mathematical Society.
Dr Salmon has just published a new and enlarged edition of his Higher Algebra[2]; he is no longer a professor of Mathematics at Dublin. To our great regret he has accepted a Professorship in Theology. I hear from Townsend that he will have to give up mathematics to some extent but I hope not altogether. I have read with great pain the controversy between Chasles and De Jonquières[3]. I am afraid there are faults on both sides. Chasles is far too sensitive about his own reputation. His position in Science is so high and incontestable that one cannot but regret that he is not more generous toward younger Geometers. On the other hand de Jonquières ought, I think, to have shown more gratitude to the Tutor to whom he owes so much

[1] Nel 1866 Cremona, su invito di Brioschi, assunse la cattedra di Statica grafica e Geometria superiore presso il Reale Istituto Tecnico Superiore (successivamente Politecnico) di Milano.
[2] G. Salmon, *Lessons introductory to the modern higher algebra*, Hodges Smith and Co., second edition, Dublin, 1866.
[3] Si tratta della polemica intercorsa tra Chasles e de Jonquières sulla paternità della scoperta del "principio di corrispondenza".

He deserves credit for his introduction of the index of a system of curves and had he been silent he would have obtained it. But his conception would have been almost unfruitful had it not been supplemented by Chasles' magnificent researches. Moreover I do not think it is so original as de Jonquières supposes. I remember well that Steiner in his lectures at Berlin made use of the conception to some small extent. Remember me most kindly to Brioschi, Beltrami, Chelini, Bellavitis to M[rs] Cremona and to your children.
Believe me to be with all good and reasonable wishes
Ever your sincerely

T Archer Hirst

5

Athenaum Club
London
Feb. 17[th] 1867

My dear Cremona.
on the part of the London Mathematical Society I have to acknowledge, with many thanks, the receipt of your "Introduzione" and "Preliminari" &c. Both works will, I assure you, meet with attentive readers.
I have also to inform you on the part of the Society that I forward herewith a complete set of our Proceedings (as far as at present published) and I am desired to request you to present it to the Istituto Lombardo in the name of the President and Council of the London Mathematical Society. Although at present we can scarcely expect to be able to offer a full equivalent we should be glad in future to receive, in exchange, the published proceedings of the Istituto Lombardo
Believe me to be

Ever yours sincerely
T Archer Hirst
Treasurer

All memoirs intended for the London Mathematical Society should be addressed thus:

To the Secretary of the
London Mathematical Society
Burlington House
Piccadilly
London

6

Paris July 3rd 1868

My dear Cremona
I am on my way to Switzerland with my friend Prof. Tyndall and before I return to England again I should very much like to see you, if that can be done without travelling too far. From your last letter I learned that you were on a tour of inspection of Schools and did not expect to return to Milan for some time. On receiving this will you write a line to me at the Post Office, Chatillon (Val d'Aosta)
I expect to pass through that place in about 10 days or a fortnight and if you should then be at Milan or Genoa I would gladly come to see you for a day or two. I dined with M Chasles yesterday and told him that I should probably see you soon. He desired me to be the bearer of his best wishes and compliments to you
I will not write at greater length to day but will wait to see whether I cannot procure for myself the greater pleasure of an interview with you
With kindest regards
Believe me to be
Very truly yours

T Archer Hirst

7

London. 10th April 1869

My dear Cremona,
Prof Cayley appears to be now occupying himself with the subject of Geometrical Transformation, and he has also been the means of inducing Clifford to work at the same subject. Both have been reading your papers on the subject published in 1863 and 1865 in the Memorie dell'Academia *[sic!]* di Bologna[1] with which, and especially with the latter, they are naturally greatly delighted. Clifford made a very brief communication on the subject at the last meeting of the Mathematical Society and I understood him to say that your general transformation of a right line into a curve of the nth degree satisfying the conditions (1) and (2) p.5 of your memoir (1865) was equivalent to a repitition[2] *[sic!]* of quadric transformation with properly chosen principal points. I made the remark that in my communication to the British Association at Bath in 1864 *[nota inserita dall'autore:* Alas still unpublished*]* I had certainly shown that all the results obtained by the repitition *[sic!]* of quadric transformation were included in the general theory which you had given and that I remembered your making a similar remark to me in conversation
Indeed I am not sure whether you did not also say that the converse of this (as Clifford now asserts) was also true. At the next meeting of the Mathematical Society on May 13th I shall if possible communicate more fully the results I obtained 5 years ago when studying the properties of quadric transformation of the most general kind. I see no prospect of being able to complete my studies and therefore I should like the results at which I then arrived to be more generally known, even though they are more incomplete than I could wish. I should be very glad if you would write to me before the 13th of May to tell me whether or not my memory of our conversation on the subject has failed me. If you have the time and think fit to supplement your beautiful memoir of 1865 with any results at which you may possibly have arrived since its publication it would be a most appropriate time to do so, and I should be proud to communicate them to the Society. It is a great pleasure to me to see English Geometers turning their attention to this subject

[1] L. Cremona, "Sulle trasformazioni geometriche delle figure piane", *Mem. dell'Acc. delle Scienze di Bologna*, 2, 1863, pp. 621-631, *Giorn. di mat.*, 1, 1863, pp. 305-311.
L. Cremona, "Sulle trasformazioni geometriche delle figure piane. Nota II", *Rend. dell'Acc. delle Scienze di Bologna*, 1864-65, pp. 18-21, *Mem. della R. Acc. Nazionale dei Lincei,* 5, pp. 3-35, *Giorn. di mat.*, 3, 1865, pp. 269-280, 363-376.
[2] repetition.

and a still greater pleasure to hear the universal high commendation bestowed upon your researches by every one who studies them. I can remember the day when I was the only Englishman acquainted with your geometrical writings; where as[1] *[sic!]* now there is not one of any mathematical reputation who has not more or less studied, and profited by them. This is but the partial fulfilment of that prediction of mine which I expressed to you long ago in our walks about Bologna; and I need scarcely assure you that in my present powerlessness to continue my own studies I find most consolation in watching the steady fulfilment of my prediction.
Do not forget that we all hope to see you in August next at the meeting of the British Association
With kindest regards to Prof. Brioschi and to Mrs Cremona

Believe me to be
Yours very sincerely
T. Archer Hirst

P.S. I received your letter to the Editor of the Giornale di Matematica[2]
There is much truth in what you say of Wilson's book and there is very little difference in our views of what the Elementary Geometry we now require for schools should be.

8

Exeter. Aug 26, 1869

My dear Cremona.
The bearer of this letter certainly needs no introduction from me to ensure a hearty welcome from you, for he is the Spottiswoode whose name is already well known to you as a mathematician and whom you also know to be one of my most esteemed friends.
He has kindly offered to place in your hands the Memoire Courronnè [3]*[sic!]* which I have at length obtained from Prof. Smith[4]. The author desires me to state that there are only six copies of this memoir in existence. Three sent to Berlin, one in my hands, the one sent to you and one retained by prof. Smith himself. No further copies will be issued until you have decided upon its publication in the Annali and actually published the number containing it.
Mr Spottiswoode will also give you a list of Errata in Prof. Smith's last paper in the Annali.
The meeting of the British Association here has just terminated but I have not yet leisure enough to write you a longer letter. Wish kind regards to Brioschi & Casorati whose acquaintance I hope Mr Spottiswoode will make I remain

Yours very sincerely
T Archer Hirst

[1] whereas
[2] F. Brioschi, L. Cremona, "Al signor Direttore del Giornale di Matematiche ad uso degli studenti delle Università Italiane – Napoli", *Giornale di Matematiche ad uso degli studenti delle Università Italiane*, Vol. VII, 1869, pp. 51-54.
In risposta al discorso di Wilson pubblicano in forma di articolo sullo stesso giornale l'anno precedente:
J.M. Wilson, "Euclide come testo di Geometria Elementare", *Giornale di Matematiche ad uso degli studenti delle Università Italiane*, Vol. VI, 1868, pp. 361-368.
[3] Couronné.
[4] H.J.S. Smith, "Mémoire sur quelques problèmes cubiques et biquadratiques", *Annali di Matematica Pura ed Applicata*, Ser . II, vol. III, pp. 112-165 , 218-242.

9

*[su carta intestata:*
Athenaeum Club*]*
Athenaeum Club.
Oct.$^{r}$ 4$^{th}$ 1870

My dear Cremona.
M. Tchebichef[1], who has been our visitor for the last two weeks, has just told me that in half an hour he will leave to prepare for his journey to Italy. Although the notice is short I profit by the opportunity of conveying to you my best wishes, and of introducing you to an eminent Russian Mathematician whom I know you will be ready to welcome.
Since I last wrote I have been very ill but I am now well again and I hope to have more leisure this year for doing a little Geometrical work on my own account.
I must not forget to tell you that I have never yet received news of the safe arrival of that Microspectroscope which at Brioschi's request I gave order to be sent to Milan. I hope it did arrive without injury. I shall be delighted to have news from you, I have already received from you one or two memoirs which prove that you are still industriously increasing our knowledge.
With kind regards to M$^{rs}$ Cremona, to Brioschi Casorati and other friends Believe me to be
Yours very sincerely
T Archer Hirst

[1] Pafnutij L'vovič Čebyšëv.

10

University of London
Burlington Gardens
London
Dec[r] 15th 1871

*[di fianco all'indirizzo l'autore ha scritto:* Please use this address in future.*]*

My dear Cremona.
The last letter of yours which I can find is dated May 6th and the first words in it are "I blush to have postponed writing to you for so long a time". I ought therefore for still stronger reasons to blush for having neglected to reply for so long time and the more because you expressly wished me not to take any revenge in this manner
I can only say that I have often wished to write and that whenever I have done so some trivial occupation or other have prevented me.
Besides since your letter arrived I have been at Edinburgh at the British Association and afterwards in the north of Scotland on an excursion, and after that I was very *[busy]* with official matters here.
It is some consolation to me however to find that I am at length not only able to write but to communicate what I trust will be acceptable news
You were yesterday unanimously elected a Foreign Member of the London Mathematical Society. You will shortly receive your Diploma, and henceforth you will be entitled as Foreign Member to a copy of our Proceedings.
Hitherto Chasles was the sole Foreign Member of our Society. The Council decided however to elect five more and, without any intention of doing so, the six foreign members now elected are equally distributed between the three principal continental nations. They are as follows:

Italy
Cremona and Betti
Germany
Hesse and Clebsch
France
Chasles and Hermite

I know not whether to congratulate you or the Society most. It is to me a very great satisfaction to see my predictions of by gone years becoming steadily fulfilled and on the other hand it is a pleasure to know that England is not behind other nations in recognising your well merited claims to stand in the front rank of living Geometers.
Notwithstanding the interruptions of work and health I still continue to do something in Geometry
I am in the middle of a large and (I trust you will think) a not unimportant investigation
I have determined to say no more about it however until it is ripe for publication
When that will be I cannot tell, but I will take care that it shall be placed only at your disposition whenever it shall be ready.
The new number of the Annali has reached me today (the last number of the volume) and the expenses of transit have amounted to one shilling. I mention this merely because it has never happened before and I am at a loss to know what has caused the change.
I hope to pay a visit to Paris next February
I wish there was a chance of meeting you there.
Wish kind regards to M[rs] Cremona and all other friends at Milan

Ever your sincerely
T. Archer Hirst

11

University of London
Burlington Gardens
Dec$^{r}$ 26$^{th}$ 1871

My dear Cremona

I received your letter yesterday, and to day I wrote to Cayley, to Maxwell and to the Secretary of the Mathematical Society conveying your requests to them. I have no doubt you will shortly receive the Memoirs you wish to possess as well as N° 34 of the Proceedings of the Mathematical Society.

I reply to your letter without delay because I am interested in the news you give me that the Italian Government has offered a prize of 2500 lire for the best text book on Geometry after the manner of Euclid. On the 12$^{th}$ of January I shall have to preside at a meeting of the Association for the improvement of Geometrical Teaching, and I should very much like to communicate to them authentic information on this interesting and important step on the part of the Italian Government. I should feel greatly obliged to you therefore if you would send me every thing that has been published on the subject, as well as a copy of those programmes of 1867 to which you allude, as speedily as you can

Please give me also your opinion (in confidence) of the Geometry of Sannia and D'Ovidio[1] and tell me what kind of reception it has met with in Italy.

I am delighted to hear that you are about to publish a book on Projective Geometry; I shall look forward to its arrival with great interest.

Accept my best wishes, my dear friend, for your happiness and for that of M$^{rs}$ Cremona and your children, during the new year which is now close at hand, and believe me to be ever yours sincerely

T. Archer Hirst

P.S. I conclude from the fact that two papers by Weyr are published in the Annali t IV, Fas. 4. and that both are dated Milano Marzo 1871 that their author must have been paying you a visit.

I do not know him personally but I read all his papers with great interest and pleasure; so much so that I have taken some pains to collect all his memoirs. Is he young or old? Has he gone back again to Prague or is he still in Italy?

[1] E. D'Ovidio, A. Sannia, *Elementi di Geometria*, Stabilimento Tipografico delle Belle Arti, Napoli, 1869.

12

University of London
Jan 31th 1873

*[in testa alla lettera l'autore ha aggiunto queste righe:* what a loss mathematical science has sustained by the death of Clebsch! The regret here and in Scotland and Ireland was profound.*]*

My dear Cremona
I have never yet written to thank you for your last two esteemed presents. I mean your Elementi di Geometria Projettiva[1] and your elegant little work Per le Nozze di Camilla Brioschi[2], both of which duly reached me and have been read with the greatest interest. With the former, long expected work, I have been especially pleased; so much so that I have often wished that a translation of it into English could be published and introduced into our Colleges and Universities. I shall keep this desirable object in view and it is not improbable that I may, before long, have a better opportunity of carrying it out than I have at present; for I have just accepted from our Government the important Post of Director of Studies in the new Naval College which is about to be instituted at Greenwich near London. Hitherto the Scientific Training of the officers of our Navy has been but a very imperfect one as-length, however, the necessity of improving it has been acknowledged and our Admiralty have decided to found a Naval College and to provide it with a thoroughly competent staff of Professors of the various Theoretical and Technical subjects which are more immediately connected with the Naval Department.
At the head of the College there will be an Admiral (Admiral Key) and next to him a Director whose duty it will be to organize and superintend the whole of the courses of study to be conducted by the various Professors. This is the position which has been offered to, and accepted by me. I shall resign my present functions at the University of London at the end of February, and enter upon my new duties on the 1st of March. (my address will be Athenaeum Club as before) One of my first duties will be to pay a visit to the Naval Establishments at Home (our Dockyards &c) and afterwards to inspect the Naval Colleges on the Continent. I shall certainly visit the German one at Kiel[3] as well as the French, Austrian and possibly Italian Institutions for training officers of the Navy. It is not improbable, therefore, that before very long I may have the great pleasure of seeing you again. When you next write I should be very glad if you could give me some information about the Naval Institutions of Italy or at all events indicate in what manner I might best procure such information. Judging from the excellence of the Istituto Tecnico di Milano I am led to hope that in Naval Matters also I may be able to learn much from Italian Institutions. You will not be surprised to hear that these and the manifold other duties in which I am always involved interfere very seriously with my own studies I have amongst my papers a great quantity of Geometrical Material, but the leisure to coordinate it and put it in a form for publication has not been granted me. For some time to come I shall be even more occupied than I have been; but once the New Naval College has been organized, I may look to have more leisure. Ever sincerely yours

T Archer Hirst.

[1] L. Cremona, *Elementi di geometria projettiva*, Paravia, Torino, 1873.
[2] L. Cremona, *Le figure reciproche nella statica grafica*, Hoepli, Milano, 1872. (La prima edizione venne pubblicata il 1 giugno 1872 in occasione delle nozze di Camilla Brioschi, figlia di Francesco, con Costanzo Carcano).
[3] Kiel – capoluogo dello Schleswig-Holstein, stato federato tedesco.

13

Bel Alp[1]
Pres Brigue
Canton de Valais
La Suisse
August 11th 1873

My dear Cremona.
I am now separated from you merely by a ridge of mountains, which ridge I would gladly pass over if I felt sure I should be able to find you on the other side. It is my intention to remain here at least a week longer and probably ten days. I have to return to Greenwich on the 1st of September, but whether I shall return by Geneva or by Milan and Turin will depend solely upon your reply consequently remained unanswered
If I come now it will be a short visit, made purely for the sake of grasping your hand once more, and renewing our friendly intercourse
Pray remember me very kindly to Madame Cremona, and believe me to be ever yours sincerely
T. Archer Hirst

14

Arona
Augst 23rd 1873

My dear Cremona
I propose paying you a short visit at Rapallo on Monday next. I shall leave Genoa by the 10.30 A.M. train and arrive at Rapallo at 12.10.
As I shall have to leave Genoa (for Paris) on Wednesday morning at the latest I shall have at most two days during which I can enjoy your company.
Rapallo is so near Genoa that it will scarcely be worth the trouble of changing my hotel (from Genoa) for the sake of passing one night at Rapallo.
Hoping to see you shortly
I remain
Ever yours sincerely
T. Archer Hirst

[1] Belalp nel distretto di Briga in Svizzera.

15

Royal Naval College
Greenwich
May 17th 1874

My dear Cremona
Nine months have passed away since we saw each other at Rapallo and not a single letter has been interchanged. It is true you sent me your new address at Rome and I concluded from the fact of your leaving Milan, that an improvement in your position[1], of which we were speaking when I last saw you, has at length taken place
I shall be very glad to hear from you that my long cherished hopes respecting a proper recognition of your position as a Geometer and man of Science have been fulfilled.
I have been terribly busy since I saw you.
Last Easter, however, I had a week rest and employed it in preparing, finally, for the press the paper on the correlation of two planes[2] of which I spoke to you
I send a proof by this post for your acceptance
I have had the paper printed at my own expense in order to save time. On Thursday last I communicated it to the London Mathematical Society, in whose Proceedings it will probably be published eventually. The publication, however, cannot take place for some months, so that if you think the paper an appropriate one for the "Annali" it is quite at your disposal, and I shall deem it an honour to be ranked at length as one of your contributor
I trust you will be able to carry out your intention of paying us a visit this year
The Secretary of the British Association also expressed to me the hope, a few days ago, that you would be one of the visitors at the Belfast Meeting (in Ireland). I assured him that I would use all my influence to prevail upon you to come. The meeting will be hold about the 17th of August. If this should be too early for you I still hope you will come to England in September. Whenever you do come, remember that my house is to be your home.
With kindest remembrances to Mrs Cremona, your daughters and your son I remain
Ever yours sincerely
T. Archer Hirst

[1] Nel 1873 il ministro della Pubblica Istruzione chiamò Cremona a sovrintendere la Regia Scuola d'applicazione per gli Ingegneri a Roma.
[2] T.A. Hirst, "On the correlation of two planes", *Annali di Matematica Pura e Applicata*, 1873-1875, s. II, VI, pp. 260-297;
T.A. Hirst, "On the correlation of two planes", *Proceedings of the London Mathematical Society*, 1873-1874, V, pp. 40-70.

## 16

Wiesbaden
August 2nd 1874

My dear Cremona.

I received your letter shortly before I left England and was sorry to learn from it that you have postponed your long promised visit to England. I can quite understand your reasons for doing so, however, for I have been in similar circumstances myself. Next year you will I trust have no such excuse and at length we shall see you. My health has been so bad during the last year that I thought is my duty to follow the advice of medical friends and pay a visit to Germany instead of to the British Association at Belfast. Hence it is that I answer your letter from this place where I propose to remain for a week or two to drink the waters, bathe, and pass idle days. Immediately after the receipt of your letter I sent you, by post, another copy of my paper on correlation which contains a note or two which did not appear in the proof I first sent to you. The memoir will not be published in England until September or October, perhaps, so that there is every probability that it will appear first in your journal as I intended it should do. It is true I have sent copies, privately, to Chasles, Schröter[1], Reye, Zeuthen and a few other friends but this will not interfere with its appropriate publication in the Annali. I need not add that I am proud to be as length a contributor, even to a small extent, to your excellent journal

I shall feel obliged if <u>at my expense</u> you will cause about 100 copies to be printed, separately, and sent to me.

I sincerely wish you health and strength to enjoy the more worthy position of Director of the Engineering School at Rome which has been justly, though tardily, assigned to you.

I have long looked forward to the time when my duties will allow me to pay a second visit to Rome. I shall do so now with increased pleasure since you are there.

Give my kindest regards to Mrs Cremona and to your children and believe me to be as ever yours very sincerely

T. Archer Hirst

[1] Schröeter.

17

Royal Naval College
Greenwich
Nov^r 15^th 1874

My dear Cremona

Immediately on receiving your letter of the 31^st of October I wrote to the Secretary of the Mathematical Society to request him to forward to you the missing numbers of the Proceedings, and also to Messr^s Hodgson[1] about the volume of the Reprint which you failed to receive at the proper time. I trust that by this time these have reached you safely and that in future there will be more regularity.

I have received the Annali quite regularly. The Giornale di Mathematiche *[sic!]*, however, often fails to reach me safely. It is usually sent in a very insecure wrapper *[l'autore ha inserito una nota*: The Annali also sometimes reaches me in bad condition.*]* Not unfrequently the wrapper alone reaches me, the journal itself being lost. I wrote a short time ago, when this had happened, to the English Postmaster General to enquire if the missing numbers were in his possession, and sent him the empty wrappers (or covers) to assist him in their identification. Fortunately he was able to find them, and my set is now complete with one exception. The missing number is that for July-August 1871. Being more than three years ago I don't precisely remember whether the empty cover ever arrived. I should be very glad if you would ask the Publisher to send me the numbers more securely wrapped up in future, and if possible also the missing number above mentioned. I have received and returned the proof of part of my paper on the Correlation of two planes. I am now extending my method to correlation in space and have already opened up a wide field for investigation. Its connection with your researches and Cayley's on rational transformation between two spaces interest me very much, but of this at another time. Trusting you are well and successful with your new school

I remain ever yours T. Archer Hirst

Query

Has not Chasles established the relations which exist between the characteristics of a system of quadric surfaces satisfying 8 conditions and the numbers of degenerate quadric surfaces which the system contains?

I cannot put my hand on it but I distinctly remember its being done either by Chasles or yourself

Where was it that Chasles stated that a pair of planes requires to be supplemented by a pair of points in order to become a degenerated quadric?

[1] Editore dei *Proceedings of the Royal Mathematical Society.*

18

Royal Naval College
Greenwich
July 14th 1875

My dear Cremona
It was indeed a great pleasure to me to learn from your letter that there is a probability of our spending a few days together next month, though it is not without disappointment that I learn that those few days will not be spent in England. You will no doubt receive an invitation from the officers of the British Association
The meeting this year is to be held at Bristol and to commence, I believe on the 25th of August
If so it will not be possible for you to attend it as well as the meeting of the Association Française at Nantes.
August being the only month during which I have holiday it was my intention to spend it in Germany and not to go to Bristol. I have received an invitation to go to Nantes, and was about to decline it also when your letter arrived
The opportunity of seeing you there, however, was irresistible and I have accordingly accepted the invitation of the Association Française. Unless something unforeseen prevent me, therefore, from carrying out my intentions we shall meet at Nantes.
I have received to day an invitation from your Società Italiana to go to Palermo[1]; but I have been obliged to decline it, for my time will not permit me to be absent so long from England. I propose to leave home on the 1st of August and to go to Darmstadt where I have promised to spend a short time with Sturm. From Darmstadt I will pass into France, pass through Paris to see Chasles, and then meet you at Nantes on the 19th of August. You do not state in your letter whether you propose to go to Paris. I should be very glad to know precisely what your plans are, for I need not say that I would endeavour so to alter mine, if need be, that we may spend the longest time possible together.
I regret to say that since last January I have not been able to do anything at "Correlation in space" and I regret still more to have to add that I fear I shall be unable to continue my studies for some time to come so overburdened am I with official duties.
Give my kindest regards to Mrs Cremona and all your family. I hope your duties as Director of the School for Engineers are becoming lighter and that you like your appointment
Believe me to be

My dear Cremona
Ever yours sincerely
T. Archer Hirst

[1] Si tratta del XII congresso della Società Italiana per il Progresso delle Scienze, che si tenne a Palermo dal 29 agosto al 7 settembre 1875.

19

Royal Naval College
Greenwich
July 26th 1875

My dear Cremona
I lose no time in replying to your letter
I shall be in Paris on the 10th of August at the latest
This will give me a clear week there, for on the 17th about (or 18th al latest) we ought to proceed to Nantes
I generally go to the Hotel de Lille et d'Albion, Rue St Honorè at Paris. At all events I will take the precaution to leave my address there in order that you may speedily find me.
I cannot tell you how pleased I am at the prospect of seeing you again soon.

Ever your sincerely
T Archer Hirst

20

Greenwich, May 31st, 1876

My dear Cremona.
Prof. Blaserna, to my great regret, was unable to deliver your letter of the 6th of May personally; but he sent it to me by post together with one of his own in which he gave reasons for his inability to pay me a visit. My occupations at Greenwich were so numerous that I was not able to be present at any of the conferences held at South Kensington in connection with the loan Exhibition of Scientific Instruments[1], and thus I lost the opportunity of making Prof. Blaserna's acquaintance
Pray express to him my great regrets thank him cordially for the very friendly letter he wrote to me, and assure him that I share his hope that on another occasion, not far remote I trust, we may be more fortunate.
It was precisely at this very busy period to which I have alluded that our Paris acquaintance Signor Gullo paid me a visit. I was merely able to introduce him to one of his countrymen Martorelli, who is studying Naval Architecture at Greenwich, and to give him an introduction to my friend Siemens, well known to all Engineers and Scientific Men. I can only trust that both these introductions were useful, as well as agreeable to him. By this time you will have received the Richards Indicator, as well as the requisite information from Elliot-Brothers, the makers, as to the most convenient way of transmitting money to them. With respect to your other commission concerning the Report of the Royal Commission in Science and Scientific Instruction[2], I have not been able yet to make the necessary enquiries
I only know that these Reports are <u>very voluminous</u> as "Blue Books" generally are. They contain a large quantity of matter in which neither you nor I will be interested, and therefore before sending you any I should wish to examine them. I have ordered them for the Library of this College, and should I find anything of interest enough to warrant my doing so, I will send you what I find. You will however have an opportunity of looking at them here for I am delighted to know that you have returned to your original project of paying us your long promised visit
The Glasgow Meeting will commence in September 6th, much earlier therefore than I had anticipated. I shall not be able to leave Greenwich for my Holidays until early in august and then I propose to pay a visit to Germany and Switzerland, to which places I am under a promise to take Miss Hirst, my niece. We shall return from our Tour early in September. We shall certainly be at home by the time the meeting at Glasgow is over but I have my doubts about being able to attend the meeting itself. In any case we shall look forward to receiving you at Greenwich on your return from Glasgow if not before. But could you not pay a visit to Germany while I am there? Think of this and let me have your opinion.
On Geometrical matters I have no news to send you. The little spare time I have had during the Session I have devoted to detailed work in connection with correlation in space; on which subject I have now a good deal of matter accumulated, but it wants arrangement and the time for methodical arrangement is just what I cannot secure, consistently with my present duties.
I hope you are not equally unfortunate; but I fear you are
Wish kind regards to Mrs Cremona
Believe me to be

Ever yours sincerely
T Archer Hirst

You and Prof. Betti will receive the usual invitation to the Glasgow Meeting from the Office of the Association.

[1] "Special Loan Exhibition of Scientific Apparatus" – esposizione internazionale inaugurata il 13 maggio 1876, aperta al pubblico dal 15 maggio al 30 dicembre 1876 presso il South Kensington Museum.
[2] Royal Commission on Scientific Instruction and the Advancement of Science – Commissione Reale del Regno Unito, attiva dal 1870 al 1875.

## 21

Athens
Jan$^{y}$ 31$^{st}$ 1877

My dear Cremona
Your letter of Dec$^{r}$ 26 has been forwarded to me here by my niece who has informed me that the photographs have arrived safely. They will be sent by her to the person for whom they are intended. You will probably have heard from Prof. Jung that in consequence of a serious illness I had to leave Greenwich last november[1] *[sic!]* and that I have been in Egypt. I came here 10 days ago and on Sunday next I leave again on my way to Rome where I hope to arrive in about a fortnight. I will tell you then all that has happened to me since you were at Greenwich
My health is much better now than it was, but I am still far from well. I know you will desire to show me great hospitality and I shall be glad to have your company during the time I remain in Rome. I shall of course be frequently at your house but I must at once tell you that in consequence of the state of my health and of the many arrangements I have to make to suit it I shall not be able to accept your hospitality as few as to take up my abode under your roof. You must allow me therefore to <u>live</u> at my hotel but to regard your house as the one where I shall always find a welcome during the days I spend in Rome. I cannot now say on what day I shall reach Rome neither have I yet decided to what Hotel I shall go
But I shall seek you without fail immediately after my arrival
Until then goodbye Ever yours

T. Archer Hirst

P.S. Give my very kind regards to Mrs Cremona and all necessary explanations
I send my greeting also to your daughters and son

## 22

Albergo del Quirinale
April 15$^{th}$ 1877

My dear Cremona
It will not be possible for me to come and dine with you to day. I have only just now left my bed and my malady is still so acute that I must remain still all day.
Tell Schroeter that if I am well enough I will go on to Florence on Tuesday
I shall descend at the Hotel de la Paix.
And if he will kindly call and leave his address there I will call upon him without fail.
Ever yours sincerely

T. Archer Hirst

[1] November.

23

Hotel de Lille et d'Albion
Paris May 2nd '77

My dear Cremona.
I left Turin at 9 P.M. on Monday and arrived here at 6 P.M. yesterday. Before I left Italy I received your affectionate letter, and the one from Sturm which you kindly forwarded
Tomorrow I propose to continue my journey homewards, but before I leave Paris I must reply, by a few lines, to your letter; even if it be but to assure you that its perusal gave me intense pleasure and satisfaction
Your friendly words found a perfect echo from me. From the time I first made your acquaintance in Cremona, I have felt that between us there was more than sympathy of scientific tastes that the bond between us had its points of attachment not only in our common study, but in the deeper sympathies, moral and intellectual, which we had in common. Although I have been long convinced of this, the expression on your part of the friendly feelings you entertain for me were none the less welcome – these expressions as I have already said found in me perfect reciprocity. Your brotherly solicitude for my health and welfare, during the whole of the time I was in Rome, filled me with gratitude, and will never be forgotten.
Your house, as I predicted it would be, was as a home for me; and with deep but still satisfaction I felt that I was treated by every member of your family as an old friend, who had come to them, broken in health, to be cheered and comforted by their kindly sympathy. I need not say that I found all I sought from them, and far more. Should I ultimately recover from the nervous depression under which I have so long and so painfully suffered, you may rest assured that I shall ever attribute such recovery, in great measure, to the kindly and active sympathy of you and yours. Accept then, my dear friend, my warmest and deepest thanks.
You asked me about my niece.
The news which reached me just before I left you agitated me deeply, as you saw. Her welfare and prospects have filled my whole thoughts ever since almost daily, since I left you, I have either heard from, or written to her. I have been intensely anxious to render to her that liberty of action and decision which was so manifestly impaired by her tender and deep attachment to me. I have sternly habituated myself to the thought of separation from her, and at length I have been brought to exercise my powers of persuasion, and to advise her to accept the offer that has been made to her that she has at length done and thus the matter is virtually settled. I propose now to go to Yorkshire, and to spend the remaining part of my holiday with her; the sight of her happiness will be cheering to me, and I shall have an opportunity of making the purposed acquaintance of the man to whom, henceforth, her welfare is to be entrusted
All accounts I have received agree in representing him as a man of high moral character and integrity. Thus it has happened that my forced absence, on account of illness, has been turned to hers good, and one ground for fear (viz[1] that her devotion to me might impair her own prospects in life) has been happily removed from me. I am now beginning to rejoice at this result, and to overlook the pain that separation from her will necessarily cause me.
I spent a very pleasant day in Pisa with Betti Bertini and Dini; and another almost equally pleasant one with Genocchi and D'Ovidio in Turin
Besides the two latter I also saw that young Professor you mention (an artillery official), but whose name I have for the moment forgotten.
I delivered your kindly greeting to all three.
And now dear friend Goodbye. Remember me ever kindly as I shall you

Affectionately yours
T. Archer Hirst

P.S. I hope you received your Guide Book safely.
I posted it in Turin.

[1] *Viz* sta per *videlicet*, cioè "vale a dire".

24

Shafton
near Barnsley
Yorkshire
May 10th '77

My dear Cremona

I send you a line to inform you that I reached London a week ago safely, and on Monday last came on here to spend a week or ten days with my niece, before returning finally to Greenwich.

Her marriage is now practically settled and I have every reason to believe that she has a happy life in prospect. She desires me to thank you for your little present. She will keep it as a souvenir of your pleasant visit to us at Greenwich, she is greatly pleased to know that she was remembered so gracefully by you. On passing through London I saw the Tyndalls and Debus; they listened with pleasure to my account of the brotherly kindness you showed me in Rome, and they desired to be kindly remembered by you

Give my best wishes to Mrs Cremona and all the other members of your family

They will be glad to hear that I am better both on health and spirits than I was when I bid them adieu

Ever yours sincerely
T. Archer Hirst

P.S. Please forward the enclosed letter to Prof. Jung for me.

25

Darmstadt
August 13th 1877

My dear Cremona
Your letter of the 27th of July reached me in Greenwich when I had no time to answer it so busy was I with official work. Not until the evening of the 11th was I able to get away from Greenwich and come here. I arrived only yesterday evening
I presented myself at once at Sturms where I met with a hearty welcome. It is my intention to remain here very quietly and if possible to continue my geometrical work at which I have done absolutely nothing since I saw you. Of course I should prefer working at your side but you are too far distant and travelling disturbs my health far more than it used to do. Moreover I fear the climate even in northern Italy would be too relaxing for me and lastly I must return to Greenwich at or about the time when the meeting at Munich will commence. At the same time if my health should improve towards the end of this month I should be very glad to spend a week in your company, my dear friend, for as you know there is no one whose companionship I appreciate and desire so much Perhaps we might meet in the neighbourhood of Salzburg which is very beautiful where I could remain with you up to the period of the meeting at Munich which city is not far distant. Let me hear from you at all events before the end of month (address your letter to Sturm's care) and tell me what your plans are. I was delighted to hear there was a possibility of your going to Pisa and glad to know that you had spoken with Betti on the subject. I knew he would be prepared to further your wishes for I had a long conversation with him about you when I was in Pisa and I at once saw that he sympathized with your desire to devote yourself exclusively to Science. Pisa is one of my favourite places It is so quiet, and so *[…]* to undisturbed scientific work. I shall be greatly interested to hear how your prospects of leaving Rome are becoming realized. Do not forget that you are to come to England next year from Paris. I am the more anxious that you should decide early to do so because Spottiswoode is to be President of the British Association whose meeting is to be held at Dublin. I expect it will be a specially good meeting for Mathematicians. I have promised Spottiswoode to go and support him and to try and induce as many continental Mathematicians as possible to be present on the occasion
In reply to your kind enquiries about my health I am glad to be able to tell you that I am decidedly better than I was in Rome. There is, it is true, no radical change in my symptoms but my general health is better and my nervous disorder less distressing. My niece I am sorry to say has been very ill with a severe cough which threatened to attack the lungs. She came to see me at Greenwich before I left and I left her there with her mother. I have hopes that the change of air will be of service to her. Her marriage has been postponed in consequence of this illness, but I trust after her recovery it will not be long delayed. As my own health declines I shall look forward with more and more satisfaction to her settlement in life.
Give my very cordial remembrances to Mrs and Miss Cremona as well as to Vittorio and Itala. I picture to myself the great pleasure which you have secured for your children by bringing them to the beautiful alpine scenery in northern Italy. Sturm and Mrs Sturm speak in the kindliest terms of you and desire to send their hearty greetings.
Receive the same from me dear friend. Ever yours

T Archer Hirst

I *[…]* mention that I have not yet received the parcel of books which you kindly undertook to send from Rome to London

## 26

Darmstädter Hd
Darmstadt
Sep. 2_ 1877

My dear Cremona
I have postponed my reply to your letter of the 24th until I could decide upon the possibility or not of meeting you at Salzburg. Much as I should have enjoyed passing a few days in your company I hesitate in my present state of health to undertake so long a journey <u>for so short a time</u>. If I were to attempt it I know that I should scarcely have recovered from the effects of the fatigue before I should have to commence my long journey to Greenwich. Under the circumstances, therefore, I feel reluctantly obliged to give up the project. I look forward with some hope to an improvement in my health before you come to England next year. Although the improvement is very slow still I can say that my nervous system is decidedly less disordered than it was in Rome. The air and, above all, the quietude of this place suit me very well. I propose to remain here about a week longer and then to return by easy stages to Greenwich. Good bye then my dear Friend. May your hopes about Pisa be speedily fulfilled, and may you enjoy yourself thoroughly in Munich
Sturm will bring you my greeting and the most recent report of my doings. With kind regards to Mrs Cremona and your children
Ever your sincerely.

T. Archer HIrst

My niece's health is much better
I trust nothing more will now delay her marriage

## 27 [1]

Greenwich
Octr – 3rd/ 77

My dear Cremona
The parcel of Books you were kind enough to forward for me to London has just arrived. Pray accept my best thanks. The College having just been opened I am very much engaged. I must postpone, therefore, writing you a letter. Let me hear from you soon. Did you enjoy your visit to Munich? Kind regards to Mrs Cremona and your family

Ever yours T Archer Hirst

[1] Cartolina postale indirizzata a "Prof. Cremona – R. Scuola degli Ingeneri *[sic!]* – S. Pietro in Vincoli – Roma – Italy".

28

Royal Naval College
Greenwich
Nov. 26th 1877

My dear Cremona
Your letter of the 7th of November interested me greatly. I was especially gratified to hear of the arrangements that had been made (and made in so satisfactory a manner) for diminishing your work at Rome. I rejoice to think that you will now have more leisure to pursue your studies.
Brioschi was here last week and was speaking to me about you and your prospects in a manner which gave me great pleasure and showed great appreciation on his part I dined with him (and Menabrea)at Spottiswoodes on Friday last. I expected a visit from both tomorrow but I heard yesterday that Brioschi had departed suddenly for Paris. I gave Spottiswoode your message and asked Mrs S.[1] to send you the promised portrait. I enquired about your communication of January last to the British Association. I found that it had been neglected (amongst many other things) by the Secretary. I wrote to him to request that your letter might be brought before the Council at their very next meeting so that I trust you will receive a favourable reply very soon. Mr Griffith is about to resign and I hope his successor will attend more punctually to correspondence
I have asked the Secretary of the Mathematical Society to forward to you the missing numbers of our proceedings. With respect to my own Journals I find that I have all the numbers of the Giornale di Matematiche up to the present time. With respect to the Annali my collection is complete with <u>one</u> exception Fascicolo 3° of Tom[2] *[sic!]* vii is missing. The last number I received was Facicolo 3°, Tom *[sic!]* viii. The last number of the "Notizie degli Scavi" which I received was the one for July. The only missing one is that for January 1877.
You do not say anything about the Quarterly Journal of Mathematics. 37
The last number I received was N° 55 issued in Oct. 1876! I have just written to ask if any have appeared since and if they have to request that they may be sent to us both without delay.
I come next to the question of my own health about which you enquire so kindly. I am decidedly stronger and better than when I left Rome
I am still unable to use my right hand and leg properly and I begin to fear I may never regain their proper use. I am however far less distressed with that nervous disorder under which I suffered during the first half of the present year.
My niece is at length married and is very happy. I have not seen her since last August but she writes to me frequently and often enquires about the friends at Rome who treated me with such *[…]* kindness when my health was so broken. Give my sincerest good wishes to Mrs Cremona and your family congratulate them from me that they are now able to remain at Rome under more favourable circumstances
Since my niece left me my house here has entirely changed in character. I spend very little time in it. All day I came at work in my office and I spend my evenings either at the Atheneaum[3] *[sic!]* Club or at the Tyndalls, who by the bye are both well and happy and return your kindly greeting with warmth. Cayley comes to stop with me for a couple of days every month when the Math. Society meets
I told him of your communication at Munich. We shall both be glad to see it when it is printed. I received your memoir on Pascals Hexagram[4] and sent the other copy to the Math Society With kind regards to Battaglini Dini Cerruti & others

Yours sincerely
T. Archer Hirst

---

[1] Mrs Spottiswoode.
[2] Tomo.
[3] Athenaeum.
[4] L. Cremona, "Teoremi stereometrici dai quali si deducono le proprietà dell'esagrammo di Pascal", *Mem. della R. Acc. Nazionale dei Lincei,* (3), 1, 1876-77, pp. 854-874.

29

Greenwich
May 17th 1878

My dear Cremona
I was very pleased to receive your letter yesterday for you have been in my thoughts very frequently since Easter, when I paid a short visit to Paris; and saw our friends Chasles Mannheim Halphen, Fouret and Darboux. I am glad to hear that although you have been too busy to write to me your occupations have been of so much more congenial nature than they were last year. Science, I feel sure, will profit as much from the change as you yourself have done The subjects you name as having occupied your attention are indeed wide in range and most fascinating in character. Of many of them - for instance of Klein's investigations - I know absolutely nothing, so completely do my administrative duties absorb my powers of study and leisure for reading. When we next meet, however, I shall trust to you to give me an idea of the nature and scope of these new researches.
My health and spirits are decidedly better than they were a year ago, but I am still far from being strong, and it is only by the observance of the greatest care in diet, in gentle exercise and in avoidance of over-fatigue and mental worry that I have been able to accomplish the work I have had necessarily to perform. As to my own geometrical work I can make no progress whatever; the utmost I can do is to prevent my own subject from becoming a foreign one to me, - to retain my hold upon it, in short, so that when more leisure shall be granted me I may take a few more steps in advance.
June and July will be very busy months for me; I expect to be greatly fatigued by them and to be far from possessing the bodily vigour which I should like to be able to place at your disposal in August. Nevertheless my intention is to go to Dublin, and the fact that you have decided to come tends strongly to strenthen[1] *[sic!]* my purpose. Give me early indication of your plans and try and arrive in London as early as you can in August, so that we may spend some pleasant days together before we proceed to Dublin Of course you will take pity on my loneliness and come to me entirely at Greenwich. You shall have my best suite of rooms all to yourself, and come and go with perfect freedom. How welcome you will be to me my dear Friend! Is there no chance of your bringing Mrs Cremona or Miss Cremona with you. *[punto di domanda]* The Association Française have arranged to hold their meeting in Paris immediately after that of the British Association shall have terminated; will you go to it? I shall not do so. My present intention is to spend what remains of my short vacation after the Dublin meeting with my niece in Yorkshire
I have not seen her since she was married. She is very well however and, as far as I can judge, very happy. I send kindest greeting to all members of your family and to all mutual friends

Ever yours sincerely
T Archer Hirst

[1] Strenghten.

30

Shafton
Near Barnsley
Yorkshire
July 3rd 78

My dear Cremona
I was very sorry indeed to receive your last letter announcing the death of your Brother and the consequent abandonment of your intention to come to England next month. I never saw your brother, I think, but you have often spoken to me about him and I can well understand that his death must be a source of great grief to you
Your change of plan will possibly cause a change in mine . I feel no great desire to go to Dublin
In fact the prospect of having you for a companion and the desire to give my support to Spottiswoode were my sole reasons for going, the meetings of the British Association are, I confess, no relaxation for me I would far rather spend my vacation elsewhere and in more congenial company
One of my two reasons for going to Dublin being withdrawn I am not yet sure that the remaining one will be strong enough to overcome my natural reluctance to sacrifice a couple of weeks of my short vacation to so unprofitable a purpose. However I am unable to decide at present and must wait for further motives. The session at Greenwich closed on the 29th of June and I came here at once for a few days rest and for the purpose of paying my first visit to my niece since her marriage. She is very well and very happy; but of course her occupations and interests have been greatly changed since she left me and I have not yet quite reconciled myself to the change. She thanks you for so kindly remembering her and send you her best wishes
I am glad to hear that Miss Cremona is making such a good progress with English, and still more glad that you entertain the hope of bringing her with you to England next year My house is not as well suited as it once was to receive her as a guest, but if she will accompany you and regard it as the house of an “old uncle” where she can have perfect freedom and a hearty welcome she will, by making it her home, give me the very greatest of pleasures.
I have no mathematical news to send you I am looking forward to being able to resume my geometrical work for a short time before next october[1] *[sic!]*, but my work at Greenwich will not be quite finished until about the 8th of August.
May I ask if you have received all numbers of the Quarterly Journal of Mathematics up to N° 60 inclusive? I ask because I have recently changed my bookseller and am not sure that he has executed my orders. With kindest regards to all your family Ever yours

T. Archer Hirst

[1] October.

31

*[su carta intestata:*
Royal Naval College
Greenwich S.E. *]*
25th March 1879

*[in testa alla lettera l'autore ha aggiunto queste righe:* You have no doubt heard of Clifford's death. I heard of poor prof. Chelini death with great sorrow*]*

My dear Cremona
It is a long time since I had any news of or from you. The reason I trust is simply that, like myself, you have been too much occupied
Between your occupations and mine, however, I hope there has been one important difference; viz: that whereas my occupations have been official ones chiefly, you have found leisure to devote yourself to science.
I break silence to day, however, for the purpose of congratulating you on the fulfilment of a wish I have long entertained
You will have probably received from our excellent President Spottiswoode intimation of your election as Foreign Member of the Royal Society of London. Strictly speaking perhaps I ought to say you have been selected by the Council for Membership, for I am not quite sure that the formal election by the Fellows of the Society took place on Thursday last or is to take place on Thursday next. This however is a mere formality. I trust, my dear friend, that this recognition, on form part, of the high value of your geometrical researches, will be a source of pleasure and satisfaction both to yourself and to all your friends in Italy.
My health during the severe winter we have had has been far from good; but still I am much stronger than I was when I last saw you in Rome. During the few weeks of holiday I had in September last in Yorkshire I began to write my memoir on Correlation in space, but as soon as the session opened in October I had to put it aside and I have not been able to advance it since. During any brief period of leisure I have merely been able to look at an isolated question or two more or less intimately connected with the main one Amongst other things I have had occasion quite recently to read again your memoir in Crelle[1] (tome 60) "Sur en surfaces gauches du troisiénne degré[2]" By the bye is there not some mistake in Cayley's construction of the cubic scroll whose single and double directrices coincide, as given by you on the first page of the memoir I have just referred to?
If I am not mistaken, in fact, you have yourself shown in your memoir on Quartic Scrolls[3] published by the Academy of Bologna about the year 1868 or 9 that the surface to which the above construction by Cayley leads is of the fourth degree. The matter, it is true, is of no importance. The error, if error it be, does not at all affect the rest of your memoir. With kind remembrances to all members of your own family and to all friends in Rome
Ever yours sincerely

T. Archer Hirst

[1] Il giornale *Journal für die reine und angewandte Mathematik* fondato nel 1826 da August Leopold Crelle divenne presto conosciuto con il nome di "giornale di Crelle".
[2] L. Cremona, "Sur les surfaces gauches du troisième degré", *J. für Math.*, 60, 1862, pp. 313-320.
[3] L. Cremona, "Sulle superficie gobbe di quarto grado", *Rend. dell'Acc. delle Scienze di Bologna*, 1868-69, 96-97, *Mem. dell'Acc. delle Scienze di Bologna*, (2), 8, 1869, pp. 235-250.

32 [1]

Hotel Mirabeau
Rue de la Paix
A Paris
6th April 1879

My dear Cremona
I would gladly contribute to your collection of papers in Memory of Chelini but I have nothing ready, and I should require time; the more so since I have come here (tired) for a holiday. Tell me; how long you can give me Will it suffice if I send you a small note in the course of a month? I assume that the note will serve your purpose if written in English. I was surprised to receive a letter from you written in such good English. You have indeed made good progress since we last met. Remember, however, that although it may be a good exercise for you to write in English, it will be a loss for me to be deprived of your letters written in Italian, letters which I always prized greatly. Clifford was only about 34 years of age
His loss is indeed a great one. I shall be delighted if I can meet you in Switzerland next Autumn.
I hear with pleasure of your election as a Senator.
Give my kind remembrances to all the members of your family, and all friends. Ever sincerely yours
T. Archer Hirst

33

15 June 1879

My dear Cremona
You would pity me, I am sure, if you knew how I have been struggling for the last two months under unfavourable circumstances to send you something for Chelini's Memorial volume. I told you from Paris that I had nothing ready; but I set to work at once to complete a little investigation I had often thought about "on the complexes generated by joining the conjugation points of two correlative planes"[2] unforeseen little difficulties presented themselves, I did not like to send you anything so long as their remained points I had not perfectly cleared up, and my Greenwich work deprived me of the leisure and intellectual freshness & vigour necessary for the completion of my task. It is, however, at last, virtually, done but it is too extensive, I fear, for your purpose
It would occupy more than 32 octavo pages. Accordingly a doubt has arisen in my own mind which I should be glad if you would settle. Would it not be better that I should send you an abstract of my paper containing results but no demonstrations? This could be done in about ten pages and I would publish the complete paper in the Proceedings of the Math. Society subsequently. If you would prefer having the complete paper however you shall have it early in July. I could not promise to send you anything earlier for the session at Greenwich being near its close I am overwhelmed with work.
I am still looking forward to meeting you in Switzerland.
I have no choice of place but I will endeavour to come to you whenever you may go.
With kindest remembrances to all your household believe me to be
Yours sincerely
T. Archer Hirst

[1] Cartolina postale indirizzata a "M. le Prof. Cremona – R. Scuola degli Ingeneri *[sic!]* – S. Pietro in Vincoli – Roma l'Italie".
[2] T.A. Hirst, "On the complexes generated by two correlative planes" in: E. Beltrami, L. Cremona (a cura di), *Collectanea Mathematica (In Memoriam Dominici Chelini)*, U. Hoepli, Milano, 1881, pp. 52-73.

34

Shafton
near Barnsley
Yorkshire
11th July 1879

My dear Cremona
By this post I send a registered letter containing my contribution to the Chelini Memorial volume. I have been able to finish it during the weeks[1] *[sic!]* rest I have enjoyed here since the College Session ended.
I return tomorrow to London for a months[2] *[sic!]* hard work and then I hope to be able to start on my real holiday. Let me know before the 6th of August if you decide upon going to the Val d'Aosta and I will try to join you there
I have endeavoured to make the M.S.S.[3] of my paper legible but I should be glad to have proofs sent to me for correction
You will notice that I intend all symbols for complexes and congruences to be in thick letters, which we call block type. If I remember right you used such characters in one of your own papers on complexes
I have numbered the foot notes consecutively and carefully distinguished them from the text, there is only one place in fact (page 35) that can cause any confusion to the printer, but you will certainly understand it
Yours very sincerely

T Archer HIrst

I remain here about 10 days

35 [4]

Munster
Friday Aug. 15/79

I have just received your card. Tomorrow Sturm and I go to Bremen to meet Schubert. On Monday I shall travel to Frankfurt. On Tuesday I shall proceed on my journey towards Friedrichshaven[5] on the Lake of Constance. Thence I shall go to Coire[6] at which place I may arrive on Thursday the 22nd but certainly not before that date
Please send me a Post card to Coire (Poste Restant[7] *[sic!]* and inform me with respect to your movements
Prof. and Mrs Sturm desire to be kindly remembered to you. Sturm joins me also in a friendly greeting to Casorati and Beltrami
We have just received the news of Jung's marriage
Ever yours sincerely

T. Archer Hirst

[1] week's rest.
[2] month's hard work.
[3] M.S.S. – manoscritti.
[4] Cartolina postale indirizzata a "Professor Cremona – Schweiz - [*incomprensibile*] – (Poste Restante)". Dal timbro postale si può desumere che fosse indirizzata a Davos.
[5] Friedrichshafen.
[6] Coira.
[7] Si tratta di "Poste Restante" che indica, in francese, il servizio di "fermo posta".

36 [1]

Axenstein
Sep.$^{r}$ 15/79

My dear Cremona
I received your letter yesterday and shortly afterward Clausius came and persuaded me to go for a day or two to the Naturforscher Reunion at Baden Baden which commences on Thursday. I told him I should wait at Lucern until you came and that if possible I would persuade you to accompany me to Baden Baden
I shall expect you at Lucern, therefore, on Thursday morning and at all events we will spend that day together. I propose to stay at the Lucerner Hof. Kind regards to Geiser and Frobenius Yours ever
T Archer HIrst

37

Greenwich
15 Oct 1879

My dear Cremona
By this time you will, no doubt, have returned to Rome invigorated, as I find myself, by your holiday in Switzerland
After leaving you at Vitznau and returning to Lucern I telegraphed to Reye to ask if he was at home and received a reply in the affirmative. Accordingly I went to Strasburg and spent a couple of days very pleasantly there. During my stay the "Naturforscher" from Baden Baden visited the city and I had the pleasure of meeting Schroeter (from Breslau) Nöther, Clausius, Christoffel and others.
At Cologne on my way home I was pleased to encounter the Wolff family and I travelled to London with them via Flushing
On the whole my holiday this year in Switzerland was a more successful and agreeable one than any I have enjoyed for many years, and as I *[sic!]* consequence I find my health greatly improved. It was a great delight to me to have your companionship and I hope we may be able to meet again next year.
The only Mathematicians I have seen since I returned are Spottiswoode Smith and Glaisher. All three were glad to receive my favourable account of you and desired to be kindly remembered when I write to you. I hope Casorati and Beltrami reached home safely and well; please convey my good wishes to them when you have an occasion for so doing.
I have received a proof of a note I communicated to the Math. Society last June. I should like in it to refer to my memoir in the Chelini collection. What title do you propose to give to the forthcoming volume, and in what manner can I most conveniently refer to it?
With kindest regards to M$^{rs}$ Cremona and all your family
I remain, my dear friend,

Yours very sincerely
T. Archer Hirst

[1] Cartolina postale indirizzata a "Prof Cremona bei Prof. Geiser Zeltweg – Escherhaüser, 15 – Zurich".

38 [1]

[2]

My dear friend
I take the advantage of Mr Guccia's coming to England, to procure him the pleasure of your acquaintance. Mr Guccia is one of my cleverest and most assiduous scholar; he is especially fond of pure geometry.
This summer I cannot go abroad; in my next letter I will inform you more exactly about my plan for the holidays, and the condition of my family.
My wife is not yet quite well.

Yours sincerely
L Cremona

39 [3]

Oct. 1880

My dear Hirst
I am still your debtor for your kind letter of the 30 August from Venice and your postcard of the 17 Sept. from Paris; by which I was very glad to learn that you found a successful subject of geometrical investigations in that which you kindly call Cremonian congruences. I look forward, as an *[sic!]* miserable Tantalus, to the day when I shall be able to retake my own studes[4] *[sic!]*, without any distraction. But I am very sorry that the bad weather spoiled all your travelling plans; and mainly I deplore that my own extraordinary duties made impossible our meeting when you were in Italy. My wife and children are still at the Bagni di Lucca, all in good health, as myself. I paid them a visit from time to time; but in the lump I was obliged to spend the summer in Rome. Only about the middle of September I took the opportunity of an excursion to the Etna, made by the Alpine Club, in order to enjoy a little rest, united to an energical *[sic!]* exercise. In this manner, I spent 17, 18 and 19 Sept. on the great volcano; and a few days before and after in Cata~~nia~~, Taormina, Aci Reale, Catania and Siracusa: the cities of Sicily that I had never seen. Aci Reale has a very excellent hotel (Grand Hotel des Dains); it must be a very delicious and salubrious stay by the winter season. At the 24 Sept. I was again in Rome.
My extraordinary duties are, I hope, not far off an end by laying down the powers of the Commissariat I must deliver up a Report and I fear that such a Report will be for me a serious work.
In order that you may understand with what sort of *[…]* I have to do, I will send you a copy of the Report of the Inquiry Commission on the library V.E. The inquiry took place before my nomination at a plenipotentiary Commissary, as you will see from my letter being a proem to the lately named Report.
Also in my daughter's name, many thanks for your friendly wishes regarding our prospect of coming to England next year. We warmly desire that nothing may prevent the realization of this plan.

*[nota a matita:* non finita e non spedita*]*

[1] Lettera di Cremona a Hirst.
[2] Dalla corrispondenza tra Guccia e Cremona (Istituto Mazziniano di Genova) si può dedurre che la lettera sia stata scritta all'inizio dell'estate del 1880: nella lettera datata 22 luglio 1880 (054-13058), infatti, Guccia avvisa Cremona di non essere andato a Londra contrariamente a quanto annunciato in precedenza e poi il 10 ottobre dello stesso anno (054-13061) scrive che si trova a Londra da 10 giorni in compagnia di Hirst. Si veda C. Cerroni (a cura di), *Il carteggio Cremona-Guccia (1879-1900),* Mimesis, Milano, 2014.
[3] Lettera di Cremona a Hirst.
[4] studies.

40

*[su carta intestata:*
Athenaeum Club
Pall Mall *]*
18 May 1881

My dear Cremona
I cannot say I was surprised, but I may truly say that I was gratified, to receive from you the intimation that your eldest daughter, Elena, whom I have known from infancy, was married
Pray convey to her my warmest congratulations and my best wishes. I trust that the marriage is one which gives entire satisfaction both to yourself and to M^rs^ Cremona; but even should this be so you will both, I am sure, have felt her departure from your home to be a momentous occurrence The first rupture in a family circle must necessarily be so. Although I never had a child of my own I can never forget how deeply I was affected by the intelligence of the approaching marriage of my Niece, which reached me when I was last your guest at Rome. Well, although the brightness of my home departed with her, I have come to regard my loss as her gain. She is happy, and her future is secured, whatever may happen to me. I sincerely trust that this experience of mine will be repeated in your case.
I fear that one effect of your daughter's marriage, however, will be to alter your plans for paying us a visit during the Meeting of the British Association at York. I remember that it was your intention to bring your daughter with you. What do you propose to do now? I trust you will not abandon your own project
The meeting is fixed for the 31st of August. Early in that month I shall probably leave London and go to Yorkshire to pass a short time in the neighbourhood of my niece, who may possibly accompany me to York. Should you come we might possibly arrange to spend a little time together in the North of England, and after the meeting, which will terminate on the 7th of September, I hope you will return with me to Greenwich and make my house your home during the remainder of your visit
In any case I hope to hear from you before long and to be made fully acquainted with your intention
I trust Mrs Cremona's health is better, and that your own continues good. I see your friend Sella is trying to form a ministry. Should he succeed he may possibly retain your services in some capacity or other. On this point I have but two wishes to express: <u>first</u> that he will not interfere too seriously with your scientific activity as a Geometer and <u>secondly</u> that he will not change your views as to paying us your promised visit.
Ever sincerely yours
T. Archer Hirst

41

*[su carta intestata:*
Royal Naval College*]*
24 July 1881

My dear Cremona
Your long expected letter reached me a few days ago, when I was on my tour of inspection of Dockyard Schools. It was very welcome, as all your letters are, though it did contain the disappointing intelligence that you had decided not to come to York. I regret this the more since I shall be unable to leave England this year, and therefore our meeting in Switzerland is out of the question. I was deeply interested in the account you gave me of Madame Cremona and all your children. To the Signora Perozzi Cremona, my former friend Elena, please give my warmest good wishes and kind regards. Thank her also for the photograph she has kindly sent me
It is the third, of her, that I have placed in my album
The first of the three represents her as I saw her at the School to which she went in Milan (I think)
The most satisfactory intelligence I found in your letter was that concerning the improved health of M$^{rs}$ Cremona
I sincerely trust it may continue to improve
Next to that, in point of interest for me, was the assurance that you will shortly be able to return to your Geometrical studies. I have no doubt you have rendered valuable services to your country during the past year, apropos of that "Biblioteca Vittorio Emanuele", but I can assure you that your scientific friends in Europe consider that the sacrifices you have had to make have been unduly great. They, and I in particular, will welcome you back again to Science
Many thanks for the "Monografia di Roma" it will have very great interest for me
I enclose the last photograph which has been taken of myself. Please convey it to the Signora Elena with my heartiest congratulations

Ever sincerely yours
T Archer Hirst

42

Devonport
5.Jan.1882

*[in testa alla lettera l'autore ha aggiunto queste righe:* I am glad to hear better accounts of M$^{rs}$ Cremona; pray remember me to her most kindly. The news that you are so soon to be Grandfather startled me. I shall be anxious to hear of the happy consummation of the Mothers[1] *[sic!]* hopes
Pray convey my very best wishes to her *] [continua in una nota a piè pagina:* My niece is well and happy but is yet childless. She was with me at York*]*

My dear Cremona
The receipt of your letter yesterday was a source of the greatest pleasure to me and although I am on voyage I will at once reply to it and send my hearty New Years to yourself and all the members of your family now around you.
I was disappointed in August last that not one of the Mathematicians who were present at the Parpan[2] Congress in 1879, made their appearance at York, as they had induced me to hope they would. Zeuthen, it is true, accepted the invitation; but his mothers[3] *[sic!]* illness, which terminated fatally I believe, prevented him at the last moment from coming.
The only foreign mathematicians there were Sturm, from Münster, and Halphen and Stephanos from Paris, the former of whom was accompanied by his friend and interpreter Chemin, of the Ponts et Chaussées[4]. Of English Mathematicians there was a fair collection. Cayley (since gone to America to give a course of Lectures at Baltimora) Smith, Spottiswoode, Glaisher Merrifield, Price (of Oxford) Ball (of Dublin) Harley and a few younger men from Cambridge & Oxford; not to forget Miss Scott, who is a Professor of Mathematics at Girton College (?) near Cambridge.
It fell to my lot to preside over the Mathematical Department, which sat for three of the six days, and received more communications than on any previous occasion. I had good reason to feel disappointed, therefore, that my old and dear friend Cremona was not present
I am very sorry to learn from your letter that that episode of the Victor Emmanuel Library ended so unfortunately for your peace of mind. I regretted its intervention from the very commencement, as you will remember: but it is a satisfaction to me to reflect that it is now past, and that there is nothing more than your ordinary duties of Director and Professor to divert you from your own researches
May you resume them with your former success and with more matured, if less energetic thought
I still occupy my rare moments of leisure with my geometrical studies on Congruences, but I cannot procure the necessary amount of continuous leisure to collect and publish my results
I still hope that less distracted days may dawn upon me before I lose the power of utilising them.
Ever sincerely your friend T Archer Hirst

[1] Mother's hopes.
[2] Ex comune svizzero del Cantone dei Grigioni.
[3] mother's illness.
[4] L'École des Ponts Paris Tech fondata nel 1747 con il nome di "Scuola reale di ponti e strade" (École nationale des ponts et chaussées).

43

*[su carta intestata:*
Royal Naval College
Greenwich S.E. *]*
22 July 1882

My dear Friend
The melancholy intelligence of the death of M[rs] Cremona reached me yesterday, and has filled me with sorrow. I grieve not only for the loss of one for whom I have long entertained a sincere esteem, but also for the irreparable loss which you, my dear friend, have sustained
The sad news did not take me by surprise
Even since I received your letter of the 28[th] of May I have felt that any post might bring the news that the sufferings of your poor wife had terminated fatally
Knowing how intense these sufferings must have been, I have even felt that it would be a mercy if death soon arrested them, seeing that it is beyond the power of medical science to do so.
Although no word of mine can now diminish your own affliction I trust that time will soon allow you to derive some consolation from the thought that you have my deepest sympathy as well as that of all your friends
When I remember the salutary effect that change of scene had on me, when I had to pass through a trial similar to your own, I cannot refrain from advising you to travel as soon as circumstances will allow you to leave Rome. A journey northwards could not fail to be of great service to you now, and I hope you will be induced to make it
I have not yet decided when I shall spend my vacation. If there were any chance of meeting you in the Alps I should be strongly tempted to go. I expect to be able to start in about a fortnight. Give my kindest regards to Elena, Vittorio and Itala and believe me to be, as ever, your sincere friend

T. Archer Hirst

P.S. (July 23) I will answer the question you ask in your letter of May 28. I have seen no notices of the Mathematicians in the “Biograph” than Harley and myself. I did not send you the last Report of the Association for the Improvement of Geometrical Teaching. Previous reports are probably now out of print *[l’autore ha inserito una nota:* I will ask the Secretary to send you all he can*]* I have not seen the two “libretti” of Prof Pincherle to which you allude, but I should be glad to read them.
My official work and delicate health have prevented me from making much progress during the past year with my geometrical studies.

44

*[su carta intestata:*
Royal Naval College
Greenwich S.E. *]*
5 Aug$^{st}$ 1882

My dear Cremona
Your letter of the 31 July has just reached me. I learn from it with the greatest satisfaction that you are at Zuz in the Engadine and that Vittorio and Itala are with you. The change of air and scene cannot fail to be of great service to you. I am also very glad to hear that Professors Beltrami and Casorati are with you; pray give my kind regards to them. I did not tell you in my last letter that I was recovering from a surgical operation to which I had to submit myself on the 3$^{rd}$ of July. My recovery is progressing favourably but my surgeon advises me not to travel too far; so that I fear I must give up my intention of meeting you in the Engadine. On Monday I propose to leave Greenwich and cross over to Dieppe solely. Should I regain my strength sufficiently I may go further South but I shall not quit France.
I read your letter with the deepest interest and sympathy
You have my warmest good wishes
Ever yours
T. Archer Hirst

45

*[su carta intestata:*
Athenaeum Club
Pall Mall*]*
28 July 1883

My dear Cremona
I have no less than three letters from you to which I ought to reply. The first is dated June 16 and was written in answer to a post card of mine in which I first informed you of poor Spottiswoode's serious illness
The next was written on the 30th of June after you had received my second card announcing the fatal termination of that illness
So great was my grief at that catastrophe that I had no heart to write more about it. I felt as if I had lost a brother, as well as a friend whose friendship I could never replace. I have still cause to feel his loss every day, and shall long continue to do so. Mrs Spottiswoode is bearing her still greater grief with great courage and resignation. At present she is with her eldest son at Combe Bank, which they will probably let until that son has finished his ministry education.
Your sincere sympathy was a great consolation to her throughout her terrible trial.
I shall look forward with melancholy interest to your notice of our friend in the proceedings of the Accad. dei Lincei[1].
With respect to your last letter written, of the 24th of this month, I was glad to hear that you were at Portomaurizio with Itala. I hope the Sea Bathing will refresh you both and that you will have a pleasant meeting with your other children in September. I have asked the Librarian of this club, who is a very competent authority, to answer your friend questions about the Newspaper Room at the British Museum. I enclose his answers which I trust will be found satisfactory. I am very busy (and very weary) arranging all my affairs at Greenwich preparatory to to *[sic!]* my retirement on the 31st of August I do not yet know, precisely, what I shall do afterwards. If I am well enough I will attend the meeting of the British Association in mid-September at Southport (in Lancashire) act of respect to Cayley who is to be President 50
Immediately afterwards I expect my niece will become a Mother
It will be her first child, although she has been married nearly five years; and her health being far from good I am anxious, I may say nervously anxious, about the result. I cannot leave England until that result is known; but should the crisis be happily passed I may possibly pass the winter in a warmer climate than ours. But all this is a mere project and very uncertain.
I must write to you again as to my probable movements of course I have thought of Rome, but I have formed no definite plans yet. My retiring pension, in fact, is not yet determined and upon that, of course, my mode of life hereafter must greatly depend. I have no cause to complain much of my present state of health; it is in fact better than it was earlier in the year
When once I am relieved of the busy tedious and time-absorbing duties of my position here, I hope to turn my attention more continuously to the geometrical studies which have been interrupted for so long a time. My first effort will be to publish the researches on (Cremonian) Congruences of the second order with which I have been engaged at broken intervals during the last few years. But the task will be no light one for as you observed in one of your recent letters I shall have to study again subjects with which I was once conversant but have now almost forgotten

Ever my dear friend
yours sincerely T.Archer Hirst

[1] L. Cremona, "Cenno necrologico di W. Spottiswoode", *Atti della R. Acc. Nazionale dei Lincei,* (3), 7, 1882-83, pp. 308-309.

46

*[su carta intestata:*
Athenaeum Club
Pall Mall S.W.*]*
19, XI, 83.

*[in testa alla lettera l'autore ha aggiunto queste righe:* M^rs Spottiswoode is fairly well.
My own movements are uncertain
I must write to you again.*]*

My dear Cremona
The day before yesterday I sent you a copy of "Nature", feeling sure that it would interest you to see that the Royal Society has awarded to me one of the Royal Medals of the year. I have to day to inform you of another event, which happened long before this award, and which has greatly diminished, not only my interest in what the Royal Society have been good enough to do for me, but my interest in life altogether.
My dear Niece whom, as you know, I loved as if she had been my own daughter is dead! Last September I went to see her in Yorkshire, and spent ten happy days in her company. She was about to become a Mother; and it was beautiful to me to see her happiness at the prospect of her maternity, and her utter fearlessness of the ordeal she was about to pass through. It is true I had my own misgivings, knowing, as I did, how feeble her health had been during the past year. Early in October, shortly after I had bid her adieu, she was confined, and gave birth to a fine boy. Her sufferings, however, appear to have been intense; for 58 hours she was in agony, and after her delivery she sank, exhausted, and live no more!
Her boy still lives, but she was hardly conscious even that it had survived the ordeal from which she was never to rally!
You remember how her betrothal, when I was last in Rome, almost broke my heart. I remember well how you, and poor M^rs Cremona tried to cheer me under the loss I then sustained. My present grief is not more acute than that one was. This one, however, is irremediable. The former one became tempered by witnessing her own happiness.

Ever dear friend
Yours affectionately
T Archer Hirst

*[sul retro della lettera si trova una nota a matita non scritta da Hirst:*
Dal Nature, 15 nov. 1883
La Società Reale ha aggiudicato
a Royal Medal al Prof. Hirst
for his researches in pure mathematics.*]*

47

Sanremo
29 Feb$^{y}$ 1884

My dear Friend
I am making my way slowly towards Rome, where I hope to arrive during the course of the month on which we shall enter tomorrow. I left London on the 8th of December, but was detained in Paris for more than a fortnight in consequence of a curios accident which befel[1] *[sic!]*me. I was on my way, just after dusk, from my Hotel, Rue de la Paise, to the Lyons Railway Station, along the Quais, when some clever thieves climbed up behind the cab, cut the straps by which my portmanteau was secured on the top, and carried it off. The road being very rough, and the noise of the cab being great, I was not aware of what had happened until, on my arrival at the Gare, I found myself destitute of luggage
of course I went at once to the Police; but to my great disgust they took no active steps to apprehend the thieves, and I had to return to my Hotel to provide, next day, for my most pressing necessities, in the way of clothing, and to wait patiently until another portmanteau, with a fresh supply, arrived from England. I have heard nothing since of my stolen property, and do not expect to do so. Fortunately I never carry money in my portmanteau, but amongs other things I had some Geometrical M.S.S. in it which I never hope to recover or replace
This accident so discouraged me that I gave up my projected journey to Algiers, and decided to spend the next month or two on the Riviera, - at Marseilles, Nice, Mentone and here. I do not like to return home without shaking your hand; but I shall not go further than Rome, - and shall return probably in april[2] *[sic!]*. I shall probably pass through <u>Genoa</u> in the course of a week or ten days
I should be grateful if you would send me a line there (Poste Restante) and let me know <u>how</u> you are, and if I shall be sure to find you on my arrival in Rome.
You will probably receive a parcel, by Book Poste, for me before I arrive. If so please return it with kind greetings to <u>La Signorina Itala</u>
I remain Ever yours affectionately

T. Archer Hirst

[1] befell.
[2] April.

48

Rome
5 April 1884

My dear Cremona
I dined at my hotel last night, and I am sorry I did so; for I should have much preferred the dinner at the "Circolo", and a walk with you afterwards. I will dine with you tomorrow at 6.30, when I hope to hear the result of the Election at the Lincei.
Many thanks for the "Popolo Romano" containing the account of your election as the representative for your University at Edinburgh. I will give you all the information I can relative to English money tomorrow night. I propose to leave Rome for Florence on Thursday next the 10th.

Ever yours sincerely
T. Archer Hirst

49

*[su carta intestata:*
Athenaeum Club
Pall Mall *]*
19 April/84

My dear Cremona
You can have a bed room at 17 Pall Mall (Chandler's) just opposite this club.
The charge will be five shillings a night. You can take all your meals at this club, of course, since it is open from 8.30 A.M. to Midnight.
This would be the simplest and most economical arrangement for you I have often availed myself of it when I had occasion to pass one or two nights only in Town.
If you would like to find a quiet and reasonable hotel near here you had better consult "Robert", the Head Waiter in our Coffee Room, should I not be at home when you arrive. In any case you had better speak to Robert before you seek lodgings

Your sincerely
T. Archer Hirst

P.S.
I write to you in Edinburgh by this post. I leave this here, in case my letter should not have reached you

50

*[su carta intestata:*
Athenaeum Club
Pall Mall *]*
19 April/84

My dear Cremona
I received your letter of the 14th immediately on my arrival yesterday. I am afraid this may not arrive in Edinburgh in time to reach you before you leave, but I will send it nevertheless.
You can have a bed room at N° Pall Mall (just opposite this Club) at the reasonable rate of five shilling a night (I say reasonable, because in this neighbourhood everything is expensive) and you can take <u>all</u> your meals at this club. This would be the most economical and convenient arrangement for you.
I have occasionally adopted it myself for a few nights. The bed room will not be an elegant one, (it will be on the third floor) But the people are known at this Club and may be trusted
You can have a bed room <u>and</u> sitting room in the same house for £ 4, a week, but probably you will prefer the bed room alone and will find it sufficient. This club, I may add, is open from 8.30 A.M. to Midnight.
With respect to your second question I will speak to you when you arrive.

Ever yours sincerely
T. Archer Hirst

Should you prefer being at some Hotel in this neighbourhood, "Robert", the head waiter in the Coffee Room of this Club, will give you useful information. You had better consult him in every case if I am not here when you arrive.

## 51

7 Oxford α Cambridge Mansions
Marylebone Road
London (N.W.)
22 April/84

*[in testa alla lettera l'autore ha aggiunto queste righe:* (the above is my private address)*]*

My dear Cremona
I have received your post cards of the 19th and 20th and I write a line to you at Glasgow as you requested. When you last wrote to me you had evidently not received the letter I wrote to you on the 19th about the possibility of finding lodgings in London on your return. If you arrive before 6 P.M. on Thursday you will probably find me at the club; at 6.30 I have to dine at the Royal Society Club, and I should be delighted if you could accompany me as my guest *[l'autore ha inserito una nota:* no evening dress required*]* Should you arrive later, however, "Robert", the Head Waiter in the Coffee Room of the Athenaeum Club, will give you what assistance he can in finding lodgings or a quieter Hotel than the Charing Cross one
I shall return from the R.S. Club dinner to the Athenaeum Club as about 8.30 P.M. and remain there for an hour, at least. Should you arrive still later and I be unable to see you on Thursday I shall call at the Athenaeum between 12 and 1 on Friday and ascertain if you have arrived.
I had news of you from Sylvester, on Sunday after his return from Edinburgh. He has now returned to Oxford.
I may possibly find a letter or post card from you when I call at this club to day and if necessary I will write again by to day's post; but should I not do so I hope very shortly to welcome you in person on your arrival in London

Ever, dear friend,
Yours sincerely
T. Archer Hirst

P.S.
Debus showed me the graceful, and well written letter he received from you on his return from Germany

52

Athenaeum Club
24 April 1884

My dear Cremona
I have just received your postcard, written yesterday from the head of Loch Lomond[1]. You do right to make the best use of your visit to this country, and you will be very welcome whenever you reach London
I will call at 17 Pall Mall and tell them that I expect you about the 29th
Give my kindest regards to Dr Salmon and to Prof. Townsend, if you see him, and believe me to be ever yours sincerely

T. Archer HIrst

P.S. I dine with Tyndall tomorrow. He expected you would accompany me; but I have just sent word to him that your return to London is postponed for a few days.

53

7 Oxford and Cambridge Mansions
Marylebone Road (N.W.)
28 April 1884

My dear Cremona
I am glad to learn, from the letter I have just received from you, that you will arrive in London on Wednesday night. I will take care to give the requisite notice of your arrival at 17 Pall Mall, as well as at the Club. I am sorry to say you will not find me in as good health as I enjoyed *[…]* you left me in Rome. The unusually cold weather (with East Wind) we have been suffering from here has had its bad effect upon me. This morning it is dark, *[…]* rather gloomy, and wet. I have been longing to be back again in your own sunny land. Although the sun ought to be in our meridian at this moment.
I am writing with a lamp on my table.
I will endeavour to see you at the Ath. Club on Thursday on the evening of that day I hope you will dine with us at the *X*-Club[2]. Since Spottiswoode died the only members are Hooker, Tyndall, Huxley, Spencer Busk, Frankland, *[…]* Lubbock, and myself. All old friends, about whom I have often spoken to you. We dine together on the first Thursday of every month
Wishing you a favourable passage across from Ireland
I remain

Yours sincerely
T. Archer HIrst

[1] Lago della Scozia meridionale.

[2] L'X-Club era un circolo informale, un cenacolo di nove uomini di scienza legati da amicizia e da interessi comuni. Si riuniva ogni primo giovedì del mese (eccetto nei mesi di luglio, agosto e settembre) a cena. La prima seduta ebbe luogo a novembre del 1864, l'ultima a marzo del 1893. Ne facevano parte eminenti ed eclettici scienziati: George Busk (1807-1886) chirurgo, zoologo e paleontologo, Edward Frankland (1825-1899) chimico, Thomas Archer HIrst (1830-1892) matematico, Joseph Hooker (1817-1911) botanico, Thomas Huxley (1825-1895) biologo, John Lubbock (1834-1913) archeologo e biologo, Herbert Spencer (1820-1903) filosofo, biologo e antropologo, William Spottiswoode (1825-1883) matematico e fisico e John Tyndall (1820-1893) fisico. Si veda R.M. MacLeod, "The X-Club a Social Network of Science in Late-Victorian England", *Notes and Records of the Royal Society of Londo,* 24 (2), 1970, pp. 305-322.

54

*[su carta intestata:*
Athenaeum Club
Pall Mall S.W.*]*
28 April 1884

My dear Cremona
I have just had a telegram from D$^r$ Salmon suggesting that you should remain at the Euston Hotel on Wednesday night, instead of at 17 Pall Mall. His apprehension was that you would not be able to obtain any food on your arrival.
This, however, is not the case; you can obtain food at this club at any time up to midnight, as that there is really no necessity for you to remain at the Euston Hotel for one night.
I should recommend you to call first at the Athenaeum Club, and order your supper, and then go across to your lodgings. I have informed the proprietor of the latter *[l'autore ha inserito una nota:* M$^r$ Chandler*]* that you will arrive at about 10.30 P.M. on Wednesday.
Your sincerely
T. Archer Hirst

55

7 Ox. α Cam. Mansions
Marylebone Road (N.W.)
9 May/84

My dear Cremona
Here is another Copy of the Catalogue of Articles sold at the Stores.
Au Revoir Lundi!
I hope I shall be in a less miserable condition, to bid good bye to my friend and brother
Ever yours
T Archer Hirst

56

7 Ox. α Cam. Mansions
Marylebone Road (N.W.)
9 May/84

My dear Cremona
In my hurried note of to day I forgot to mention one or two things: viz.
Your invitation to the Athenaeum Club has already been prolonged. Pray consider yourself as an Honorary Member as long as you may be able to remain in London, - the longer the better for us all!
The proofs of your interesting communication to the Mathematical Society shall be forwarded to you in Rome as you desire
I enclose a printed statement of the conditions under which you can be provided with separate copies. You have merely to write on the proofs which will be sent to you, how may separate copies you desire to have
Should this number exceed 25 (to which you are entitled) an account for the extra charge will be forwarded to you by the Printer, when he sends you the copies in question. If you will return this account to me I will at once settle.
The few affectionate words with which you closed your letter of to day touched me as deeply as those of a brother could. I cannot tell you what a bitter disappointment I have experienced at being so ill during your visit to London.
I had long been looking forward to your visit, and had formed many plans for your entertainment and my gratification
To be compelled to renounce them all is a real affliction from which it will take me long to recover.
I still hope to spend Monday evening with you at the Club, but I do so with misgivings.
In any case I will send you a line on Monday morning

Ever, dear friend,
Yours affectionately
T Archer Hirst

57

7 Ox. & Camb. Mansion
Marylebone Road (N.W.)
11 May/84[1]

My dear Cremona
The four days imprisonment to which I have subjected myself have had for effect a decided diminution of the inflammation and swelling of my legs. I propose to drive down to the Athenaeum Club. Tomorrow (Monday) afternoon. I shall be there at about 3 P.M. and I will accompany you to the stores where we shall find all the things you wish to purchase.
After purchasing them we will return to the Club and at 7.30 dine together.
This being the last opportunity I shall have, for some time to come, of dining with you, you must if you please be my guest
I have also invited Debus and Spencer to dine with me to meet you. Tyndall I regret to say is not accessible. He is at his country house
Hoping you have enjoyed your visit to Oxford and that nothing will interfere with my plans for tomorrow *[l'autore ha inserito una nota:* should anything unforeseen happens I will of course send you word at once.*]*
I remain

Ever yours sincerely
T Archer Hirst

[1] Hirst scrive la data "11 maggio 1874", ma, considerando le informazioni a noi note, indicate di seguito, si deve dedurre che ha sbagliato anno e che si tratta, in realtà, dell'11 maggio 1884. Infatti:

I. Cremona si recò in Inghilterra la prima volta nel settembre 1876 e la seconda dall'8 aprile al 13 maggio 1884;
II. Hirst scrive che il giorno seguente sarà lunedì: l'11 maggio 1884 è una domenica mentre l'11 maggio 1874 è un lunedì;
III. nella lettera del 27/07/1884 Cremona fa riferimento allo stato di salute di Hirst durante il suo recente soggiorno in Inghilterra (si veda L. Nurzia (a cura di), "La corrispondenza tra Cremona e Hirst", in: *Per l'archivio della corrispondenza dei Matematici Italiani – La corrispondenza di Luigi Cremona (1830-1903)*, Vol. IV, Quaderni P.RI.ST.EM.);
IV. nella lettera a Cremona del 17/05/1874 (n. 15) Hirst lamenta di non avere sue notizie da ben nove mesi, da quando, cioè, si sono visti a Rapallo.

58

*[su carta intestata:*
Athenaeum Club
Pall Mall*]*
2 July 1884

My dear Friend
It is nearly two months since you left me, and although, during that time, I have only received one brief post.card from you I know, through M$^{r}$ Tucker, that you have arrived safely in Rome. I am looking forward, consequently, with pleasure to the arrival of a letter from you, which will give me fuller information concerning your journey home, your occupations since you arrived there and your projects for the summer.
Has Brioschi entered on his duties as President of the Lincei ?
Has anything further transpired regarding his election and *[…]* preparations for removal to the Palazzo Corsini proceeding?
I shall be especially interested to hear how Itala and Vittorio are, and whether they were pleased with the numerous purchases you made for them when you were here
I am glad to be able to tell you that I now enjoy good health. All the ills from which I suffered so much when you were here, and which interfered so distressingly with my power of enjoying your visit and of administering to your own enjoyment of it; all these ills I am glad to say have now disappeared, and my health and strength have not only returned, but have much improved relatively to what they were when I was in Italy. I told you, I think, that I had lost weight during my residence abroad, although I was far below my proper weight on leaving England in December last. I am glad to tell you however that during the last two months I have not only recovered the 14 pounds I had lost in Italy, but have even reached 156 pounds (English); a weight beyond any to which I have approached for the last three years. This news can only interest you as confirming any statement that my health has really improved greatly and the sign of it is that I have now resumed my geometrical work that and my duties on the Senate of the University of London, and on the Council of University College completely occupy me at present 60
Nothing of much interest, for you, has occurred here since you left. Sylvester is in town, his lectures at Oxford being concluded for this session. I have not heard whether they were deemed to be successful
Cayley has communicated a paper to the Royal Society on <u>non-euclidean plane geometry</u> which will interest Beltrami and yourself more than it does me. Tyndall rarely comes to town the weather is too hot for him here. Lord Claud Hamilton, M$^{rs}$ Tyndall’s Father, died very suddenly about a month ago. Lady Claud, Mrs Tyndall’s Mother, was at Heidelberg at the time undergoing a critical operation for cataract which has not been quite successful
She is still there, and it is feared that total blindness is inevitable, this is very sad, and has cast a gloom over the Tyndalls which their departure for Switzerland alone can brighten temporarily. Let me hear from you soon. Ever dear friend yours sincerely

T Archer Hirst

59

Athenaeum Club
Pall Mall S.W.
August 8th 1884

My dear Friend

Your long-expected letter reached me about a week ago and from its perusal I derived great pleasure and consolation – the consolation of finding my brotherly regard for you so affectionately returned. I am glad you are enjoying, with Itala, the fresh alpine breezes of your native land. You did right to avoid crossing over the frontier; for the absurd disinfesting and quarantine arrangements one hears so much about are not to be undergone except in cases of most extreme necessity. I trust the cholera panic will have subsided before next october[1] *[sic!]*, otherwise I fear I must relinquish my project of passing the winter in Sicily. But there are already signs of improvement and I have not yet given up hope.

This cholera scare has changed the plans of many English Travellers; - to the great loss I fear, of many continental pleasure providers, and the great gain of English Hotels $\alpha$ Irish ones. I gave your message to Sylvester, who is still here, his mind possessed with his theory of Matrices to the conclusion of all thought on his forthcoming Inaugural Lecture at Oxford. Tyndall has gone alone to his chalet at the Bel Alp, leaving his wife to take charge of her sick mother. I received the Milanese paper containing an account of the trial in which the honour of Brioschi and, as you truly say, the credit of the Lincei is implicated

Surely many members now regret that the choice of a President did not fall on you! I read with pleasure that you are contemplating the resumption of your Geometrical work. I hope to live to see the intention realized

When you next write you must explain to me more in detail what "notizie" relative to the university of London, and to university college you desire to possess. I can easily send you the most recent volumes of the Calendars of these two institutions as well as of King's College London. These I think will contain all the information you desire to possess

Give my loving remembrances to Elena, Itala and Vittorio

Ever yours sincerely

T. Archer HIrst

[1] October.

60

7 Oxford & Cambridge Mansions
Marylebone Road (N.W.)
23 Oct. 1884

*[in testa alla lettera l'autore ha aggiunto queste righe:* Give my kindest remembrances to Itala and Vittorio. I trust his examination will pass off satisfactorily*]*

My dear Cremona
I have bad news to send you respecting myself.
From the time you left me last May up to the 9th of this month my health continued to improve, and I was making preparations for my departure to Sicily on the 1st of November.
Now alas my prospects are completely changed, and I scarcely hope to be able to leave England at all this winter. Quite unexpectedly I found myself losing the control of my left arm and leg. I recognized at once that it was a return of that attack of partial paralysis which caused me, eight years *[ago]* to give up all work at Greenwich and go to Egypt. Then, however, the attack was caused by over work and mental anxiety, now, no such causes can be assigned. My work has been merely a pleasant occupation for me and I have had no anxiety whatever. It is purely constitutional, therefore
Sir Andrew Clark, my Physician, does not consider it a case of ordinary paralysis, proceeding from the brain. My intellect in fact, has remained up to this time perfectly clear, and in its normal condition. It would appear, however, that the low motor nerves, proceeding from the spinal column, are somehow pressed upon and rendered useless by congested capillary blood vessels. Whatever may be the cause, it is a deplorable fact that for the last week I have been unable to leave my room and even to do the simplest and most necessary things for myself. Thanks to the perfect quietude in which I have been kept, I am decidedly better to day than I was a week ago, and I have still some hope of regaining a good deal of the power of motion which I possessed ten days ago. My hopes are founded solely upon the fact that I was so well before the attack came on. For the present, however, all intention of going to Sicily must be postponed. Should I continue to improve I will write again and tell you what I propose to do.
I received the newspaper containing the news of the disgraceful collapse of the Company with which Brioschi was connected.
I agree with you that it is simply deplorable to see a man who occupies so high a position in Science soiling not only his own reputation but science itself by such unworthy connexions.
Surely the eyes of the men who raised him to the position of President of the Lincei are by this time opened!
By this time you will I trust have received copies of the Calendar 1) of the University of London 2) of University College and 3) of Kings College.
Tyndall has returned in good health from Switzerland, and Mrs Tyndall, who remained at home to comfort her mother, is perfectly well.
Huxley, the President of The Royal Society, is out of health and has gone abroad for six months rest.
Debus has returned in perfect health. When I remember that he taught me Chemistry and is still lecturing on the subject as brilliantly as ever I feel ashamed of myself.
Spencer health too has improved. He has been writing actively lately. Mrs Spottiswoode called upon me a day or two ago. She is well and wholly devoted to her two sons. I have <u>not</u> heard from Leudersdorf[1].
Ever yours affectionate

T. Archer Hirst

[1] Leudesdorf.

61

*[su carta intestata:*
Athenaeum Club
Pall Mall S.W. *]*
13 Nov. 1884

My dear Friend
Your letter reached me this morning and was very welcome. From the fact that I write to day in reply to it, and that this, my letter, is written at the Athenaeum Club, you will naturally and justly infer that my health is improved.
I remained perfectly quiet at home for a month, and gradually the use of my left leg was restored to me. My left hand, however, is still useless to me; though the pain in it has somewhat diminished. My physician, Sir Andrew Clark, has satisfied himself that it is <u>not</u> ordinary paralysis, proceeding from the brain, from which I suffer; but he frankly confesses he does not know <u>what</u> it is. The object of his treatment of me is simply to keep up my general strength, and leave nature to do the rest. He still encourages me to pass the winter in a warmer climate, and he think I might travel, especially by Sea, without danger.
But I, for my own part, hesitate to travel alone; I am so helpless that my courage fails me, and moreover I am so bad a sailor. If I could be transported safely to Rome, I should be thankful for I should have the great advantage and consolation of your help and company
But I do not yet see how this can be managed
Sicily according to your account has now become impossible for me.
It is possible I may go to Gibraltar or Malta by one of this fine steamers that ply between London and India, via the Suez canal; but as yet I am unable to decide. On any case you will hear from me again before I start.
I was glad to hear from you that your daughter Elena, (whose name I cannot now recall) with her children, is coming to visit you. It will be a great joy to Itala and yourself, as well as to her.
Pray give her my kindest regards and best wishes
The news that Vittorio has passed his examination and become a Civil Engineer is also satisfactory. Give him also my best wishes; for his future welfare and success I have every hope. All my friends here, male and female, have been exceedingly kind to me during my illness. Even good M^rs^ Spottiswoode came twice to see me, and to try to cheer me. I need scarcely say that M^rs^ Tyndall was as affectionate as a Sister to me.
Ever dear friend

Yours affectionately
T Archer Hirst

62

"Splendid" Hotel
Nice
20 Feb$^{y}$ 1885

My dear Cremona
It is time that I sent you a few lines to tell you of my wandering since I last wrote to you; which, I believe, was in November, when I was recovering from that illness which so unexpectedly prevented me from leaving home at the beginning of November and going to Italy to spend the winter. Gradually the new attack subsided, and by the 14$^{th}$ of December I was able to leave England once more. I was afraid, however, to undertake a long journey, so I gave up my intention of going to Sicily, and contented myself with the Riviera again. I spent seven weeks at Cannes where, notwithstanding the cold, and often wet, weather, I continued to gain strength and courage. On Saturday last I came here, where I shall remain, if all go well with me, for a month at least. After that I may go on to some other place on the Riviera, and gradually, approach Genoa. Thence I shall probably visit the Italian lakes and, as the warm weather returns, approach England once more by the S$^{t}$ Gothard, Lucern and Brussels. I don't think it probable that I shall turn southwards to Pisa, Florence, or Rome. My plans, however, are not yet fixed I cannot say with certainty what I may do.
The desire to occupy myself a little with Geometry is taking possession of me again and I may possibly find it desirable to place myself within reach of a good mathematical library.
Up to the present I have had no books near me to consult, and the time may come when pen and paper will no longer suffice. Yesterday I had the great pleasure of meeting Zeuthen at Monaco (Monte Carlo) and his conversation on Geometrical matters was very refreshing to me. He gave me some news of you, too; which was very welcome to me. He left you in January , he said, both well and busy; - busy, he told me, with the report on University matters which you are preparing for the Senate. This, no doubt, will be good and important; but I could not help reflecting, not so interesting to your Scientific friend as Geometry. I do not forget, however, that even in Geometry you have produced more than many of us during the year. I have no news whatever, concerning our scientific friends, to communicate to you. 64
In England they appear to look upon me, already, as an emigrant who is half lost to them. The last I heard of Sylvester was that he was writing English verses and challenging Oxford undergraduates to translate them into Latin or Greek! Of Cayley I have heard absolutely nothing since I left. Tyndall, Debus, $\alpha$ Spencer are well, and Huxley, when I last heard of him, was leaving Naples for Florence. He may possibly have passed through Rome by this time
I trust your son and daughters are all well. Are they still with you? If so give them all my very kindest remembrances
Let me soon hear from you and believe me my very dear friend to be yours affectionately

T. Archer Hirst

63 [1]

Splendid Hotel
Nice
(19;3;85)

My dear Friend
At the end of this month I propose to leave this place and to travel by Genoa and Milan towards the S^t Gothard. I should like to call upon Tardy and Jung. Can you send me their present addresses? My health continues fairly good. I have just finished a paper on the Congruence (3,3) treated by Kummer in 1878 in the Monats Bericht[2]. How are you and all around you? Give my kindest greetings to all. Did you receive the letter I sent to you on February 20^th?
Ever yours sincerely

T. Archer Hirst

64

*[su carta intestata:*
Splendid Hotel
Nice
G. Laurent*]*
30. March 1885

My dear Cremona
I received your letter last night. On Friday last I sent you a Post Card to say that I propose to leave here on the 4^th or 5^th of April and to spend a short time (perhaps a week) at Bordighera (Poste Restante)Thence I shall go direct to Genoa, where I shall probably remain one night only, unless there should be a possibility of meeting you there. At any rate I shall call at the Post office in Genoa for my letters
From Genoa I shall go on to Milan and thence, I think, to Stresa; but of this I am not yet quite certain. I shall, however, travel slowly towards the St Gothard and thence, by Brussels, to England.
I shall be glad to be kept informed of your movements for I have not yet given up the hope that we may meet on the Riviera.
I remembered at once that it was Itala's birth day on the 28^th When you see or write to her give her my hearty congratulations. Are you going to bring her to England this summer?
I have only seen *[…]*, vague reports of the disturbances, in Turin Padua and Rome, caused by students. I hope your enquiry will soon terminate and the disaffection likewise.

Ever yours affectionately
T. Archer Hirst

1 Cartolina postale indirizzata a "M le Professeur Cremona – R. Scuola degl'Ingeneri *[sic!]* – S. Pietro in Vincoli – Roma l'Italie".
2 T.A. Hirst, "On the Cremonian Congruences which are contained in a Linear Complex", *Proc. London Math. Soc.,* s1-17 (1), 1885, pp. 287-296.

## 65 [1]

G$^{d}$ Hotel de Bordighera
10 April 1885

I have just received your Post card of the $7^{th}$. I shall remain here probably until Wednesday next ($15^{th}$) and then go to Genoa. It would be a great delight to me to be able to see you, and the Signorina Itala, either here or, in any case, on passing Porto Maurizio. On receiving a post card from you I will write again giving a more precise statement of my movements

Ever yours
T. Archer Hirst

## 66

G$^{d}$ Hotel de Bordighera
13 April 1885

My dear Cremona
I received your telegram on Saturday, in reply to mine of the same day; and yesterday came your Post Card, giving me full and sufficient explanations of your inability to come to Bordighera. I need not tell you what a bitter disappointment it was to me to learn that, after being so close to you, I should have to leave Italy without seeing you and Itala
I felt the disappointment the more keenly because I was conscious of not having fully carried out the programme I had formed and communicated to you when your letter from Turin informed me that you would return to Rome on the $2^{nd}$ I thought it highly improbable that you would be able to come to Porto Maurizio as soon as you actually did *[…]* so. Moreover the weather at Nice became so bad that, fearing it would be still worse in North Italy, I decided to postpone my departure for a few days. Accordingly it was only on the $10^{th}$ that I received your letter written on the $7^{th}$, at Rome I immediately sent a Post Card in reply to Genoa (Posta Restante) telling you I was here and expressing the hope that I might either receive a visit from you and Itala, or at any rate meet you on my passage through Porto Maurizio on my way to Genoa. Next morning, that is to say on the $11^{th}$, I received your first post card from Porto Maurizio and answered it, as you know, by a telegram. I am sorry now I did not take the train and go at once to Porto Maurizio to see you. Had I known that you were compelled to return to Rome so speedily I certainly should have done so; for it would have been a great satisfaction to me if I could have grasped your hand once more, and have had even an hour's friendly talk with you.
Regrets are useless, however, and it is better to look to the future to repair the mistakes and disappointments of the past. I hope you intend to carry out the plans you formed last year, and to come with Itala to England. In a month's time I expect to be once more at home. I will send you a Post Card on my arrival and, in reply, I shall look for a letter from you, giving me full details concerning your projects for the Summer.
Zeuthen intended to come to Bordighera to meet you. I have just written to him at Menton (Pension S$^{t}$ Maria) to say it is useless. He too will share the disappointment of your constant friend.

T. Archer Hirst

[1] Cartolina postale indirizzata a "Al Professore Cremona (Senatore del Regno) – Poste Restante – Genova".

# 67

*[su carta intestata:*
Athenaeum Club
Pall Mall S.W.*]*
13 July 1885

My dear Friend
It is now three months (if not more) since I wrote to you. I have often thought of you and often been on the point of writing when something or other interfered with my intention. At length a friendly remark of yours, forwarded to me by M$^{r}$ Tucker, has decided me to delay no longer.
After leaving Bordighera I travelled through Porto Maurizio and halted for a week at Alassio, where I had beautiful weather. Thence I went to Genoa, where the weather changed for the worse. I went on to Milan, where it was no better. After waiting there for a few days, in the company of our friend Jung, I went on to Lugano, hoping that the weather would recover itself and permit of my enjoying the Italian lakes. This, however, it would not do; the rain fell incessantly until at last, weary of waiting for fine weather, I crossed the S$^{t}$ Gothard and halted again for a day or two at Lucern
People assured me that the wind, which brought rain to Lombardy, left sunshine in Switzerland, but this was not the case and I left Lucern to go, without further delay, to Brussels
There it was cold as well as wet, and after a days[1] *[sic!]* halt I left Brussels and came home, where I have remained perfectly quiet ever since.
I have had plenty of occupation; for after five months absence from home, and after neglecting so many social and other duties, many arrears had accumulated. I have moreover been quietly working at my Geometry, and the two months that have elapsed since my return have passed away with great rapidity.
My health has not been seriously interrupted though I have rarely been free from neuralgic pain. Mannheim tells me that you are not going to Switzerland as he and I expected, and the fact that you have given up this project, as well as the prior one of coming to England, causes me to fear that private affairs 67
may have influenced you unfavourably. I hope my surmise is baseless; but I should be glad to be assured by yourself that it is so
My best wishes will go with you and *[la frase continua in testa alla lettera:* Itala whenever you may spend your vacation. I hope you will find time soon to send me a line to assure me that you are both well and happy*]*

Evers yours
T Archer Hirst

*[l'autore ha inserito una nota:* Tyndall has been very poorly; he has just left for his place in Switzerland. Cayley, a few days ago, was as well and calm as ever Sylvester is at present in Paris.*]*

[1] day's halt.

68

14 August 1885

My dear Cremona
I was very glad to receive your letter of 20th July in which you assured me that my fears were quite groundless relative to the reasons which had induced you to change your plans for passing your vacation. It was especially satisfactory to me to find that one of your reasons for remaining on the Riviera was that you had purchased a house at Porto Maurizio; for it cannot fail to be of great advantage both to you and to Itala to be able to retire there occasionally to enjoy quietude and good air. I was sorry to hear, however, that you were somewhat anxious about Itala's health. The last time I saw her at Rome she appeared to me to be so strong, so full of animation and healthy vigour, that I cannot now imagine that anything serious threatens her physical welfare
At her age palpitations, although distressing enough so long as they persist, are often merely indicative of transient disturbances which pass away with change of air and mode of life. I remember that in my own early days I suffered from disturbances of this kind, thorough I have reason to believe that there never was anything organically defective about my heart
Such I trust will prove to be the case with Itala, and that a long, happy and healthful life is before her. Pray remember me to her affectionately, and when you next write to Elena ( I always forget her present name) tell her that I hear with joy of her domestic happiness. I am beginning to feel solitary in London now but I am thankful to say that my health remains fairly good, and that I can occupy myself daily with geometrical thought and reading. Tyndall is at Bel Alp, near Brieg, whiter[1] his wife has now also gone. Debus is at Wilhelmshöhe, near Cassel[2], his native place. Both enjoy good health, and are evidently amusing themselves
Spencer is in the lake district of Westmoreland[3]. Cayley is in Wales and Sylvester, after being in Paris and Switzerland, has now returned to London
In October next (my health and the Cholera permitting) I hope *[la frase continua in testa alla lettera:* to go again to Italy if possible I will spend the winter in Sicily 68
In any case I hope to be more fortunate in meeting you than I was last april[4] *[sic!]*

Ever my dear friend
yours affectionately
T. Archer HIrst*]*

*[nota inserita dall'autore:* I am glad also to hear that Vittorio has begun his practical career as an Engineer convey to him my best wishes for his future health and prosperity.*]*

[1] Termine arcaico che significa "dove".
[2] Così si chiamava fino al 1926 la città di Kassel, in Germania.
[3] Contea storica del Regno Unito, dal 1974 fa parte della contea amministrativa di Cumbria.
[4] April.

69

Hotel Quirinal
24 Nov. 1885

My dear Cremona.
At 9 P.M. yesterday I concluded that I should not have a visit from you, and I went out for a walk until bed time. On my return I was much disappointed to find that you had called shortly after I had left the Hotel.
On returning from my walk I had the misfortune to be caught in a violent storm of rain. After taking shelter for some time under a portico, I hailed an omnibus. I ran after it in the wet; but it would not stop. In despair, I called a cab; but on stepping into the wretched, low vehicle, my hat rolled off into the middle of the street. The driver, safe under his big umbrella, would not stir from his dry seat to recover my hat for me, and I was compelled to do so myself. The rain poured down in torrents on my bare head, and my hat, (itself full of liquid) instead of protecting me gave me a veritable shower bath.
I returned, as you may readily imagine, in a wretched condition, - wet both internally and externally. I had to go to bed at once. This morning, as a natural consequence of my adventure, I have a cold. I do not intend to remain at home all day, however; weather permitting, I will try and restore myself by exercise. But, this evening, I fear I must not go out. Tomorrow at 6.30 P.M. I trust I shall be sufficiently recovered to accept your hospitality
Kindest remembrances to Itala and Vittorio

Ever yours sincerely
T Archer Hirst

70

Hotel du Quirinal
29th Novr 1885

My dear Cremona
on my return to this Hotel to day I found the letter of which you spoke. It was the letter you wrote to me on the 27th of April last, and which arrived at Genoa, probably, after I had left that City.
One item only of its contents was news to me, and very good news. It was that in which you expressed your intention of bringing Itala to England this year. I hope nothing will occur to prevent your doing so and above all I hope I shall be there to receive and welcome you

Until tomorrow at 6.30 P.M.
Ever yours sincerely
T Archer Hirst

71

Hotel du Quirinal
1 : 12 : 85

My dear Cremona
M^r^ Potter, the English member of Parliament whose acquaintance I have made here, would like, (with his Daughter), to visit that Exhibition of Drawings of S. Marco which we saw on Sunday last, as well as the Statuary that pleased us so much at the Baths of Diocletian. If you could lend me a ticket which would admit the two, I should feel greatly obliged to you, and they would be very pleased.
I gave them to day my copy of Rossi's description of the Temple of Vesta. They are reading it to night and will visit the Forum tomorrow. Both they and I went to day to the Farnesina Palace and enjoyed it very much

Ever yours sincerely
T Archer Hirst

72

Hotel Quirinal
15 : 12 : 85

My dear Cremona
Herewith I send you a list of the principal memoirs of Froude with which I have been furnished by a very competent Authority Professor Cotterill of Greenwich.

Guccia desires to be kindly remembered to you. He adds:
"Veuillez rappeler le Circolo Matematico di Palermo au bon, souvenir de M. Cremona!" – a message which you will probably understand better than I do.
Hoping to see you tomorrow evening at 9 o'clock
I remain

Yours sincerely
T Archer Hirst

73

17 : 12 : 85

My dear Cremona
Many thanks for Schur's memoir. I have just written to Cotterill about the Transactions of the Institute of N. Architects and asked him to communicate the result of his enquiries to you directly
I will gladly dine with you on Saturday
Ever yours

T. Archer Hirst

74

Grand Hotel
Naples
[2 : 1 : 86]

My dear Cremona
I send you a line to convey to you my very best wishes for the happiness of Itala, Vittorio and yourself during the year which has now commenced. Long before it expires I hope to able[1] *[sic!]* to welcome you and Itala, and perhaps Vittorio also, to England. I received, a day or two ago, the Report of the Speeches you made to senate during my sojourn in Rome. I shall retain the Pamphlet as a souvenir of that pleasant visit.
Thus far I have derived both pleasure and profit from my visit to Naples; though the weather, I confess, has been far from favourable
Dini and Battaglini have each paid me a visit, and the latter has kindly lent me the memoir of Schur which I wished to consult. He informs me that Guccia has already reached Rome, and that I may expect to see him here in the Course of a week's time. Sir Henry Thompson has changed his Hotel from the Bristol to this one. I see him daily; but he leaves on Wednesday for Monte Carlo. I trust your visit to Porto Maurizio was a pleasant one. Ever yours sincerely

T. Archer Hirst

[1] to be able.

75

*[su carta intestata cancellata:*
Circolo Bellini
Palermo *]*
Hotel des Palmes
Palerme
26th Jan$^y$ 1886

My dear Cremona
Thanks to our good friend Guccia I am now most comfortably settled here. I have access, as you see, to a very elegant club, where I find not only all the conveniences I could desire, so far as Journals and newspapers are concerned, but also very pleasant and hospitable acquaintances.
I arrived on Saturday morning last, from Naples, where I spent three weeks very profitably, so far as my health is concerned; the weather however, was sometimes very cold, and often very wet.
I saw Battaglini, Caporali and one or two other Geometers; but neither Gasparis[1] nor Fergola. I left Naples on the evening of a beautiful day, Friday last, and had a smooth and most favourable passage, by Steamer, to this place. On the very day of my arrival, however, the weather changed again and with few intervals it has rained ever since; so that I have not been able to make a single excursion outside this town.
I am hoping, from day to day, that a more favourable weather will compensate me for all my disappointments
I have not yet presented the letter you very kindly procured for me to the Principe di Scalea; but the moment I see an opportunity of profitably visiting Sicilian antiquities I shall avail myself of your introduction.
I have received a letter from Sturm, wherein he tells me that he has enquired of you concerning my travels. I have already answered it; so that, on my account alone, you need not trouble yourself to do so.
Guccia has already put the rooms, the journals, and the memoirs of the Circolo Matematico di Palermo at my disposal
He really deserves great credit for the disinterested zeal he has manifested in establishing this useful little institution When once warmer weather arrives I hope to do some work there
He desires to be very kindly remembered to you
I was amused the other day to read, in the Italian papers, that Modena and Sienna *[sic!]* (I think) were about to apply to be raised, like Genoa, Messina and Catania, to towns of first rank relative to their Universities
This results might have been foreseen, and I am curious to see what consistent reply M. Depretis can give them. I hope you continue well, under your increased duties, and that you will sometimes remember
Your affectionate friend
T. Archer Hirst

My warmest remembrances to Itala and Vittorio

[1] Probabilmente si riferisce ad Annibale De Gasparis.

76

Grand Hotel
Napoli
{15;4;86}

My dear Friend
I have just received your letter of yesterday, and I hasten to reply that I shall certainly not leave Naples until Sunday next indeed, since probably I may have the pleasure of seeing you here, I am not sure now that I shall halt at Rome at all although I had written to the Quirinal to secure a room for Sunday next, I will not decide on leaving here until I have seen you. Should any change in your plans take place, therefore, please telegraph to me at once, as you propose.
I have had bad weather ever since I came to Naples on the 8th; so much so that I have not been able to pay my intended visit to Sorrento and Capri. It is now uncertain whether I shall be able to do so. On one point only I am quite decided, - I will not run the risk of failing to see you before I return home. My health has not been as good as it was before I arrived here; - a fact that I attribute, chiefly, to the very great change which has taken place in the weather. After being bright, dry, and warm, in Taormina, Cava, and Amalfi, it has become wet and stormy and sometimes cold.
Since my return to Naples, a week ago, I have done nothing, and seen nobody. I have been watching the weather, and living from day to day in the expectation, or rather in the hope of a favourable change therein.
Should I stop at Rome a day I shall not fail to pay a visit to S. Pietro in Vincoli in order to bid good bye to Vittorio and Itala.
Trusting to seeing you here on Saturday, however, I will write no more at present, except to wish you a pleasant journey southwards and success to your mission, whatever that may be

Yours affectionately
T. Archer Hirst

77

16 July 1886

My dear Friend
Since I arrived home on the 8th of May – now more than two months ago! – I have often thought of you and of yours, and I have often entertained an intention of writing to you. But you know how such vague intentions become dissipated; especially when, in consequence of half a year's absence from home, so many arrears have accumulated. The receipt yesterday, however, of a large black-edged card announcing the death of Caporali, gave me so great a shock, that I at once decided that my long silence must be broken.
When I last saw Caporali, in the Cour d'Honneur of the Grand Hotel at Naples, you were present. He had, I remember, a severe cold, and I thought his manner less genial and joyous than usual But apart from this he looked strong, and in robust health; I should not have hesitated to ascribe to him a long and useful life.
What was it that caused his life to terminate so prematurely? In vain I read, many times and carefully, the black-edged card; but it gave me no information.
You will do so, I trust, when you next write. I had a great regard for Caporali; not merely as a geometer, but as a man. He had a great esteem for you; and it was a pleasure to me, whenever I conversed with him, to notice this.
Just before leaving Rome, in fact the evening before, I made the acquaintance of Kronecker. He would no doubt tell you of this and I will merely add that our interview interested and gratified me much.
At Paris I saw Halphen twice, - once at the Academy, and once at his own dinner table. Mannheim had a sick house. His children were attacked by scarlatina, and Madame was, in consequence, invisible. I saw Mannheim himself, however, several times, as well as de Jonquiéres, Bertrand and others.
Tyndall, when I reached home, was recovering from his severe illness; since then, however, he has had a relapse from which he is now once more recovering. He will probably soon leave for the Bel Alp. I have not yet seen Cayley. I hear he has lost his youngest Sister and is much distressed thereby.
I have often seen Sylvester. He is well and lively. At the present moment he is at St Andrews in Scotland (he has returned to day! unexpectedly)
Immediately after my return to London I made a communication to the Mathematical Society "on the Cremonian Congruences which are contained in a linear complex". I have just finished writing it. You will recognize in it a development of the subject on which you were speaking when in Rome last December.
I trust Itala and Vittorio are both well; pray remember me to both kindly. I send the letter to Rome thought I expect you are no longer there. Wherever you may be, however, - at your house at Porto Maurizio or on the more breezy alps – my good wishes for your health and happiness will be with you

Ever yours sincerely
T. Archer Hirst

78

*[su carta intestata*
Athenaeum Club
Pall Mall*]*
13 May 1887

My dear Cremona
After an unusually long silence I find myself able, once more, to send you my friendly greeting from my customary place in the Athenaeum Club. I returned a week ago from my long and eventful visit to Spain.
I have been absent from home for nearly six months, during which time I have travelled over more than 4000 miles through a country which was entirely new to me. Although I found much to interest me in Spain, I cannot say that my journey through it was a successful one or that I shall ever wish to revisit it. The climate, with a few exceptions, was a trying one for me, and my health suffered much in consequence of the many trials I had to undergo. These trials had their origin not merely in the climate; but also in the food and entertainment which with a few exceptions, were either bad or quite unsuitable for me.
One consequence of all this is that I have returned home very thin, and very feeble.
In a short time, however, I trust these defects will be repaired. One of the most marked deficiencies in Spain is the absence of all scientific life and interest. I did not make the acquaintance of a single Spaniard who occupied himself with Mathematics or Physics or Astronomy or indeed with Science in any form! This dearth of scientific interest imparted itself to me even for I did nothing in Geometry during the whole six months!
It was a perfect holiday for me; or rather a perfect change of thought and occupation. The Architectural and artistic features of the Country, its historical associations, and the peculiarities and customs of its inhabitants occupied my whole attention and, I ought to add, repaid it fully. I entered Spain at Barcelona; followed the Eastern Coast down through Tarragona and Valencia to Alicante, and then crossed the Continent to Cordoba, Seville and Cadis. Thence I went by sea to Gibraltar; then across to Tangier, and finally after other cruises, I reentered Spain at Malaga. From the last mentioned place I went to Granada, and finally to Madrid, Toledo, and the Escorial. At the last dismal and chilling place I caught a serious cold (the third from which I had to suffer in Spain), and I had in consequence to terminate my wanderings abruptly. I took refuge at Biarritz, in France, - a place I had visited 30 years before – and remained there until I recovered, I then returned home, via Paris where I saw the Mannheims and from them received news of you.
I was sorry to find that our friend had once more tried, unsuccessfully, to enter the Academy of Sciences, and moreover that he had lost, by death, his excellent mother, whom I knew well and esteemed highly.
From Madame Mannheim I heard of Vittorio's departure for America, after some disappointment connected with his love affair. The incident, I trust, will have no serious consequences, and that he will return, before long, to resume his professional work in Italy, near you and Itala.
Of your own proceedings during the session which is now drawing to a close, I have heard nothing whatever; but I hope to have a few lines from yourself before long. Here there are not many changes
I saw Sylvester yesterday for a few moments and he appeared to be in great spirits. Cayley, Glaisher tells me, is very well but very busy with university work. Tyndall, as you may have heard, has resigned his Professorship at the Royal Institution, and is now, like Huxley, Hooker Frankland and myself, unfettered by official duties.
Lord Rayleigh has become Tyndall's successor at the Royal Institution, as well as Stokes' successor, as Secretary of the Royal Society. But I must close my letter with kindest regards to Itala
Ever yours sincerely
T. Archer Hirst

*[l'autore ha inserito una nota:* I trust your house and property at Porto Maurizio did not suffer from the recent tremblements de terre.*]*

79

*[su carta intestata*
Athenaeum Club
Pall Mall S.W.*]*
27 July 1887

My dear Cremona
I do not propose, in this letter, to reply to the very welcome one I received from you a few weeks ago. I write principally to ask if there is any truth in the report which has found currency here, (in consequence of a paragraph which appeared in Nature) that you are going to be present at the meeting of the British Association which is to be held at Manchester during the first week of the month of September
When I first saw the paragraph to which I have alluded, I at once assumed that it must have been published by mistake; but I have just been speaking to Sir Henry Roscoe, the President Elect, and he assures me that you are really expected
In your letter to me you expressed the opinion that you would be unable, this year, to go further than to the Italian Alps
I shall indeed be delighted to find that you have altered your plans, and are really coming to pay us a visit. Should you come, I shall certainly alter <u>my</u> plans also, and go to Manchester to meet you. Pray let me know, therefore, at your earliest convenience, if I may indeed expect to see you in London, on your way to Manchester, and when you will arrive.
I will postpone writing more until your reply reaches me, and I will conclude my letter by sending my kindest remembrances to Itala. Will she not come with you? I hope she may; but from what you said in your last letter to me I fear she will not.
I hope you have good news of your son Victor *[…]*
Ever yours sincerely
T Archer Hirst

80

Hotel Mont – Fleuri
Cannes
Alpes Maritimes
France

3 Jan$^{y}$_1888

My dear Cremona
It is a long time since I wrote to ask you if the report was true that you were about to honour the British Association by your presence at their Manchester Meeting! The reply which I promptly received from you fully realized my anticipations. The report, I felt sure, was based on a misunderstanding. Since then I have been occupied incessantly, though it is with difficulty that I can now recall anything of importance that I have done.
I left London however on the 7$^{th}$ of November, spent a few days in Paris, where I had the pleasure of seeing our old friend Mannheim and his family, and of exchanging a friendly greeting with Guccia, just before he returned to Palermo.
On passing through Rome he called on you, probably, and gave you news of me.
At that time I was very undecided with respect to my own intentions for the winter.
I thought, at first, of going to Algiers; but on reaching Marseilles the weather was so stormy that I yielded to my fears (I am such a bad sailor!) and came here instead. I now propose to remain here sometime longer and ultimately to find my way home again, in the spring.
I scarcely hope, now, to pay a visit to Rome; much as I should be rejoiced to grasp your hand again, and to greet, affectionately, Itala. But the journey is so long, and every year travelling becomes more and more toilsome and even wearysome to me. Moreover I am not without some hope of being able to see you, on the Riviera, before I return home. I will keep myself in communication with you, at all events, in the event of your having occasion to pay a visit to your house at Porto Maurizio.
I trust your daughter Elena, and Itala are both well, and that you have good news of Vittorio The Emperor of Brazil, as you know, is here. He occupies the adjoining Hotel (Beau – Sejour), in fact, and I see him frequently.
I do not know whether, under your existing Government you are occupied more or less than usual with political matters; but I trust your occupations, whatever their nature may be, do not deprive you <u>wholly</u> of the opportunity to pursue our favourite study. Since I came here, I have been trying to complete work commenced more than ten years ago!

Ever sincerely yours
T. Archer Hirst

81

4 June 1888
(posted June 5)

My dear Cremona
Itala's pretty little volume, "Le Alpi", reached me the day before yesterday, and gave me great pleasure, especially since it contained, also, your own card, and thus proved to me that I was not forgotten by my real friend; although he certainly did appear to have lost all memory of me, for a long time past. Pray thank Itala for me and tell her I congratulate her heartily, on her authorship.
It was in January that I received your last letter. It reached me in Cannes. In March I sent you a post card telling you that I had just returned from Corsica, in very bad health, and that the season on the Riviera had been a most disastrous one for me. I wrote at the same time, and to the same effect, to our mutual friend Mannheim in Paris. In his reply to my letter, he asked me if I had heard of your marriage; which of course I had to answer in the negative. Since then indeed the same intelligence has reached me from other sources; but I have waited in vain for any confirmation of the news from yourself.
I can, however, no longer entertain any doubt about the truth of the report, and I have no longer any ground for postponing my congratulations, and sending you my very best wishes.
I am glad to be able to supplement my present brief communication by the intelligence that, at length (I am assured by my medical advisers) the crisis of my own malady has been passed; and that I may hope, soon, to recovery my lost health and strength.
This recovery, however, has not as yet manifested itself visibly. I am still very weak and very thin. Nevertheless *[continua in testa alla lettera:* I trust I shall have a better account to send you soon, and in the mean time I repeat my congratulations and best wishes.
Ever yours sincerely
T. Archer Hirst*]*

82

26 June 1888

*[In una nota a matita, probabilmente scritta da Cremona:* scr. 10 Sett. 1888 da Silvaplana*]*

My dear Cremona
Although I could not account for your silence towards me, when I was first informed of your marriage, I felt sure that it could not have arisen from any omission on your part to inform your old english[1] *[sic!]* friend of so important an event in your life.
Your last letter, enclosing as it did the one which you sent to me at Cannes, simply verified this conviction of mine. I repeat my congratulations, therefore, and heartily wish you and M^rs^ Cremona every happiness. I remember very well meeting her at your house, when I was last in Rome, and I trust she has not forgotten me.
I hope your visit to Bologna was a pleasant one. I ought to have been there myself officially; the Senate of the University of London having requested me, in conjunction with D^r^ Pole, to represent them on the occasion.
Their letter was sent to me in Ajaccio; but arrived there after my sudden departure, - a departure which was caused as you know by my illness, - and thus never reached me. Had it reached me in time, however, I should have been unable to accept the honourable office assigned to me. The weather, according to the newspaper account, was extremely hot in Bologna. Your King, it would appear suffered from the heat; but you, I hope, bore it with impunity
I am glad to hear that you and M^rs^ Cremona propose to spend part of your vacation in Switzerland after visiting your property in the Riviera. May you both return invigorated and refreshed by your holiday!
I hope you have still good news of Vittorio, and of Elena, as well as of Itala
Remember me affectionately to all them, whenever you are able to do so.
My health improves but it does so very slowly

Ever yours sincerely
T. Archer Hirst

[1] English.

83

*[su carta intestata*
Athenaeum Club
Pall Mall*]*
26th Sep 1888

My dear Cremona
The letter which I received from you when you were at Silva Plana gave me great pleasure.
I know the place well and, in former days, I paid many visit to it, to St Morits, and to Pontresina.
During the time you were there my friend Huxley was at the Maloya. I wish you had known it; for it would have been a pleasure to both if you could have met and renewed acquaintance with one another. He arrived there, this year, in very bad health, but he improved as far as to be able to walk five or ten miles without fatigue – a feat which was far beyond his strength when he left England. I was glad to learn that your own visit to the Engadine had been equally successful, and I trust you have now returned to Rome with health quite re-established. My recovery has been a slow one; but without any relapse, and I shall attempt to pass the coming winter in England, - probably at Bournemouth in the South of our little island. Tyndall and his wife are at the Bel Alp. They are both in good health and have had, recently, excellent weather. They will not return until the snows of October drive them down. Spencer's health is still very poor. He has spent the summer in the country; but not far from London. He will probably pass the winter, also, at Bournemouth. Sylvester is here. I see him almost daily. He spends his time, chiefly, in writing sonnets, several of which he has published in "Nature". Goethe it is said, was prouder of his botanical researches, and of his theory of light, than of his poetry. Sylvester I believe esteems his poetry to be of greater merit, even, than his mathematics! It is almost a year since I saw Cayley I hear, however, that he has spent the summer with his son in Switzerland
I am very glad to hear that you have good news of Vittorio, as well as of Elena and of Itala. Pray give my kindest remembrances to all three, when you have any opportunity and please, give my best regards also to the Signora Cremona
*[continua in testa alla lettera:* Ever yours, with unalterable friendship, T. Archer Hirst*]*

84

*[su carta intestata*
Athenaeum Club
Pall Mall S.W.*]*
12 Jan$^{y}$ 1889

My dear Cremona.
The letter which I received from you on the 10$^{th}$ gave me great pleasure. I heartily reciprocate your good wishes, for the new year! May it bring health and happiness to you, to M$^{rs}$ C. and to Elena, Vittorio and Itala! In all form please convey my kindest regards.
I was sorry to hear of your indisposition at Silva Plana, and to learn that it had returned after your arrival in Rome. May I express to you a word of warning relative to the cold douch you are taking every morning. The almost universal experience of men of our age is that cold bath, except in very warm climates, are prejudicial to health; often indeed accompanied with serious risks. At all events their effects require to be closely watched. Watch them my dear friend.
I am glad to say that my health, although far being good, has suffered no serious relapse during the last nine months. We have not had a severe winter, as yet, and I have been able, consequently, to remain in London, within reach of medical advice, should such again prove to be needful. I shall probably go to Bournemouth in early spring.
I occupy myself occasionally with geometrical studies; but I have produced little, during the last year, worthy of publication
Cayley came to see me about a month ago. He was in good health; greatly interested in new studies of his only son (who is now a student of Trinity College Camb.) and fully occupied with the republication of his voluminous researches. It is expected that they will fill ten quarts volumes
Sylvester, for some unknown reason, is not friendly with me just now; so that I can give you no news of him.
Spencer is still very ill. We never see him. *[l'autore ha inserito una nota:* He lives at Dorking*]* – His best friends, he says, are now his worst enemies. They excite him nervously. Poor Fellow!
Huxley's health has greatly improved since his return from the Engadine; and Tyndall and M$^{rs}$ T. are both in excellent condition. They live, however, at a distance from London, and we rarely see them. (Haslemere)
I have procured, from the best Billiard Table Makers in London, copies of the rules which apply to the three games that are usually played in England. I send them to day by book post.
And now, my dear friend, accept my heartiest good wishes. I trust your lectures on the "Theory of Imaginaries in Geometry" will be brilliantly successful.
Ever yours affectionately

T. Archer Hirst

85

7 Oxford & Cambridge Mansions
Marylebone Road, London, (N.W.)
27 Nov$^{r}$ 1890

*[In una nota a matita, probabilmente scritta da Cremona:* scr. 17-7-91*]*

My dear Cremona
Your letter, commenced on July 26 and continued on Nov. $^{r}$ 20, pleased me very much.
The part of it was especially gratifying which inform me of the approaching marriage of Itala.
The fact that Elena had already arrived to take part in your family rejoicings, caused me to conclude that the marriage was about to take place at once. Accordingly, I sent you a telegram without delay, expressing in three words ("Every happiness Itala") my hearty good wishes.
I hope you received it.
Until your letter of July 26 arrived I had heard nothing of the telegram, to Professor Stokes, which you sent on the evening of the Royal Society Conversazione.
Notwithstanding this miscarriage, I thank you now, sincerely, for your kindly greeting.
You ask if there is no hope of my coming, to pass the winter in Rome. Alas no!
I fear I shall never see Rome again.
The malady from which I suffer, and shall suffer, I fear, for the rest of my days, is a slow form of paralysis. It had already begun to manifest itself when I was in Rome in 1877. I recovered a little afterwards; but it has now nearly disabled my left-leg. I walk with difficulty; but I still continue to walk a little daily; - from my home to the Athenaeum Club. My head remains clear; I have never had any, so called, "paralytic stroke". My medical advisers describe the case as one of "creeping paralysis"
These however are names merely. The simple fact remains that I grow daily more infirm, and more unable to travel; - my former delight.
Tyndall, too, is far from well. But his malady is very different from mine. He still retains his wonderful activity. He has only recently returned from his mountain perch on the "Bel Alp". But, poor fellow, he cannot sleep; insomnia pursues him every where[1] *[sic!]*. How it will end I know not.
Cayley is fairly well as is also Spencer. Debus has left London. He now lives with his invalid Sister in Cassel, his native place. Huxley has also left London. He lives at Eastbourne and is quite well. Stokes is about to retire from the Presidentship of the Royal Society. Sir William Thomson will become his successor.
I hope you have good news of Vittorio. Remember me to him kindly when you write. Give my best regards also to Mrs Cremona and to Helena. Remember that a letter from you is always welcome

Ever yours affectionately
T. Archer Hirst

[1] everywhere.

86

*[su carta intestata*
Athenaeum Club
Pall Mall S.W.*]*
9 January 1892

My dear Cremona
The receipt of your last letter gave me very great pleasure, - as indeed your occasional letters always do, and always did.
Let me first answer your kind enquiry concerning my health: I become more and more infirm.
But I do so very gradually, and, at the present time, I am not worse than I expected to be (on Jan^r^ 1^st^ 1892)
I become more and more unable to travel. In August last, however, I succeeded in spite of many difficulties, in reaching Paris. I saw none of our friends there. All were absent. My sole object in going was to settle, finally, certain arrangements connected with the grave, at Montmartre, of my dear, lost wife. This done, I returned home, without delay, and have remained here ever since.
The same difficulty, of locomotion, will, I fear, prevent me from going to Dublin; although I have two sisters-in-law, who live there, and would welcome me. I hope, however, to be here when you and Mr^s^ Cremona pass through London. I need not say that I shall be delighted to see you both.
Your news of Elena (and of her son 10 years odl[1]!), of Itala, and of Vittorio was very welcome to me. Remember me most kindly to all three when you next write to them.
I mourned, greatly, the death of Casorati. I did not know that Kronecker had also been called away; but I received, two days ago, the melancholy intelligence of the death of my friend, and fellow student Heinrich Schröter.
I still remember, well, the welcome you gave to him, and Mrs Schröter, when they were last in Rome with me.
I hear with sorrow, of the illness of Betti and Battaglini; but I am glad to know that Beltrami and Brioschi are both still well, and ever young! Remember me kindly to both.
Tyndall has been very dangerously ill during the year; but I am glad to tell you he is now convalescent.
Cayley too is ailing; but he works hard at the publication of his collected papers. W. Thomson has just received a Peerage[2]!
I do not yet know by what title he will be known henceforth[3]
Hoping to shake your hand once more, before long, I remain as ever,
Yours sincere friend

T. Archer Hirst

*[l'autore ha aggiunto una nota:* I am glad to hear that you still continue to lecture. Long may you be able to do so!*]*

[1] Ettore Perozzi. Si veda la lettera del 1/1/1892 n. 126 di Cremona a Hirst in: L. Nurzia (a cura di), "La corrispondenza tra Cremona e Hirst", *Per l'archivio della corrispondenza dei Matematici Italiani – La corrispondenza di Luigi Cremona (1830-1903)*, Vol. IV, Quaderni P.RI.ST.EM.
[2] Titolo di "Pari del Regno".
[3] Lord Kelvin.

## Tabella - i dati delle lettere

Nella tabella sono indicati i dati relativi alle lettere del carteggio dell'Archivio Mazziniano di Genova oggetto di questa pubblicazione e, in corsivo, quelli delle lettere di Cremona, dall'Archivio della London Mathematical Society, in risposta o alle quali rispondono le missive di Hirst qui trascritte.
Le lettere di Cremona (tranne la n. 38 e 39) sono tratte dalla corrispondenza pubblicata in: L. Nurzia (a cura di), "La corrispondenza tra Cremona e Hirst", *Per l'archivio della corrispondenza dei Matematici Italiani – La corrispondenza di Luigi Cremona (1830-1903)*, Vol. IV, Quaderni P.RI.ST.EM.
Delle lettere qui trascritte si riporta, nella prima colonna, il numero d'ordine riportato in testa a ogni lettera, delle altre si riporta il numero d'ordine assegnato dalla curatrice nella pubblicazione già citata (Quaderno P.RI.ST.EM.).
Le righe colorate evidenziano quei gruppi di lettere consecutive nelle quali gli autori rispondono uno alla missiva dell'altro senza interruzioni nel dialogo. Questa informazione si ricava molto facilmente dal testo delle lettere dato che entrambi i matematici riportano quasi sempre la data della missiva alla quale sono in procinto di rispondere, anche quando le risposte sono molto tardive.

| N° | Quaderno P.RI.ST.EM | Mittente | LUOGO | DATA | SEGNATURA FILE E CARTACEA | LINGUA |
|---|---|---|---|---|---|---|
| **1** | | Hirst | - | 24/08/1865 | 054-13095 (9778) | Inglese |
| **2** | | Hirst | St. John's Wood (UK) | 24/10/1865 | 054-13094 (9777) | Inglese |
| **3** | | Hirst | St. John's Wood (UK) | 05/04/1866 | 054-13096 (9779) | Inglese |
| | *40* | *Cremona* | *Bologna* | *21/05/1866* | | *Italiano* |
| | *44* | *Cremona* | *Milano* | | | |
| **4** | | Hirst | Bristol (UK) | 26/12/1866 | 054-13097 (9780) | Inglese |
| | *45* | *Cremona* | *Milano* | *12/01/1867* | | *Italiano* |
| **5** | | Hirst | Londra (UK) | 17/02/1867 | 054-13098 (9781) | Inglese |
| | *48* | *Cremona* | *Milano* | *16/06/1867* | | *Italiano* |
| **6** | | Hirst | Parigi (F) | 03/07/1868 | 054-13099 (9782) | Inglese |
| **7** | | Hirst | Londra (UK) | 10/04/1869 | 054-13100 (9783) | Inglese |
| | *60* | *Cremona* | *Milano* | *07/05/1869* | | *Italiano* |
| **8** | | Hirst | Exeter (UK) | 26/08/1869 | 054-13101 (9784) | Inglese |
| **9** | | Hirst | Londra (UK) | 04/10/1870 | 054-13102 (9785) | Inglese |
| | *64* | *Cremona* | *Milano* | *08/01/1871* | | *Italiano* |
| **10** | | Hirst | Londra (UK) | 15/12/1871 | 054-13103 (9786) | Inglese |
| | *67* | *Cremona* | *Milano* | *21/12/1871* | | *Italiano* |
| **11** | | Hirst | Londra (UK) | 26/12/1871 | 054-13104 (9787) | Inglese |
| | *68* | *Cremona* | *Milano* | *06/01/1872* | | *Italiano* |
| **12** | | Hirst | Londra (UK) | 31/01/1873 | 054-13105 (9788) | Inglese |
| **13** | | Hirst | Bel Alp (CH) | 11/08/1873 | 054-13106 (9789) | Inglese |
| | *70* | *Cremona* | *Rapallo* | *16/08/1873* | | *Italiano* |
| **14** | | Hirst | Arona | 23/08/1873 | 054-13107 (9790) | Inglese |
| **15** | | Hirst | Greenwich (UK) | 17/05/1874 | 054-13109 (9792) | Inglese |
| | *71* | *Cremona* | *Roma* | *21/07/1874* | | *Italiano* |
| **16** | | Hirst | Wiesbaden (D) | 02/08/1874 | 054-13110 (9793) | Inglese |
| | *72* | *Cremona* | *Roma* | *31/10/1874* | | *Italiano* |
| **17** | | Hirst | Greenwich (UK) | 15/11/1874 | 054-13111 (-) | Inglese |
| | *73* | *Cremona* | *-* | *09/07/1875* | | *Italiano* |
| **18** | | Hirst | Greenwich (UK) | 14/07/1875 | 054-13112 (9795) | Inglese |
| | *74* | *Cremona* | *Roma* | *21/70/1875* | | *Italiano* |
| **19** | | Hirst | Greenwich (UK) | 26/07/1875 | 054-13113 (9796) | Inglese |
| | *80* | *Cremona* | *Roma* | *06/05/1876* | | *Italiano* |
| **20** | | Hirst | Greenwich (UK) | 31/05/1876 | 054-13114 (9797) | Inglese |
| | *81* | *Cremona* | *Roma* | *05/07/1876* | | *Italiano* |

| N° | Quaderno P.RI.ST.EM | Mittente | LUOGO | DATA | SEGNATURA FILE E CARTACEA | LINGUA |
|---|---|---|---|---|---|---|
| | *84* | *Cremona* | *-* | *26/12/1876* | | *Italiano* |
| **21** | | Hirst | Athens (GR) | 31/01/1877 | 054-13115 (9798) | Inglese |
| | *85* | *Cremona* | *-* | *07/02/1877* | | *Italiano* |
| **22** | | Hirst | Roma | 15/04/1877 | 054-13117 (9800) | Inglese |
| | *86* | *Cremona* | *-* | *28/04/1877* | | *Italiano* |
| **23** | | Hirst | Parigi (F) | 02/05/1877 | 054-13118 (9801) | Inglese |
| **24** | | Hirst | Shafton - Barnsley (UK) | 10/05/1877 | 054-13116 (9799) | Inglese |
| | *87* | *Cremona* | *S. Martino di Castrozza* | *27/07/1877* | | *Italiano* |
| **25** | | Hirst | Darmstadt (D) | 13/08/1877 | 054-13119 (9802) | Inglese |
| | *88* | *Cremona* | *S. Martino di Castrozza* | *24/08/1877* | | *Italiano* |
| **26** | | Hirst | Darmstadt (D) | 02/09/1877 | 054-13120 (9803) | Inglese |
| **27** | | Hirst | Greenwich (UK) | 03/10/1877 | AAA-1361 (1361) | Inglese |
| | *89* | *Cremona* | *Roma* | *07/11/1877* | | *Italiano* |
| **28** | | Hirst | Greenwich (UK) | 26/11/1877 | 054-13121 (9804) | Inglese |
| | *90* | *Cremona* | *Napoli* | *13/05/1878* | | *Italiano* |
| **29** | | Hirst | Greenwich (UK) | 17/05/1878 | 054-13122 (9805) | Inglese |
| | *91* | *Cremona* | *Roma* | *25/06/1878* | | *Italiano* |
| **30** | | Hirst | Shafton - Barnsley (UK) | 03/07/1878 | 054-13123 (9806) | Inglese |
| **31** | | Hirst | Greenwich (UK) | 25/03/1879 | 054-13125 (9808) | Inglese |
| | *92* | *Cremona* | *-* | *01/04/1879* | | *Inglese* |
| **32** | | Hirst | Parigi (F) | 06/04/1879 | 054-13176 (9859) | Inglese |
| | *93* | *Cremona* | *Roma* | *10/04/1879* | | *Inglese* |
| **33** | | Hirst | - | 15/06/1879 | 054-13124 (9807) | Inglese |
| | 94 | Cremona | Roma | 24/06/1879 | | Inglese |
| **34** | | Hirst | Shafton - Barnsley (UK) | 11/07/1879 | 054-13126 (9809) | Inglese |
| | 95 | Cremona | Roma | 15/07/1879 | | Italiano |
| **35** | | Hirst | Münster (D) | 15/08/1879 | 054-13178 (9861) | Inglese |
| **36** | | Hirst | Axenstein (CH) | 15/09/1879 | 054-13175 (9858) | Inglese |
| **37** | | Hirst | Greenwich (UK) | 15/10/1879 | 054-13127 (9810) | Inglese |
| | *97* | *Cremona* | *Roma* | *26/11/1879* | | *Inglese* |
| **38** | | Cremona | - | 1880 | 052-11940 (8627) | Inglese |
| **39** | | Cremona | - | 10/1880 | 052-11938 (8625) | Inglese |
| **40** | | Hirst | Londra (UK) | 18/05/1881 | 054-13128 (9811) | Inglese |
| | *104* | *Cremona* | *-* | *17/07/1881* | | *Italiano* |
| **41** | | Hirst | Greenwich (UK) | 24/07/1881 | 054-13129 (9812) | Inglese |
| | *105* | *Cremona* | *Roma* | *31/12/1881* | | *Italiano* |
| **42** | | Hirst | Devonport (UK) | 05/01/1882 | 054-13130 (9813) | Inglese |
| | *106* | *Cremona* | *Roma* | *28/05/1882* | | *Italiano* |
| | *107* | *Cremona* | *Roma* | *17/07/1882* | | *Italiano* |
| **43** | | Hirst | Greenwich (UK) | 22/07/1882 | 054-13131 (9814) | Inglese |
| | *108* | *Cremona* | *Zuz* | *31/07/1882* | | *Italiano* |
| **44** | | Hirst | Greenwich (UK) | 05/08/1882 | 054-13132 (9815) | Inglese |
| | *109* | *Cremona* | *Roma* | *16/06/1883* | | *Italiano* |
| | *110* | *Cremona* | *Roma* | *30/06/1883* | | *Italiano* |
| | *111* | *Cremona* | *Portomaurizio* | *24/07/1883* | | *Italiano* |
| **45** | | Hirst | Londra (UK) | 28/07/1883 | 054-13133 (9816) | Inglese |
| **46** | | Hirst | Londra (UK) | 19/11/1883 | 054-13134 (9817) | Inglese |
| | *112* | *Cremona* | | *24/11/1883* | | *Italiano* |

| N° | Quaderno P.RI.ST.EM | Mittente | LUOGO | DATA | SEGNATURA FILE E CARTACEA | LINGUA |
|---|---|---|---|---|---|---|
| **47** | | Hirst | Sanremo | 29/02/1884 | 054-13135 (9818) | Inglese |
| **48** | | Hirst | Roma | 05/04/1884 | 054-13138 (9821) | Inglese |
| **49** | | Hirst | Londra (UK) | 19/04/1884 | 054-13139 (9822) | Inglese |
| **50** | | Hirst | Londra (UK) | 19/04/1884 | 054-13140 (9823) | Inglese |
| **51** | | Hirst | Londra (UK) | 22/04/1884 | 054-13141 (9824) | Inglese |
| **52** | | Hirst | Londra (UK) | 24/04/1884 | 054-13142 (9825) | Inglese |
| **53** | | Hirst | Londra (UK) | 28/04/1884 | 054-13143 (9826) | Inglese |
| **54** | | Hirst | Londra (UK) | 28/04/1884 | 054-13144 (9827) | Inglese |
| **55** | | Hirst | Londra (UK) | 09/05/1884 | 054-13145 (9828) | Inglese |
| **56** | | Hirst | Londra (UK) | 09/05/1884 | 054-13137 (9820) | Inglese |
| **57** | | Hirst | Londra (UK) | 11/05/1884 | 054-13108 (9791) | Inglese |
| **58** | | Hirst | Londra (UK) | 02/07/1884 | 054-13146 (9829) | Inglese |
| | *113* | *Cremona* | *Caprile* | *27/07/1884* | | *Italiano* |
| **59** | | Hirst | Londra (UK) | 08/08/1884 | 054-13147 (9830) | Inglese |
| | *114* | *Cremona* | *Roma* | *07/10/1884* | | *Italiano* |
| **60** | | Hirst | Londra (UK) | 23/10/1884 | 054-13148 (9831) | Inglese |
| | *115* | *Cremona* | *Roma* | *10/11/1884* | | *Italiano* |
| **61** | | Hirst | Londra (UK) | 13/11/1884 | 054-13149 (9832) | Inglese |
| **62** | | Hirst | Nizza Marittima (F) | 20/02/1885 | 054-13150 (9833) | Inglese |
| **63** | | Hirst | Nizza Marittima (F) | 19/03/1885 | 054-13179 (9862) | Inglese |
| **64** | | Hirst | Nizza Marittima (F) | 30/03/1885 | 054-13151 (9834) | Inglese |
| **65** | | Hirst | Bordighera | 10/04/1885 | 054-13180 (9863) | Inglese |
| **66** | | Hirst | Bordighera | 13/04/1885 | 054-13152 (9835) | Inglese |
| **67** | | Hirst | Londra (UK) | 13/07/1885 | 054-13154 (9837) | Inglese |
| **68** | | Hirst | - | 14/08/1885 | 054-13155 (9838) | Inglese |
| **69** | | Hirst | Roma | 24/11/1885 | 054-13156 (9839) | Inglese |
| **70** | | Hirst | Roma | 29/11/1885 | 054-13157 (9840) | Inglese |
| **71** | | Hirst | Roma | 01/12/1885 | 054-13158 (9841) | Inglese |
| **72** | | Hirst | Roma | 15/12/1885 | 054-13159 (9842) | Inglese |
| **73** | | Hirst | Roma | 17/12/1885 | 054-13160 (9843) | Inglese |
| **74** | | Hirst | Napoli | 02/01/1886 | 054-13161 (9844) | Inglese |
| **75** | | Hirst | Palermo | 26/01/1886 | 054-13162 (9845) | Inglese |
| **76** | | Hirst | Napoli | 15/04/1886 | 054-13164 (9846?) | Inglese |
| **77** | | Hirst | - | 16/07/1886 | 054-13163 (9847) | Inglese |
| | *116* | *Cremona* | *Portomaurizio* | *28/07/1886* | | *Italiano* |
| **78** | | Hirst | Londra (UK) | 13/05/1887 | 054-13165 (9848) | Inglese |
| | *117* | *Cremona* | *Roma* | *05/06/1887* | | *Italiano* |
| **79** | | Hirst | Londra (UK) | 27/07/1887 | 054-13166 (9849) | Inglese |
| | *118* | *Cremona* | *Cogne* | *03/08/1887* | | *Italiano* |
| **80** | | Hirst | Cannes (F) | 03/01/1888 | 054-13167 (9850) | Inglese |
| | *119* | *Cremona* | *Roma* | *31/01/1888* | | *Italiano* |
| **81** | | Hirst | - | 04/06/1888 | 054-13168 (9851) | Inglese |
| | 121 | Cremona | Roma | 07/06/1888 | | Italiano |
| **82** | | Hirst | - | 26/06/1888 | 054-13169 (9852) | Inglese |
| | *122* | *Cremona* | *Silvaplana* | *10/09/1888* | | *Italiano* |
| **83** | | Hirst | Londra (UK) | 26/09/1888 | 054-13170 (9853) | Inglese |
| | *123* | *Cremona* | *Roma* | *07/01/1889* | | *Italiano* |
| **84** | | Hirst | Londra (UK) | 12/01/1889 | 054-13171 (9854) | Inglese |
| | *125* | *Cremona* | *Roma* | *20/11/1890* | | *Italiano* |
| **85** | | Hirst | Londra (UK) | 27/11/1890 | 054-13172 (9855) | Inglese |
| | *126* | *Cremona* | *Roma* | *01/01/1892* | | *Italiano* |
| **86** | | Hirst | Londra (UK) | 09/01/1892 | 054-13174 (9857) | Inglese |

## Indice dei nomi citati nelle lettere

MT – collegamento ipertestuale alla biografia su *Mac Tutor History of Mathematics Archive* della School of Mathematics and Statistics, University of Saint Andrew, Scotland
http://www-history.mcs.st-andrews.ac.uk/index.html

TR - collegamento ipertestuale alla biografia su *TRECCANI*
http://www.treccani.it/

SR - collegamento ipertestuale alla biografia nell'Archivio storico del Senato della Repubblica
http://notes9.senato.it/Web/senregno.NSF/SenatoriTutti?OpenPage

| Nome | N° lettera/e | | Link/breve biografia |
|---|---|---|---|
| Cerruti Valentino | 28 | SR | http://notes9.senato.it/web/senregno.nsf/96ec2bcd072560f1c125785d0059806a/3dbf7b7e77af9c9e4125646f005a113b?OpenDocument |
| Chasles Michel | 2, 3, 4, 6, 10, 16, 17, 18, 29 | MT | http://www-gap.dcs.st-and.ac.uk/~history/Biographies/Chasles.html |
| Chelini Domenico | 3, 4, 31, 32, 33, 34, 37 | TR | http://www.treccani.it/enciclopedia/domenico-chelini/ |
| Chemin Jean Charles Octave | 42 | | (Péronne 1844–1930) Ingegnere, docente all'École nationale des Ponts et Chaussées. Traduttore di testi scientifici. Cavaliere della Légion d'honneur. |
| Christoffel Elwin Bruno | 1, 37 | MT | http://www-gap.dcs.st-and.ac.uk/~history/Biographies/Christoffel.html |
| Claud Lady Elizabeth (Proby) | 58 | | Suocera di John Tyndall (madre di Louisa Hamilton). |
| Claud Lord Hamilton | 58 | | Suocero di John Tyndall (padre di di Louisa Hamilton). |
| Clausius Rudolf Julius Emmanuel | 36, 37 | MT | http://www-gap.dcs.st-and.ac.uk/~history/Biographies/Clausius.html |
| Clebsch Rudolf Friedrich Alfred | 2, 10, 12 | MT | http://www-gap.dcs.st-and.ac.uk/~history/Biographies/Clebsch.html |
| Clifford William Kingdone | 7, 31, 32 | | http://www-gap.dcs.st-and.ac.uk/~history/Biographies/Clifford.html |
| Cotterill James Henry | 72, 73 | | (1836 – 1922). Professore di Matematica applicata al Royal Naval College di Greenwich dal 1873 al 1897. |
| Crelle August Leopold | 31 | MT | http://www-gap.dcs.st-and.ac.uk/~history/Biographies/Crelle.html |
| Cremona Antonio Luigi Gaudenzio Giuseppe | 10 | MT | http://www-gap.dcs.st-and.ac.uk/~history/Biographies/Cremona.html |
| Cremona Cozzolino Itala | 15, 16, 21, 25, 39, 43, 44, 45, 47, 58, 59, 60, 61, 64, 65, 66, 67, 68, 69, 70, 74, 75, 76, 77, 78, 79, 80, 81, 82, 83, 84, 85, 86 | | (Bologna 1865 – Genova 1939). Terza figlia di Luigi Cremona ed Elisa Ferrari. |
| Cremona Perozzi Elena | 15, 16, 21, 25, 29, 30, 39, 40, 41, 43, 59, 61, 68, 80, 82, 83, 84, 85, 86 | | (Pavia 1856 - ?). Primogenita di Luigi Cremona ed Elisa Ferrari. |
| Cremona Tranquillo | 30 | TR | http://www.treccani.it/enciclopedia/tranquillo-cremona_(Dizionario-Biografico)/ |
| Cremona Vittorio | 15, 16, 21, 25, 39, 43, 44, 58, 59, 60, 61, 68, 69, 74, 75, 76, 77, 78, 79, 82, 83, 84, 85, 86 | | (Bologna 1861 - ?). Secondogenito di Luigi Cremona ed Elisa Ferrari. |
| Darboux Jean Gaston | 29 | MT | http://www-gap.dcs.st-and.ac.uk/~history/Biographies/Darboux.html |
| De Gasparis Annibale | 3, 75 | SR | http://notes9.senato.it/Web/senregno.NSF/c1544f301fd4af96c125785d00598476/49abcf55023b1618c1257069003186a4?OpenDocument |

| Nome | N° lettera/e | | Link/breve biografia |
|---|---|---|---|
| De Jonquières Ernest Jean Philippe Fauque | 2, 4, 77 | MT | http://www-gap.dcs.st-and.ac.uk/~history/Biographies/Jonquieres.html |
| De Morgan Augustus | 1, 4 | MT | http://www-gap.dcs.st-and.ac.uk/~history/Biographies/De_Morgan.html |
| Debus Heinrich | 24, 51, 57, 60, 62, 68, 85, | | (Wolfhagen 1824 – Kassel 1915). Chimico tedesco. Dal 1873 Insegnò Chimica presso il Royal Naval College di Greenwich. Fu eletto Membro della Royal Sociaty nel 1861. |
| Depretis Agostino | 75 | TR | http://www.treccani.it/enciclopedia/agostino-depretis/ |
| Dini Ulisse | 3, 23, 28, 74 | MT | http://www-gap.dcs.st-and.ac.uk/~history/Biographies/Dini.html |
| D'Ovidio Enrico | 11, 23 | TR | http://www.treccani.it/enciclopedia/enrico-d-ovidio_(Enciclopedia-Italiana)/ |
| Durège Karl Heinrich | 1 | | (Danziga 1821 – 1893). Docente di Matematica al Politecnico e all'Università di Zurigo e in seguito professore all'Università di Praga. |
| Fergola Emanuele | 3, 75 | TR | http://www.treccani.it/enciclopedia/emanuele-fergola_(Enciclopedia-Italiana)/ |
| Ferrari Cremona Elisa | 4, 7, 9, 10, 11, 13, 15, 16, 18, 20, 21, 24, 25, 26, 27, 28, 29, 37, 38, 39, 40, 41, 43, 46, | | (Genova 05/06/1826 – Roma 16/07/1882). Prima moglie di Luigi Cremona. Il matrimonio ebbe luogo il 03/08/1854. |
| Fouret Georges François | 29 | | (1845–1921). Matematico francese. Fu presidente della Société Mathématique de France nel 1887. |
| Frankland Edward | 53, 78 | TR | http://www.treccani.it/enciclopedia/edward-frankland/ |
| Frobenius Ferdinand Georg | 36 | MT | http://www-gap.dcs.st-and.ac.uk/~history/Biographies/Frobenius.html |
| Froude William | 72 | TR | http://www.treccani.it/enciclopedia/froude-william_(Dizionario-delle-Scienze-Fisiche)/ |
| Geiser Karl Friedrich | 1, 36 | MT | http://www-gap.dcs.st-and.ac.uk/~history/Biographies/Geiser.html |
| Genocchi Angelo | 23 | MT | http://www-gap.dcs.st-and.ac.uk/~history/Biographies/Genocchi.html |
| Glaisher James Whitbread Lee | 37, 42, 78 | MT | http://www-gap.dcs.st-and.ac.uk/~history/Biographies/Glaisher.html |
| Goethe Johann Wolgang von | 83 | TR | http://www.treccani.it/enciclopedia/johann-wolfgang-von-goethe_(Enciclopedia-Italiana)/ |
| Griffith George | 28 | | (1833–1902). Direttore del Dipartimento di Scienze ad Harrow (Londra) e Vicesegretario della British Association for the Advancement of Science dal 1862 al 1878 e dal 1890 al 1902. |
| Guccia Giovanni Battista | 38, 72, 74, 75, 80 | MT | http://www-gap.dcs.st-and.ac.uk/~history/Biographies/Guccia.html |
| Gullo | 20 | | Non identificato. |
| Halphen George Henry | 29, 42, 77 | MT | http://www-gap.dcs.st-and.ac.uk/~history/Biographies/Halphen.html |
| Harley Robert | 42, 43 | | (Liverpool 1828 – Londra 1910). Matematico e pastore della Chiesa congregazionalista di Leicester. |
| Hermite Charles | 10 | MT | http://www-gap.dcs.st-and.ac.uk/~history/Biographies/Hermite.html |
| Hesse Ludwig Otto | 10 | MT | http://www-gap.dcs.st-and.ac.uk/~history/Biographies/Hesse.html |
| Hirst Emily Anna (Lilly) | 20, 21, 23, 24, 25, 26, 30, 45, 46, | | (? - 1883). Nipote di Thomas Archer Hirst. |

| Nome | N° lettera/e | | Link/breve biografia |
|---|---|---|---|
| Hodgson | 17 | | Editore dei *Proceedings of The Royal Mathematical Society*. |
| Hooker Sir Joseph Dalton | 53, 78 | TR | http://www.treccani.it/enciclopedia/sir-joseph-dalton-hooker/ |
| Huxley Thomas Henry | 53, 60, 62, 78, 83, 84, 85 | TR | http://www.treccani.it/enciclopedia/thomas-henry-huxley/ |
| Janni Vincenzo | 3 | | (Barletta 1819 – Napoli 1891). Matematico. Nel 1863 fondò, con Battaglini e Trudi, il *Giornale di Matematiche*. |
| Jung Giuseppe | 21, 24, 63, 67 | TR | http://www.treccani.it/enciclopedia/giuseppe-jung_%28Enciclopedia-Italiana%29/ |
| Key Sir Astley Cooper | 12 | | (1821–1888). Ammiraglio della Royal Navy. |
| Klein Felix Christian | 29 | MT | http://www-gap.dcs.st-and.ac.uk/~history/Biographies/Klein.html |
| Kronecker Leopold | 77, 86 | MT | http://www-gap.dcs.st-and.ac.uk/~history/Biographies/Kronecker.html |
| Kummer Ernst Edward | 63 | MT | http://www-gap.dcs.st-and.ac.uk/~history/Biographies/Kummer.html |
| Leudesdorf Charles | 60 | | (1853–1924). Matematico. Tradusse in inglese gli *Elementi di Geometria Projettiva* di Luigi Cremona (1885). |
| Lubbock John | 53 | TR | http://www.treccani.it/enciclopedia/avebury-sir-john-lubbock-barone/ |
| Manher Muller Cremona Anna | 82, 83, 84, 85, 86 | | Seconda moglie di Luigi Cremona, lo sposò nel 1887. |
| Mannheim Victor Mayer Amédée | 29, 67, 77, 78, 80, 81 | MT | http://www-gap.dcs.st-and.ac.uk/~history/Biographies/Mannheim.html |
| Mannheim (madame) | 78 | | Moglie di Mannheim. |
| Martorelli Giacomo | 20 | | (Roma 1849 - ?). Ingegnere. Fu Ispettore capo del Genio Navale e deputato del Regno d'Italia. |
| Maxwell James Clerk | 11 | MT | http://www-gap.dcs.st-and.ac.uk/~history/Biographies/Maxwell.html |
| Menabrea Luigi Federico | 28 | MT | http://www-gap.dcs.st-and.ac.uk/~history/Biographies/Menabrea.html |
| Merrifield Charles Watkins | 42 | MT | http://www-gap.dcs.st-and.ac.uk/~history/Biographies/Merrifield.html |
| Miller William John Clarke | 1 | | (1832–1903). Matematico autodidatta. Dal 1850 al 1897 fu direttore della rivista *Educational Times*. |
| Neumann Carl Gottfried | 3 | MT | http://www-gap.dcs.st-and.ac.uk/~history/Biographies/Neumann_Carl.html |
| Nöther/Noether Max | 37 | MT | http://www-gap.dcs.st-and.ac.uk/~history/Biographies/Noether_Max.html |
| Padula Fortunato | 3 | MT | http://www.treccani.it/enciclopedia/fortunato-padula_(Dizionario-Biografico)/ |
| Perozzi Ettore | 86 | | (1882 –1952). Figlio primogenito di Elena Cremona Perozzi. |
| Pincherle Salvatore | 43 | MT | http://www-gap.dcs.st-and.ac.uk/~history/Biographies/Pincherle.html |
| Pole William | 82 | | (1814–1900). Ingegnere inglese. Ebbe la cattedra di Ingegneria civile all'University College di Londra dal 1859. |
| Potter Thomas Bayley | 71 | | (1817 –1898). Politico inglese del Liberal Party, fu membro del Parlamento dal 1865 al 1895. |
| Price Bartholomew | 42 | | (1818–1898). Matematico inglese. Fu docente di Matematica al Pembroke College di Oxford dal 1845. |
| Principe di Scalea/Lanza Spinelli Francesco | 75 | SR | http://notes9.senato.it/web/senregno.nsf/7d795bf0b249d716c125711400599ff4/ba0a24bc1d0f66874125646f005cac70?OpenDocument |

| Nome | N° lettera/e | | Link/breve biografia |
|---|---|---|---|
| Prohuet Eugéne | 3 | | (1817–1867). Matematico francese. Docente (Répétiteur) presso l'École polytechnique. Fu, tra l'altro, l'editore dei *Nouvelles annales de mathématiques.* |
| Prym Friedrich | 1 | | (1841–1915). Matematico Tedesco. Insegnò all'università di Strasburgo e di Würzburg. |
| Rayleigh Lord/ Strutt John William | 78 | MT | http://www-gap.dcs.st-and.ac.uk/~history/Biographies/Rayleigh.html |
| Reye Karl Theodor | 16, 37 | MT | http://www-gap.dcs.st-and.ac.uk/~history/Biographies/Reye.html |
| Riemann Georg Friedrich Bernhard | 1, 3 | MT | http://www-gap.dcs.st-and.ac.uk/~history/Biographies/Riemann.html |
| Roscoe (sir) Henry Enfield | 79 | | (1833–1915). Chimico inglese. Fu Presidente della British Association for the Advancement of Science nel 1887. |
| Salmon George | 3, 4, 52, 54 | MT | http://www-gap.dcs.st-and.ac.uk/~history/Biographies/Salmon.html |
| Sannia Achille | 11 | SR | http://notes9.senato.it/web/senregno.nsf/3b06b7313c966b4dc125711400599aa3/cb3dd5403678a7054125646f00604f7d?OpenDocument |
| Schröeter Heinrich Eduard | 16, 22, 37, 86 | MT | http://www-gap.dcs.st-and.ac.uk/~history/Biographies/Schroeter.html |
| Schröeter Mrs. | 86 | | Clara Rodewald, moglie di Heinrich Schröeter. |
| Schubert Hermann Cäsar Hannibal | 35 | MT | http://www-gap.dcs.st-and.ac.uk/~history/Biographies/Schubert.html |
| Schur Friedrich Heinrich | 73, 74 | MT | http://www.treccani.it/enciclopedia/friedrich-heinrich-schur_(Enciclopedia-Italiana)/ |
| Scott Charlotte Angas | 42 | MT | http://www-gap.dcs.st-and.ac.uk/~history/Biographies/Scott.html |
| Sella Quintino | 40 | TR | http://www.treccani.it/enciclopedia/quintino-sella/ |
| Sidler George Joseph | 1 | MT | http://www-gap.dcs.st-and.ac.uk/~history/Biographies/Sidler.html |
| Siemens Sir Carl Wilhelm | 20 | | (1823–1883). Ingegnere, inventore e industriale tedesco naturalizzato britannico. Fu membro della Royal Society of London for Improving Natural Knowledge e della British Association for the Advancement of Science di cui venne eletto presidente nel 1882. |
| Smith Henry | 8, 37, 42 | MT | http://www-gap.dcs.st-and.ac.uk/~history/Biographies/Smith.html |
| Spencer Herbert | 53, 57, 60, 62, 68, 83, 85 | TR | http://www.treccani.it/enciclopedia/herbert-spencer/ |
| Spottiswoode William | 8, 25, 28, 30, 31, 37, 42, 45, 53, | MT | http://www-gap.dcs.st-and.ac.uk/~history/Biographies/Spottiswoode.html |
| Spottiswoode Mrs. | 28, 45, 46, 60, 61 | | Eliza Taylor Arbuthnot, moglie di William Spottiswoode |
| Steiner Jakob | 3, 4 | MT | http://www-gap.dcs.st-and.ac.uk/~history/Biographies/Steiner.html |
| Stephanos Kyparissos | 42 | | (1857 –1917). Matematico greco. Fu membro della Société mathématique de France. |
| Stokes Sir George Gabriel | 78, 85 | MT | http://www.treccani.it/enciclopedia/sir-george-gabriel-stokes/ |
| Sturm Friedrich Otto Rudolf | 18, 23, 25, 26, 35, 42, 75 | MT | http://www-gap.dcs.st-and.ac.uk/~history/Biographies/Sturm_Rudolf.html |
| Sylvester James Joseph | 2, 3, 4, 51, 58, 59, 62, 67, 68, 77, 78, 83, 84, | MT | http://www-gap.dcs.st-and.ac.uk/~history/Biographies/Sylvester.html |
| Tardy Placido | 3, 63 | TR | http://www.treccani.it/enciclopedia/placido-tardy/ |
| Thompson Sir Henry | 74 | | (1820 –1904). Chirurgo ed erudito britannico. |

| Nome | N° lettera/e | | Link/breve biografia |
|---|---|---|---|
| Thomson William (Lord Kelvin) | 85, 86 | MT | http://www-gap.dcs.st-and.ac.uk/~history/Biographies/Thomson.html |
| Tortolini Barnaba | 1 | TR | http://www.treccani.it/enciclopedia/barnaba-tortolini/ |
| Townsend Richard | 4, 52 | | (1821–1884). Matematico Irlandese. Fu membro della Royal Society of London for Improving Natural Knowledge. |
| Transon Abel | 3 | | (1805 –1876). Matematico francese. |
| Trudi Nicola | 3 | TR | http://www.treccani.it/enciclopedia/nicola-trudi/ |
| Tucker Robert | 58, 67 | MT | http://www-gap.dcs.st-and.ac.uk/~history/Biographies/Tucker_Robert.html |
| Tyndall John | 6, 24, 28, 52, 53, 57, 58, 59, 60, 62, 67, 68, 77, 78, 83, 84, 85, 86 | TR | http://www.treccani.it/enciclopedia/john-tyndall/ |
| Tyndall Mrs. | 60, 61, 83, 84 | | Louisa Hamilton, moglie di John Tyndall. |
| Zeuthen Hieronymous Georg | 16, 42, 62, 66 | MT | http://www-gap.dcs.st-and.ac.uk/~history/Biographies/Zeuthen.html |
| Weyr Eduard | 11 | MT | http://www-gap.dcs.st-and.ac.uk/~history/Biographies/Weyr_Eduard.html |
| Wilson James Maurice | 7 | | (1836–1931). Prete della Chiesa anglicana, teologo, astronomo e docente di Matematica alla scuola di Rugby, cittadino britannico. |
| Wolff | 37 | | Non identificato. |
| Young John Radford | 2 | | (1799–1885). Matematico inglese, fu docente di Matematica presso il Belfast College. |

## Bibliografia